\documentclass[12pt]{article}

\usepackage{microtype}
\usepackage{graphicx}
\usepackage{subcaption}
\usepackage{booktabs} 
\usepackage{hyperref}

\usepackage{amsmath}
\usepackage{amssymb}
\usepackage{mathtools}
\usepackage{amsthm}

\usepackage[margin=1in]{geometry}

\usepackage[capitalize,noabbrev]{cleveref}

\theoremstyle{plain}
\newtheorem{theorem}{Theorem}[section]
\newtheorem{proposition}[theorem]{Proposition}
\newtheorem{lemma}[theorem]{Lemma}
\newtheorem{corollary}[theorem]{Corollary}
\theoremstyle{definition}
\newtheorem{definition}[theorem]{Definition}
\newtheorem{assumption}[theorem]{Assumption}
\theoremstyle{remark}
\newtheorem{remark}[theorem]{Remark}

\def\rmd{\mathrm{d}}
\newcommand{\E}{\mathbb{E}}
\usepackage{multirow}
\usepackage{wrapfig}
\usepackage{placeins}
\usepackage{makecell}
\usepackage{amsmath}

\usepackage{booktabs}
\usepackage{comment}
\usepackage{amssymb}
\def\rmd{\mathrm{d}}
\newcommand{\N}{\mathbb{N}}
\newcommand{\R}{\mathbb{R}}
\usepackage{bbm}

\usepackage{xr-hyper}
\usepackage[textsize=tiny]{todonotes}

\title{Zeroth-Order Langevin Dynamics for Stochastic Nonconvex Optimization: Curvature Regularization and Non-Asymptotic Guarantees}

\usepackage{authblk}

\author[1,2]{Stefano Bruno}
\author[3]{Dong-Young Lim$^{\dagger}$}

\affil[1]{\textit{School of Mathematics, University of Edinburgh, Edinburgh, EH9 3FD, United Kingdom}}
\affil[2]{\textit{UNIST InnoCORE AI-Space Solar Initiative, Ulsan National Institute of Science and Technology (UNIST), Ulsan, 44919, Republic of Korea}}
\affil[3]{\textit{ Department of Industrial Engineering and Artificial Intelligence Graduate School, Ulsan National Institute of Science and Technology (UNIST), Ulsan, 44919, Republic of Korea}}
\date{}
\begin{document}

\maketitle

\begingroup
\renewcommand\thefootnote{}\footnotetext{$^{\dagger}$
Corresponding author: dlim@unist.ac.kr}
\endgroup

\begin{abstract}
Zeroth-order methods provide a way to optimize stochastic nonconvex objectives using only noisy function evaluations. However, the smoothing inherent in their gradient estimators effectively modifies the objective being optimized. In particular, Gaussian smoothing with smoothing parameter $\eta>0$ induces, to leading order, an implicit Hessian-trace regularization, which motivates the study of the curvature-regularized objective of the form $f(x) + (\eta^2/2)\operatorname{tr}(\nabla^2 f(x))$. We develop a non-asymptotic theory in the nonconvex setting for a class of algorithms called stochastic zeroth-order Langevin dynamics (ZOLD), which use only noisy function evaluations and do not have access to first- or higher-order derivative information, and derive excess-risk guarantees for the curvature-regularized objective. We first control the finite-time distribution of the ZOLD iterates in Wasserstein distance with respect to the Gibbs measure associated with the Gaussian-smoothed objective. An explicit coupling between the smoothing and inverse temperature parameters then quantifies the proximity of this measure to the Gibbs measure associated with the corresponding Hessian-trace-regularized objective. These distributional estimates in Wasserstein distance lead to excess-risk bounds for the curvature-regularized objective. Our analysis accommodates assumptions under which the stochastic function evaluations are allowed to have non-uniform dependence on the data stream. In addition, the constants in the non-asymptotic bounds do not rely on log-Sobolev or Poincar{\'e} inequalities as in prior zeroth-order Langevin analyses. To the best of our knowledge, this is the first non-asymptotic analysis of a zeroth-order Langevin method that simultaneously establishes excess-risk guarantees for a curvature-regularized objective and Wasserstein guarantees for the associated Gibbs measure in the nonconvex setting.
\end{abstract}

\section{Introduction}
We consider the following nonconvex stochastic optimization problem:
\begin{equation}\label{original_objective_f}
\min_{x \in\mathbb{R}^{d}} f(x):=\mathbb{E}[F(x,\xi)],
\end{equation}
where $f:\mathbb{R}^{d}\to\mathbb{R}$ is a  smooth function, $\xi: \Omega \rightarrow \mathbb{R}^m$ denotes a random variable associated with the probability space $(\Omega, \mathcal{F}, \mathbb{P})$, and $F:\mathbb{R}^d \times \mathbb{R}^m \to \mathbb{R}$ is a measurable function satisfying $\mathbb{E}[|F(x,\xi)|] < \infty$, for all $x \in \mathbb{R}^d$.
Zeroth-order, or derivative-free, stochastic optimization methods provide a natural approach to solving \eqref{original_objective_f} when gradients are unavailable, expensive to compute, or only (noisy) function evaluations are accessible.
Such methods have found wide applicability in black-box
and bandit optimization
\cite{flaxman_online,Shamir13}, reinforcement learning and policy search
\cite{choromanski18a,salimans2017evolution}, hyperparameter optimization
\cite{snoek2012practical}, and black-box adversarial attacks
\cite{chen2019zo,andriushchenko2020square}.
 More recently,
zeroth-order methods have also been used for memory-efficient fine-tuning
of large language models, where forward-only evaluations avoid the memory and computational costs
associated with backpropagation
\cite{malladi2023fine,zhang2024dpzero}.

A substantial literature has developed theoretical guarantees for stochastic zeroth-order optimization, including nonconvex settings, with results typically formulated in terms of stationarity, objective suboptimality, or related measures of solution quality for the original objective $f$; see, e.g., \cite{ghadimi2013stochastic,ji2019improved,balasubramanian2022zeroth}. These analyses generally view the zeroth-order estimator as an approximation to the gradient of $f$. However, they typically do not characterize how the implicit bias induced by the zeroth-order estimator affects the solutions being targeted even though the smoothing mechanism intrinsic to the zeroth-order estimator modifies the effective objective being optimized. This distinction is relevant in nonconvex
optimization, where multiple globally or nearly optimal regions may coexist
and objective suboptimality alone does not capture the mechanism by which an
algorithm preferentially explores or selects among them.

For $\eta>0$, consider the
Gaussian-smoothed objective
\begin{equation}\label{smoothed_Gaussian_function}
f_{\eta}(x)
:=
\mathbb{E}[f(x+\eta U)], \qquad U \sim \mathcal{N}(0, I_d),
\end{equation}
where the expectation is taken with respect to $U$, a standard $d$-dimensional Gaussian random variable. A two-point (Gaussian)
zeroth-order estimator $G_{\eta}$ is given by
\begin{equation}
\label{difference_quotient_estimator_introduction}
    G_{\eta}(x, U, \xi):=   \frac{F(x+\eta U,\xi) - F(x,\xi)}{\eta} U.
\end{equation}
This estimator is unbiased for $\nabla f_\eta$, namely, $\mathbb{E}[G_{\eta}(x,U,\xi)]
=
\nabla f_{\eta}(x)$,\footnote{The expectation is taken with respect to the random variables $U$ and $\xi$ which are assumed to be independent of each other.} rather than for $\nabla f(x)$. 
We use
the geometric effect induced by the smoothing \eqref{smoothed_Gaussian_function} constructively rather than treating it as a source of approximation error that should be made as small as possible.  Under Assumptions \ref{independence_assumption_ZO} and \ref{higher_order_regularity_objective_assumption}, an application of Taylor's theorem to
\eqref{smoothed_Gaussian_function} yields 
\begin{equation}\label{Taylor_expansion_Gaussian_smoothed_function}
f_{\eta}(x)
=
f(x)
+
\frac{\eta^{2}}{2}
\operatorname{tr}\left(\nabla^{2}f(x)\right)
+
O(\eta^{4}).\footnote{The remainder term $O(\eta^{4})$ is uniform in $x$; see Lemma \ref{remainder_lemma}.}
\end{equation}
Thus, the first nontrivial correction induced by Gaussian smoothing is a
Hessian-trace term. This naturally leads to the curvature-regularized
objective
\begin{equation}\label{flat_loss}
v_{\eta}(x)
:=
f(x)
+
\frac{\eta^{2}}{2}
\operatorname{tr}\left(\nabla^{2}f(x)\right),
\end{equation}
which provides the leading-order approximation to 
$f_\eta$ beyond the original objective $f$.
Since the trace of the Hessian measures aggregate curvature across
directions, minimizing   $v_{\eta}$ balances the objective value
against a penalty on aggregate local curvature. We refer to the minimizers of $v_{\eta}$ as
curvature-regularized minima. This notion is related to the
broader concept of flat minima studied in the machine-learning literature,
where curvature-based criteria have been used to characterize the
preference of optimization methods for solutions with lower local
curvature \cite{ahn2023escape,zhang2025zeroth,bruno2026flatnessaware,hochreiter1997flat}. Importantly, the Hessian-trace regularization is not externally
imposed on the optimization problem \eqref{original_objective_f}. It arises intrinsically from the
smoothing mechanism used to construct the zeroth-order gradient estimator. Moreover,
$ |f_{\eta}(x)- v_{\eta}(x)|=O(\eta^{4})$,
whereas replacing $f_{\eta}$ directly by $f$ incurs an $O(\eta^{2})$ error.
Consequently, for small $\eta$, the curvature-regularized objective \eqref{flat_loss} provides
a more accurate leading-order description of the effective objective induced
by Gaussian smoothing \eqref{smoothed_Gaussian_function} than the original objective itself.

The recent work of Zhang et al.~\cite{zhang2025zeroth}, which provides a key motivation for our study, represents an important step toward a theoretical understanding of implicit regularization in zeroth-order optimization. 
They show that a zeroth-order method based on the standard two-point estimator favors minimizers with smaller Hessian trace and, under convexity and sufficient smoothness assumptions, provide convergence rates to minimizers that attain the smallest Hessian trace among all minimizers. Their analysis, however, is restricted to the convex setting. Moreover, as noted in \cite[Remark~3.3]{zhang2025zeroth}, extending their definition of minima to the general nonconvex setting remains theoretically challenging.

Langevin-based algorithms provide a useful framework for analyzing nonconvex optimization problems.
Consider the overdamped Langevin diffusion
\begin{equation}\label{Langevin_SDE_equation}
\rmd \bar{Y}_t
=
-\nabla f(\bar{Y}_t)\,\rmd t
+
\sqrt{2\beta^{-1}}\,\rmd B_t,
\end{equation}
where $(B_t)_{t\geq 0}$ is a $d$-dimensional Brownian motion and
$\beta>0$ is the inverse temperature. Under suitable assumptions\footnote{Assumptions \ref{independence_assumption_ZO}, \ref{Lipschitz_assumption_full_gradient}, and \ref{strong_convexity_at_infinity_assumption_ZO} are sufficient to guarantee the existence and uniqueness of $\pi_{\beta}^{f}$, see e.g. \cite{eberle2016reflection}.} on $f$, this
diffusion admits a unique invariant measure
\begin{equation} \label{original_objective_invariant_measure}
\pi_{\beta}^{f}(\rmd x)
\propto
\exp(-\beta f(x))\,\rmd x.
\end{equation}
For sufficiently large values of $\beta$, this measure concentrates around the global 
minimizers of $f$, see, e.g., \cite{hwang1980laplace,pmlr-v65-raginsky17a}. The Euler--Maruyama
discretization of \eqref{Langevin_SDE_equation}, with the exact gradient $\nabla f$
replaced by the stochastic gradient $\nabla_x F(x,\xi)$, gives rise to the class of algorithms known as Stochastic Gradient
Langevin Dynamics (SGLD) \cite{welling2011bayesian,pmlr-v65-raginsky17a}.

The combination of the implicit regularization induced by the zeroth-order estimator and the sampling-based mechanism of the Langevin dynamics leads to
the central object of this work. We introduce stochastic zeroth-order
Langevin dynamics (ZOLD) for solving the nonconvex optimization problem \eqref{original_objective_f} and obtain finite-time excess risk guarantees for its curvature-regularized objective using only noisy function evaluations $F(x,\xi)$. For this class of algorithms, the two-point Gaussian-smoothed zeroth-order estimator
\eqref{difference_quotient_estimator_introduction} replaces the stochastic gradient in the SGLD algorithm and its update rule is given by
\begin{equation}\label{eq:ZO_introduction}
x_{k+1}
=
x_k
-
\lambda G_{\eta}(x_k,U_{k+1},\xi_{k+1})
+
\sqrt{2\lambda\beta^{-1}}\, \zeta_{k+1}, \qquad k \in \mathbb{N},
\end{equation}
where $\lambda>0$ is the stepsize, $(\xi_k)_{k\in \N}$ is an independent sequence of random variables with the same distribution as $\xi$, and
$(\zeta_k)_{k\in \N}$ is a sequence of independent standard Gaussian random
vectors. The two Gaussian perturbations in the ZOLD recursion play distinct roles. The random variable $U_{k+1}$ is related to the Gaussian smoothing underlying the zeroth-order estimator and thereby is responsible for modifying the effective objective, inducing a preference toward low-curvature regions. In contrast, the Gaussian noise $\sqrt{2\lambda\beta^{-1}}\zeta_{k+1}$ arises from the Langevin dynamics \eqref{Langevin_SDE_equation} and provides a mechanism for global exploration.
The underlying Langevin diffusion associated with \eqref{eq:ZO_introduction}, i.e. $\rmd Y_t=
-\nabla f_{\eta}(Y_t)\,\rmd t
+
\sqrt{2\beta^{-1}}\,\rmd B_t$, admits a unique invariant
measure
\begin{equation}\label{Gaussian_smoothed_invariant_measure}
\pi_{\beta}^{f_{\eta}}(\rmd x)
\propto
\exp(-\beta f_{\eta}(x))\,\rmd x.
\end{equation}
For large values of $\beta$, this Gibbs measure concentrates around the global minimizers of $f_{\eta}$. The main question is therefore whether the finite-time distribution of ZOLD \eqref{eq:ZO_introduction}
can be quantitatively related to the Gibbs measure associated with the curvature-regularized objective 
\begin{equation}\label{Hessian-trace_regularized_objective}
\pi_{\beta}^{v_{\eta}}(\rmd x)
\propto
\exp(-\beta v_{\eta}(x))\,\rmd x,
\end{equation}
 thereby connecting the ZOLD algorithm to the minimizers of $v_{\eta}$.
Our analysis establishes this connection through two steps.
First, we derive non-asymptotic Wasserstein bounds between the law of the
finite-time ZOLD iterates and the invariant measure \eqref{Gaussian_smoothed_invariant_measure}. Then, we quantify the distance
between the invariant measure $\pi_{\beta}^{f_{\eta}}$ and the Gibbs measure
$\pi_{\beta}^{v_{\eta}}$ associated with the Hessian-trace-regularized
objective. An explicit coupling between the smoothing
scale $\eta$ and the inverse temperature $\beta$ allows us to control this distance for large values of $\beta$.
After suppressing logarithmic
factors and the dependence on problem-specific parameters, the resulting
iteration complexities for the distance between the distribution of the ZOLD iterate $\mathcal{L}(x_k)$ and $\pi_{\beta}^{v_{\eta}}$ scale in Wasserstein distance as
\begin{equation*}
    \begin{split} \widetilde{\mathcal{O}}(d^2\epsilon^{-2} \bar C(\beta))
\ \text{for }  W_1(\mathcal{L}(x_k),\pi^{v_{\eta}}_{\beta} ) \le \epsilon,
\quad  \widetilde{\mathcal{O}}(d^4\epsilon^{-4}\widehat{C}(\beta) )
\ \text{for }W_2(\mathcal{L}(x_k),\pi^{v_{\eta}}_{\beta} ) \le \epsilon,
    \end{split}
\end{equation*}
where the dependence in $\beta$ in the constants $\bar C(\beta)$ and $\widehat{C}(\beta)$ is provided in Theorem \ref{Theorem1_ZO} and Corollary \ref{Theorem1_ZO_W_2}, respectively.
Finally, we translate these distributional estimates into excess-risk bounds for $v_{\eta}$.

Our analysis applies to nonconvex objectives, including multimodal potentials such as Gaussian-mixture objectives, under a
strong-convexity-at-infinity condition and first- and third-order smoothness assumptions on $F(x, \xi)$, and avoids the use of log-Sobolev or Poincar{\'e} inequalities.  This distinction is important for nonconvex optimization: functional inequalities provide powerful sampling tools, but their constants can deteriorate severely with dimension and inverse temperature and may be difficult to compute for multimodal objectives. Our approach instead exploits contraction properties implied by the geometry of the objective (Assumption \ref{strong_convexity_at_infinity_assumption_ZO}), allowing us to obtain polynomial dependence on the dimension in the non-asymptotic bounds without relying on log-Sobolev or Poincar{\'e} constants. This avoids
the potentially exponential dimension dependence inherited from such constants in the nonconvex setting, which appears in prior analyses of zeroth-order Langevin algorithms \cite{roy2022stochasticLangevin,nakakita2023non,naldi2026langevin,sahin2026zerothorder,liuone}; see Remarks~\ref{rem:dim_dependence1} and~\ref{rem:dim_dependence2}. Our framework also accommodates assumptions in which the stochastic oracle is allowed to exhibit non-uniform dependence on the data, thereby relaxing the uniform derivative conditions commonly imposed in existing zeroth-order (Langevin) analyses. Lastly, our guarantees are obtained with a single two-point zeroth-order estimator per iteration, so that, unlike in \cite{roy2022stochasticLangevin,nakakita2023non,naldi2026langevin}, the batch size does not need to grow as the target accuracy improves; see Section~\ref{sec:comparison}. This is particularly relevant in large-scale machine learning applications, where each function evaluation is costly, such as zeroth-order fine-tuning of large language models, in which every evaluation requires a forward pass through a model with billions of parameters \cite{malladi2023fine}, and black-box adversarial attacks, in which the number of queries is limited \cite{andriushchenko2020square}.

To the best of our knowledge, this is the first non-asymptotic analysis of a
zeroth-order Langevin algorithm that simultaneously provides
  excess-risk guarantees
for a curvature-regularized objective and distributional guarantees for the finite-time law and Gibbs measure associated with this objective in a nonconvex setting.

\section{Preliminaries} \label{sec:preliminaries} In this section, we provide the notation used throughout this paper, the definition of curvature-regularized minimum, and the zeroth-order Langevin algorithm.

\paragraph{Notation} Fix integers $d,m \ge 1$, and let $(\Omega, \mathcal{F}, \mathbb{P})$  be a fixed
	probability space.  We denote the probability law of a random variable $Y$ by $\mathcal{L}(Y)$. We write $U\sim\mathcal{N}(0,I_d)$ for a standard $d$-dimensional Gaussian random variable.  The Euclidean scalar product is denoted by $\langle \cdot, \cdot \rangle $, with $| \cdot |$ standing for the corresponding norm. The notation $\widetilde{\mathcal{O}}(\cdot)$ is used to omit logarithmic 
factors in the relevant parameters. Let $f: \mathbb{R}^{d } \rightarrow \mathbb{R}$ be a continuously differentiable function, and we denote its gradient by $\nabla f$. Let $\mathcal{P}(\mathbb{R}^d)$ be the set of probability measures on $\mathcal{B}(\mathbb{R}^d)$ (the Borel sigma-algebra of $\mathbb{R}^d$). For $\tilde{\nu}$, $\bar{\nu} \in \mathcal{P}(\mathbb{R}^d)$, let $\mathcal{C}(\tilde{\nu}, \bar{\nu})$ denote the set of probability measures $\varrho$ on $\mathcal{B}(\mathbb{R}^{2d })$ such that its respective marginals are $\tilde{\nu}$ and $\bar{\nu}$. For any $\tilde{\nu}$ and $\bar{\nu} \in \mathcal{P}(\mathbb{R}^{d })$ with finite $p$th moments, the Wasserstein distance of order $p \ge 1$ is defined as
	\begin{equation*} %\label{Wasserstein_distance_notation}
		W_p(\tilde{\nu}, \bar{\nu})= \left(  \inf_{\varrho \in \mathcal{C}(\tilde{\nu}, \bar{\nu})} \int_{\mathbb{R}^d}  \int_{\mathbb{R}^d}  |x - z|^p \ \text{d} \varrho(x,z) \right)^{\frac{1}{p}}.
	\end{equation*}
   For any  $\tilde{\nu}$ and $ \zeta \in \mathcal{P}(\mathbb{R}^{d })$, if $\tilde{\nu}$ is absolutely continuous with respect to $\zeta$, we denote by $\rmd \tilde{\nu} / \rmd \zeta $ the Radon-Nikodym derivative of $\tilde{\nu}$ with respect to $\zeta $, and the Kullback-Leibler divergence of  $\tilde{\nu}$ with respect to $\zeta $ is given by
   \begin{equation*}
       \text{KL}(\tilde{\nu} | | \zeta) 
			 =    \int_{\mathbb{R}^d}  \frac{\rmd \tilde{\nu} }{  \rmd \zeta}  \log \left( \frac{ \rmd \tilde{\nu} }{  \rmd \zeta}  \right)  \rmd \zeta .
   \end{equation*}

\begin{definition} \label{curvature_regularized_minimum_definition}
Suppose that $f$ is twice differentiable with a well-defined Hessian matrix and ${arg\,min}_{x\in\mathbb{R}^d} v_{\eta}(x)$ is nonempty,\footnote{This is satisfied under Assumptions \ref{independence_assumption_ZO}, \ref{Lipschitz_assumption_full_gradient}, \ref{strong_convexity_at_infinity_assumption_ZO}, and \ref{higher_order_regularity_objective_assumption}.} with $v_{\eta}(x)= f(x) +\frac{\eta^{2}}{2}
\operatorname{tr}(\nabla^{2}f(x))$, for $\eta >0$. We say that $x^{\star}$ is a curvature-regularized minimum of $f(x)$ if $x^{\star} \in \operatorname*{arg\,min}_{x\in\mathbb{R}^d} v_{\eta}(x)$.
\end{definition}

\subsection{ZOLD Algorithm} \label{sec:ZO_Algorithm} In the following, we introduce the stochastic zeroth-order Langevin dynamics (ZOLD) algorithm.
For $x , u \in \mathbb{R}^d$ and $z \in \mathbb{R}^m$, we define the two-point finite-difference estimator $ G_{\eta} : \mathbb{R}^{d} \times \mathbb{R}^{d} \times  \mathbb{R}^{m} \rightarrow \mathbb{R}^d $   as
\begin{equation} \label{difference_quotient_estimator}
    G_{\eta}(x, u, z)=   \frac{F(x+\eta u,z) - F(x,z)}{\eta} u, 
    %\qquad \text{for} \ \eta >0,
\end{equation}
where $ \eta >0$ denotes the smoothing parameter.   
%the central difference
%\begin{equation} \label{symmetric_difference_quotient_estimator}
 %   G_{\eta}(\theta, u, x):=   \frac{F(\theta+\eta u,x) - F(\theta-\eta u,x)}{2\eta} u, \qquad \text{for} \ \eta >0.
%\end{equation}
In particular, 
\begin{equation} \label{biased_estimate_two_point_estimator}
    \E [ G_{\eta}(x, U, \xi)] = \nabla f_{\eta}(x),
\end{equation}
where the expectation is taken with respect to the random variable $\xi$ with probability law $\mathcal{L}(\xi)$ and the standard $d$-dimensional Gaussian random variable $U $.

Let $(\xi_k)_{k \in \mathbb{N}}$ be independent and identically distributed with $\xi_k \sim \mathcal{L}(\xi)$ for all $k \in \mathbb{N}$, and let $(U_{k})_{k \in \mathbb{N}}$ and $(\zeta_{k})_{k \in \mathbb{N}}$ be independent sequences of standard $d$-dimensional Gaussian random variables. Let $x_0$ be an $\mathbb{R}^d$-valued random variable representing the initial value of the ZOLD algorithm. We assume that $x_0$,  $(\xi_k)_{k \in \mathbb{N}}$, $(U_k)_{k\in \mathbb{N}}$,  and $(\zeta_k)_{k\in \mathbb{N}}$ are all mutually independent.
%Let $(\Xi_k)_{k \in \N}$ be a given filtration representing the flow of past information and denote by $\Xi_{\infty}$ the sigma-algebra generated by $\cup_{k \in \mathbb{N}} \Xi_k$. Let  $(\xi_{k})_{k \in \mathbb{N}}$ be an $\mathbb{R}^m$-valued $(\Xi_k)$-adapted process with $\xi_k \sim \mathcal{L}(\xi)$ for all $k \in \mathbb{N}$.
%Let $x_0$ be an $\mathbb{R}^d$-valued random variable representing the initial value of the ZOLD algorithm,  
%$(U_{k})_{k \in \mathbb{N}}$ be i.i.d copies of the standard $d$-dimensional Gaussian random variable $U$, and $(\zeta_{k})_{k \in \mathbb{N}}$ be an independent sequence of standard $d$-dimensional Gaussian random variables. We assume that $x_0$, $(U_k)_{k\in \mathbb{N}}$,  $\Xi_{\infty}$, and $(\zeta_k)_{k\in \mathbb{N}}$ are all independent.
Then the ZOLD algorithm is given by
\begin{equation} \label{eq:ZO}
\begin{split}
%\begin{cases}
%	\theta_0^{\text{ZO}} &:= \theta_0,
%    \\ 
  x_{k+1}
   = x_k
    - \lambda    G_{\eta}(x_k, U_{k+1}, \xi_{k+1})
    + \sqrt{2 \lambda \beta^{-1}} \zeta_{k+1}, \qquad k \in \mathbb{N}.
%    \end{cases}
\end{split}
\end{equation}
For the analysis of ZOLD, we use the following coupling between the
smoothing parameter $\eta$ and the inverse-temperature parameter $\beta$:
\begin{equation}\label{coupling}
    \eta = \beta^{-\frac{1+\tau}{4}}d^{-1}, \qquad  \text{for} \quad  \tau \in (0,1).
\end{equation}
This choice is motivated by the fourth-order remainder 
estimate
%Lemma~\ref{remainder_lemma} and
in the expansion
\eqref{Taylor_expansion_Gaussian_smoothed_function}. In particular, it
ensures that the Gibbs measure
$\pi_{\beta}^{f_\eta}$ associated with the Gaussian-smoothed objective remains sufficiently close for large values of $\beta$ to
the Gibbs measure $\pi_{\beta}^{v_{\eta}}$ associated with the
Hessian-trace-regularized objective  (see, e.g. Proposition \ref{Wasserstein_KL_convergence_two_measures_theorem_ZO} and Remark \ref{Wasserstein_KL_convergence_two_measures_theorem_ZO_remark}). The coupling also
reduces the number of independent tuning parameters in the ZOLD
algorithm.  Once $d$ and $\tau$ are fixed, $\eta$ is determined by $\beta$,
so that only the inverse-temperature parameter $\beta$ needs to be tuned
independently. The parameter
$\tau\in(0,1)$ is introduced only for technical purposes in the analysis.

\section{Assumptions}  \label{sec:assumptions}

We state the main assumptions of this work.
\begin{assumption}\label{independence_assumption_ZO}
Let $f: \mathbb{R}^d \rightarrow \mathbb{R}$ be a three times continuously differentiable function.  Assume that the random variables $x_0$ with $\E [|x_0|^2] < \infty$, $(U_k)_{k\in \mathbb{N}}$, $(\xi_k)_{k\in \mathbb{N}}$ and $(\zeta_k)_{k\in \mathbb{N}}$ defined in Section \ref{sec:ZO_Algorithm} are all mutually independent, and that the stochastic oracle   satisfies $\E [  F(x, \xi)]= f(x)$ and $\E [\nabla_{x}  F(x, \xi)]= \nabla f(x)$ for every $x \in\mathbb{R}^d$. 
\end{assumption}
\begin{assumption}[Lipschitzness] \label{Lipschitz_assumption_full_gradient}
There exists $\varphi : \mathbb{R}^{m} \rightarrow [1, \infty)$ with $\E [|(1+ | \xi|) \varphi(\xi)|^2] < \infty $, and $L_1,L_2 >0$ such that, for all $z$, $z' \in \mathbb{R}^m$ and $x$, $x' \in \mathbb{R}^d$,
    \begin{equation*}
    \begin{split}
        | \nabla_{x} F(x, z)  -  \nabla_{x'} F(x', z)|  & \le L_1 \varphi (z) |x - x'|, \qquad \text{and}
   \\  | \nabla_{x} F(x, z)  -  \nabla_{x} F(x, z')| &  \le L_2 (\varphi (z) + \varphi(z') )
   (1+ | x |) |z - z'|.
        \end{split}
    \end{equation*}
\end{assumption}
\begin{assumption}[Strong convexity-at-infinity]\label{strong_convexity_at_infinity_assumption_ZO}  There exist $\mu, R>0$ and a measurable function  $\rho : \mathbb{R}^{m} \rightarrow [1, \infty)$ with $\E[|\rho(\xi)|]< \infty$, such that, for all $z \in \mathbb{R}^m$ and for all $x$, $ x' \in\mathbb{R}^d$,
\begin{equation*} %\label{strong_convexity_assumption_inequality_ZO}
    \langle \nabla_{x} F(x,z)-  \nabla_{x'} F(x',z) , x - x' \rangle \geq \mu \rho(z) |x  - x'|^2,\  \text{with} \  |x - x'|\geq R.
\end{equation*}
\end{assumption}
Under Assumption
\ref{strong_convexity_at_infinity_assumption_ZO}, the nonconvex behavior of $F(x,z)$ is
allowed within regions whose relevant point pairs are separated by
distances smaller than $R$.
The non-asymptotic Wasserstein bounds between the law of the ZOLD iterates and the Gibbs measure associated with the Gaussian-smoothed objective $\pi_{\beta}^{f_{\eta}}$ are established under Assumptions \ref{independence_assumption_ZO},  \ref{Lipschitz_assumption_full_gradient}, and \ref{strong_convexity_at_infinity_assumption_ZO}.
Importantly, Assumptions \ref{Lipschitz_assumption_full_gradient} and \ref{strong_convexity_at_infinity_assumption_ZO} are crucial in proving the contraction estimates in Wasserstein distances for the Langevin dynamics \eqref{zeroth_order_Langevin_SDE} and the auxiliary process \eqref{before_continuous_time_process_started_interpolation_definition_ZO}-\eqref{continuous_time_process_started_interpolation_definition_ZO}; see, e.g., Proposition \ref{prop2_ZO_before} and Lemma \ref{Lemma7_ZO_before}.

We impose the following additional regularity assumption to control the distance between $f_{\eta}$ and $v_{\eta}$; see, e.g., Lemma \ref{remainder_lemma}.
\begin{assumption}\label{higher_order_regularity_objective_assumption}
   There exist $L_3>0$ and a measurable function $\bar{\varphi}: \R^m \rightarrow [1,\infty)$ with $\E [|\bar{\varphi}(\xi)|]< \infty$, such that for all $z \in \R^m$ and  for all $x , x' \in \R^d$,
\begin{equation*} %\label{assumption_large_values_remainder}
\begin{split}
  & | F(x',z) - F(x,z) - \langle \nabla F(x,z), x' - x \rangle  \\ & \quad - \frac{1}{2} \langle \nabla^2 F(x,z) (x' - x) , x' - x \rangle  - \frac{1}{6} \nabla^3 F(x,z) [x' - x]^3 | 
   \\ & \le \frac{L_3 \bar{\varphi}(z)}{24} | x - x' |^4,
   \end{split}
\end{equation*}
where $ \nabla^3 F(x,z) [x' - x]^3 = \sum_{i,j,h \in [d]} (\partial^3 F(x,z) / \partial x_i \partial x_j \partial x_h)(x'_i - x_i)(x'_j - x_j)(x'_h - x_h)$.
\end{assumption}
\begin{remark}
    Because the curvature-regularized objective contains second-order information of $f$, i.e., $\frac{\eta^2}{2}\operatorname{tr} (\nabla^2 f(x))$ in \eqref{flat_loss}, it is natural for our analysis to impose higher-order regularity on $f$ (Assumption \ref{independence_assumption_ZO}) and on $F(x,z)$ (Assumption \ref{higher_order_regularity_objective_assumption}). Such higher-order regularity assumptions are common in analyses of implicit regularization \cite{arora2022understanding,wen2023how} and in  complexity analyses of sharp and flat minima in machine learning \cite{ahn2023escape,zhang2025zeroth,bruno2026flatnessaware}. The ZOLD algorithm \eqref{eq:ZO} itself uses only zeroth-order information through the estimator $G_{\eta}$; the higher-order regularity assumptions on $F(x,z)$ are required only for the non-asymptotic analysis.
\end{remark}
\begin{remark}\label{non_uniform_dependence_data_stream_assumptions}
 Assumptions \ref{Lipschitz_assumption_full_gradient}, \ref{strong_convexity_at_infinity_assumption_ZO}, and \ref{higher_order_regularity_objective_assumption} impose regularity conditions with data-dependent coefficients on $\nabla_x^i F(x,z)$, $i=1,2,3$, allowing the bounds to depend on $z$. This significantly relaxes the analogous requirements in the zeroth-order literature, where the derivatives of $ F(\cdot,z)$ are typically assumed to admit bounds uniform in $z$; see, e.g., \cite{roy2022stochastic,naldi2026langevin,roy2022stochasticLangevin,balasubramanian2022zeroth}. 
\end{remark}
\begin{remark}
Our assumptions cover the case of multimodal Gaussian-mixture potentials. %suitable
%component means and mixing weights, the associated Gibbs distribution is
%non-log-concave. 
For such potentials,  the Lipschitz-gradient and strong-convexity-at-infinity
conditions follow from
\cite[Proposition~A.1]{silveri2025beyond} and
\cite[Proposition~13 and Section~3.4.2]{bruno2025wasserstein}.
Moreover, in the symmetric two-component Gaussian-mixture case, the
higher-order regularity in Assumption~\ref{higher_order_regularity_objective_assumption}
holds whenever the third derivative of the potential is globally Lipschitz;
in particular, it holds if the fourth derivative is globally bounded.
\end{remark}

\section{Main results}
\label{sec:main}
In this section, we state the main results of this work in Theorem  \ref{Theorem1_ZO}, Corollary \ref{Theorem1_ZO_W_2}, and Theorem \ref{Theorem3_ZO}. The proofs are postponed to Section \ref{sec:proofs_main_results}. 
First, we define the maximum stepsize
\begin{equation} \label{expression_lambda_max_ZO}
    \begin{split}
        \lambda_{\text{max}}=\min\left \{\frac{a}{4(4+d) L_1^2 \E[\varphi^2(\xi)]},\frac{1}{a},1 \right \},
    \end{split}
\end{equation}
with $a$ given in \eqref{constants_dissipativity_flatness_full_gradient_ZO}.  The condition \eqref{expression_lambda_max_ZO} is obtained by controlling the second moment of the ZOLD iterates; see, e.g., the proof of Lemma \ref{Lemma2_ZO}. The dependence on the dimension $d$ in \eqref{expression_lambda_max_ZO} is standard in zeroth-order optimization; see, e.g., \cite{nesterov2017random,ghadimi2013stochastic,lan2020first}.

\begin{theorem}\label{Theorem1_ZO}
Let Assumptions \ref{independence_assumption_ZO},  \ref{Lipschitz_assumption_full_gradient},
\ref{strong_convexity_at_infinity_assumption_ZO}, and \ref{higher_order_regularity_objective_assumption}  hold, let $\epsilon >0$ be the target accuracy, and let $\eta = \beta^{-\frac{1+\tau}{4}}d^{-1}$  for $\tau \in (0,1)$. Then, by choosing $\beta \ge \beta_{\epsilon}$ with $\beta_{\epsilon}$ given in \eqref{definition_beta_bar_W_1_ZOLD},  $0< \lambda < \min \left \{\lambda_{\text{max}}, \frac{\epsilon^2}{9 C_{6,1}^2} \right \}$ with $\lambda_{\text{max}}$  and $C_{6,1}$ given in \eqref{expression_lambda_max_ZO} and \eqref{constant_C_6_1_main_proof} respectively, we achieve 
\begin{equation*}
    W_1(\mathcal{L}(x_K),\pi^{v_{\eta}}_{\beta} ) \le \epsilon,
\end{equation*}
 after $K = \widetilde{\mathcal{O}} ( d^2 \epsilon^{-2}   \beta^{-3/2}  e^{\beta R^2 }  )$ iterations.
\end{theorem}

By additionally using a simple functional inequality, a corresponding result can be formulated in the Wasserstein-2 distance and used for the derivation of the  expected excess risk bound in Theorem \ref{Theorem3_ZO}.
\begin{corollary}\label{Theorem1_ZO_W_2}
Let Assumptions \ref{independence_assumption_ZO},  \ref{Lipschitz_assumption_full_gradient}, \ref{strong_convexity_at_infinity_assumption_ZO}, and \ref{higher_order_regularity_objective_assumption}    hold,  let $\epsilon >0$ be the target accuracy,   and let $\eta = \beta^{-\frac{1+\tau}{4}} d^{-1}$  for $\tau  \in (0,1)$.  Then, by choosing $\beta $ as in Theorem \ref{Theorem1_ZO}, $ 0< \lambda <   \min \left \{\lambda_{\text{max}}, \frac{\epsilon^4}{81 C_{6,2}^4} \right \}$ with $\lambda_{\text{max}}$  and $C_{6,2}$ given in \eqref{expression_lambda_max_ZO} and \eqref{constant_C_6_2_main_proof} respectively, we achieve 
\begin{equation*}
    W_2(\mathcal{L}(x_K),\pi^{v_{\eta}}_{\beta} ) \le \epsilon,
\end{equation*}
 after $ K =  \widetilde{\mathcal{O}} ( d^4   \epsilon^{-4} \beta^{- 21/4} e^{R^2\beta^{3/2}}    )$ iterations.
\end{corollary}
\begin{remark}\label{rem:dim_dependence1}
The non-asymptotic bounds derived for Theorem \ref{Theorem1_ZO} and Corollary \ref{Theorem1_ZO_W_2} are 
\begin{equation} \label{sampling_bound_two_main_results_remark}
    W_p (\mathcal{L}(x_K),\pi_{\beta}^{v_{\eta}} )  \leq C_{W_p}e^{-C_{r_p}\lambda K} \Delta_0^{(p)} +C_{6,p}\lambda^{1/(2p)} +  C_5,  \quad p \in \{ 1, 2\},
\end{equation}
where the first term on the right-hand side above is the exponential contraction term for the underlying Langevin diffusion, the second term accounts for the Euler--Maruyama discretization error, and the third term is the bound for the Wasserstein distance between $\pi_{\beta}^{f_{\eta}} $ and $\pi_{\beta}^{v_{\eta}}$.
We emphasize that the constants $C_{W_p}$ and $C_{r_p}$ do not depend on the dimension $d$ (Proposition \ref{prop2_ZO_before}). They follow from the contraction estimates in $W_p$ established under 
Assumptions~\ref{Lipschitz_assumption_full_gradient} and
\ref{strong_convexity_at_infinity_assumption_ZO}. Consequently, the dimension enters only through $C_{6,1}$, $C_{6,2}$ and $C_5$ (see, e.g., \eqref{constant_C_6_1_main_proof_asymptotic_notation}, \eqref{constant_C_6_2_main_proof_asymptotic_notation}, and Remark \ref{Wasserstein_KL_convergence_two_measures_theorem_ZO_remark}). This dependence arises through the moment estimates for the
ZOLD iterates, which are obtained from the dissipativity of $f_{\eta}$
(Remark~\ref{dissipativity_assumption_remark_ZO}) and the second moment of the zeroth-order estimator
(Lemma \ref{linear_growth_two_point_estimator_lemma}), and through the bound on the distance between $\pi_{\beta}^{f_{\eta}}$ and $\pi_{\beta}^{v_{\eta}}$ (Proposition \ref{Wasserstein_KL_convergence_two_measures_theorem_ZO}). Therefore, all constants in the above bound depend at most polynomially on $d$, whereas the prior zeroth-order Langevin analyses inherit the potentially exponential dimension dependence of log-Sobolev or Poincar\'e constants, see, e.g., \cite{roy2022stochasticLangevin,nakakita2023non,liuone,naldi2026langevin,sahin2026zerothorder}. Moreover, our constants in \eqref{sampling_bound_two_main_results_remark} also depend on $\beta$, $R$, $\mu$, $L_1$, $L_2$, and $L_3$, as is typical in nonconvex analyses.
\end{remark}
The following theorem provides the non-asymptotic expected excess risk bound with
respect to the curvature-regularized objective $v_{\eta}$. 
\begin{theorem}\label{Theorem3_ZO}
Let Assumptions \ref{independence_assumption_ZO}, \ref{Lipschitz_assumption_full_gradient}, \ref{strong_convexity_at_infinity_assumption_ZO}, and \ref{higher_order_regularity_objective_assumption} hold, let $\epsilon \in (0,1] $ be the target accuracy, and let $\eta = \beta^{-\frac{1+\tau}{4}} d^{-1}$  for $\tau \in (0,1)$. Then, by choosing $\beta \ge \beta_{\epsilon}$ with $\beta_{\epsilon}$ given in \eqref{beta_delta_expression_excess_risk_bound}, $ 0 < \lambda < \min \left \{\lambda_{\text{max}}, \frac{\epsilon^4}{81(C^{\#}_2)^4}  \right \}$ with $\lambda_{\text{max}}$ and $C^{\#}_2$ given in \eqref{expression_lambda_max_ZO} and \eqref{constants_excess_risk_bounds} respectively, we achieve
%for any $\beta\geq \max\{2/a, M^{-1}\}$ with $M$ given in \eqref{lemmaT2_constant_M_ZO}, $\lambda\in(0,\lambda_{\text{max}})$ with $\lambda_{\text{max}}$ given in \eqref{expression_lambda_max_ZO} and $K \in\mathbb{N}$, 
\begin{equation*} %\label{expected_excess_risk_g_epsilon_proof_semiconvexity_final}
     \E  [v_{\eta}(x_K)] - \inf_{x \in \mathbb{R}^d} v_{\eta}(x)   \le \epsilon, 
%   \\ & \leq  \Bigg( 2 L_1 \E[\varphi(X)] (   \sqrt{C_1\left(1+\mathbb{E}|\theta_0|^2\right)} 
%   \\ & \qquad + 2 L_1 \E[\varphi(X)]  \sqrt{1+ b a^{-1} +d (a \beta)^{-1}} 
%   \\ & \qquad + \frac{2 L_1 \E[\varphi(X)]}{  \beta^{\frac{1+\tau}{4}} d^{1/2}}  +  2 L_2  \E [\widetilde{\varphi}(X)]  + 2 |\nabla F(0,0)| \Bigg) W_2(\mathcal{L}(\theta_K),\pi_{\beta})
 %  \\ & \quad \times (C_{W_2}e^{-C_{r_2}\lambda K} \max\{ W_2(\mathcal{L}(\theta_0),\pi_{\beta}), \sqrt{W_1( \mathcal{L}(\theta_0),\pi_{\beta})} \}+C_{6,2}\lambda^{1/4}) 
  % \\  & \qquad  + \frac{d}{2\beta}\log\left(\frac{2\pi e(b+d/\beta)}{a d}\right) - \frac{1}{\beta} \log \left( \frac{2 \pi^{\frac{d}{2}}}{d \  \Gamma (\frac{d}{2})} \right) 
 % \\ & \qquad  + \frac{2}{\beta}   +    \frac{d}{\beta} \log(\beta M) 
% +  \frac{  3 L_3 \E [\bar{\varphi}(X)]}{ \beta^{1+\tau} d^2},
\end{equation*}
after $  K = \widetilde{\mathcal{O}}( d^8 \epsilon^{-4} \beta^{- 21/4} e^{R^2\beta^{3/2}}   ) $ iterations.
%\begin{equation*}
%  K = \widetilde{\mathcal{O}}( d^8 \epsilon^{-4} \beta^{- 21/4} e^{R^2\beta^{3/2}}   )    \quad \text{iterations}.
%\end{equation*}
%with $C^{\diamond}(\beta, \tau, R):= $.
% \frac{C_{W_2}^4}{ \beta^{1+ \tau} C_{r_2}\left(1-e^{-C_{r_2}/2}\right)^4}$, where the constants $C_{W_2}=\mathcal{O}(\beta^{-1}e^{R^2\beta^{3/2}})$ and $C_{r_2}=\mathcal{O}(e^{-R^2\beta^{3/2}})$ are provided explicitly in Proposition \ref{prop2_ZO_before}.
%the constant $C^{\#}_2 =\mathcal{O}\left((1+ \left(\frac{d}{\beta} \right)\frac{C_{W_2}}{1-e^{-C_{r_2}/2}}\right)$ the constant $C_{r_2}=\mathcal{O}(e^{-R^2\beta^{3/2}})$ do not depend on the dimension $d$, and is provided explicitly in Proposition \ref{prop2_ZO_before}.
%where $C_1$ is given in \eqref{eqProp2_constants},  $a$, $b$ are provided in Remark \ref{dissipativity_assumption_remark_ZO}, and $\widetilde{\varphi} = (\varphi (x) + \varphi (0))  | x|$.
%, the constants $C_{r_2}$ and $C_{W_2}$ are given in Corollary \ref{Theorem1_ZO_W_2}, and $C_{6,2}$ is given in \eqref{constant_C_6_2_main_proof}.
\end{theorem}
\begin{remark}\label{rem:dim_dependence2}
%The expected excess-risk bound in Theorem~\ref{Theorem3_ZO} is stated with
%respect to the curvature-regularized objective $v$.
The non-asymptotic bound derived for Theorem \ref{Theorem3_ZO} is
\begin{equation}
\label{sampling_bound_expected_excess_risk_main_results_remark}
        \E  [v_{\eta}(x_K)] - \inf_{x \in \mathbb{R}^d} v_{\eta}(x) \le C_1^{\#} e^{-C_{r_2}\lambda K}  +  C_2^{\#} \lambda^{1/4} + C_3^{\#},
\end{equation}
where the constants $C_1^{\#}$, $C_2^{\#}$, and $C_3^{\#}$ have a polynomial dependence in $d$ (see, e.g., \eqref{constants_excess_risk_bounds_C_1}, \eqref{constants_excess_risk_bounds}, and \eqref{constants_excess_risk_bounds_C_3}, respectively), while the constant $C_{r_2} = \mathcal{O}(e^{-R^2\beta^{3/2}})$, see, e.g., Proposition \ref{prop2_ZO_before}. 
In contrast to prior expected excess risk bounds for zeroth-order Langevin algorithms, see e.g., \cite[Corollary~4.8]{naldi2026langevin}, our analysis does not depend on any log-Sobolev
or Poincar{\'e} constants.
However, as in other analyses of Langevin-based algorithms  \cite{naldi2026langevin,zhang2023nonasymptotic}, the inverse temperature must be chosen
sufficiently large to control the Gibbs-concentration error. In Theorem \ref{Theorem3_ZO}, the choice of $\beta$ is determined by the bound on the concentration of the Gibbs measure $\pi_\beta^{f_{\eta}}$ and the distance between $f_\eta$ and $v_{\eta}$ (Lemma \ref{remainder_lemma}). 
\end{remark}

\subsection{Comparison with related work}\label{sec:comparison}

Our work lies at the intersection of three lines of research:
zeroth-order Langevin methods for sampling, zeroth-order Langevin methods for nonconvex optimization, and
implicit regularization induced by zeroth-order estimators.  Table \ref{tab:comparison_related_work} provides a comparison between our work and the results available within these three lines of research. 
\begin{table}[t]
\centering
\caption{Comparison of zeroth-order Langevin methods and curvature-regularization
approaches in zeroth-order optimization. The table indicates whether each
method uses Langevin dynamics, accommodates nonconvex objectives,
considers some form of curvature-regularization objective, provides distributional guarantees, and whether the batch size depends on the target accuracy $\epsilon$.}
\label{tab:comparison_related_work}
\small
\begin{tabular}{lcccccc}
\toprule
Reference 
& Langevin 
& Nonconvex 
& \makecell{Curvature- \\ regularized \\ target}
& \makecell{Distributional\\metric}
& \makecell{Batch \\ size \\ ($\epsilon$ dep.)}
%& \makecell{Optimization \\ target} 
 \\
\midrule
\cite[Theorem 3.2]{roy2022stochasticLangevin}
& Yes
& Yes
%& Sampling
& No
& $W_2$  
& $\mathcal{O}(\epsilon^{-2})$ \\
%& No \\

\cite[Section 3.4]{nakakita2023non}
& Yes
& Yes
%& Sampling
& No
& $W_2$ 
& $\widetilde{\mathcal{O}}(\epsilon^{-4})$ \\
%& No \\

\cite[Corollary 1]{sahin2026zerothorder}
& Yes
& Yes
%& Sampling
& No
& Squared TV
&  $\mathcal{O}(1)$ \\
%& No \\

\cite[Corollary 4.8]{naldi2026langevin}
& Yes
& Yes
%& \checkmark
& No
& KL 
& $\widetilde{\mathcal{O}}(\epsilon^{-2})$ \\
%& $f$ \\

\cite[Corollary 2]{zhang2025zeroth}
& No
& No
%& \checkmark
& Yes
& None 
& $\mathcal{O}(1)$  \\
%& \makecell{Flatness \\ criterion\footnote{Zhang et al.~\cite{zhang2025zeroth} control objective suboptimality and a Hessian-trace-related quantity.}}
%\\

Ours
& Yes
& Yes
%& \checkmark
& Yes
& $ W_1, W_2 $ 
& $\mathcal{O}(1)$ \\
%& $v$ \\
\bottomrule
\end{tabular}
\end{table}

%\subsubsection{Zeroth-order Langevin sampling}

The first line of work studies zeroth-order Langevin methods for sampling from the target distribution $\pi_\beta^f(\rmd x) $ given in \eqref{original_objective_invariant_measure}. Roy et al.\
\cite{roy2022stochasticLangevin} establish non-asymptotic Wasserstein approximation guarantees for $\pi_1^f$ for a zeroth-order Langevin Monte Carlo
(ZO-LMC) algorithm. Their analysis covers strongly convex potentials
\cite[Theorem~2.1]{roy2022stochasticLangevin} and certain nonconvex $f$ under a log-Sobolev inequality assumption for $\pi_1^f$
\cite[Theorem~3.2]{roy2022stochasticLangevin}. In the latter result,
the iteration complexity depends explicitly on the corresponding
functional-inequality constant. Such constants may have unfavorable
dependence on the dimension and inverse temperature in nonconvex
problems; see, for example, the discussion in
\cite[Section~4]{pmlr-v65-raginsky17a}.
Nakakita~\cite{nakakita2023non} develops a non-asymptotic analysis of
spherically smoothed ZO-LMC and obtains Wasserstein-$2$ guarantees under
dissipativity-type assumptions and regularity conditions on the weak
gradients. The resulting bounds involve a Poincar{\'e} constant, which may exhibit similar difficulties in nonconvex settings. 
Sahin et al.~\cite{sahin2026zerothorder} study a variance-reduced
zeroth-order Langevin method for non-log-concave sampling. Their analysis
provides non-asymptotic guarantees in squared total variation under
smoothness assumptions and a Poincar{\'e} inequality for the target
distribution. 
%Wasserstein-$2$ guarantees for zeroth-order gradient estimators in
%Markov-chain Monte Carlo methods have also been obtained under strongly
%convex assumptions; see, for example, \cite[Proposition~1]{pozza2026zeroth}.
More recently, Li and Spall~\cite{li2026zeroth} analyze zeroth-order
Langevin discretizations based on simultaneous perturbation stochastic
approximation (SPSA) for strongly convex potentials with additional
higher-order regularity assumptions.

Another line of work considers zeroth-order Langevin-based algorithms for nonconvex optimization.
Naldi et al.~\cite{naldi2026langevin} derive expected excess-risk
guarantees for ZO-LMC algorithms from
Kullback--Leibler (KL) divergence bounds. Their analysis treats the zeroth-order
estimator as an inexact approximation of $\nabla f$ and focuses on the optimization
of the original objective $f$. Under their smoothness and dissipativity
assumptions, the resulting nonconvex complexity bounds depend on a
log-Sobolev constant. A related preliminary work by Liu and Wang~\cite{liuone} establishes
non-asymptotic Wasserstein-$2$ error bounds for one-point zeroth-order
stochastic gradient Langevin dynamics under smoothness, dissipativity, and
positivity of the objective assumptions. Their bounds also depend on a
log-Sobolev constant.

A third line of work studies the geometric effect induced by the Gaussian smoothing mechanism introduced by the zeroth-order estimator. Zhang et al.~\cite{zhang2025zeroth} show that zeroth-order
optimization can favor solutions with lower Hessian trace. Their analysis
establishes an explicit connection between Gaussian smoothing and the trace
of the Hessian and provides approximate flat-minimum guarantees 
under first-, second-, and third-order smoothness assumptions together with
a convexity assumption on the original objective.
Zhang et al.~\cite{zhang2025zeroth} provide an important foundation for interpreting the smoothing
effect of zeroth-order methods as an implicit regularizer.  Our
work starts from this observation but considers a different optimization
regime. We embed the zeroth-order estimator into Langevin dynamics and we combine the
smoothing-induced curvature preference with the stochastic exploration
mechanism of Langevin dynamics. This allows us to
establish non-asymptotic distributional guarantees together with
excess-risk guarantees for the resulting curvature-regularized objective. Our analysis is based on strong convexity at infinity and
first- and third-order smoothness rather than imposing a log-Sobolev or
Poincar{\'e} assumption on the target measure. This allows us to cover the case of multimodal
objectives while yielding polynomial explicit dependence on the dimension
in the non-asymptotic bounds derived for Theorem  \ref{Theorem1_ZO}, Corollary \ref{Theorem1_ZO_W_2}, and Theorem \ref{Theorem3_ZO}.

To the best of our knowledge, there is currently no directly comparable zeroth-order Langevin result stated in Wasserstein-1 and Wasserstein-2 distances, and expected excess risk for a
curva-ture-regularized target. Therefore, we report the algorithmic complexity of  closely related methods. In this case, the  $\widetilde{\mathcal{O}}$ notation suppresses logarithmic factors and the dependence on the dimension and inverse temperature and other problem-specific parameters. For non-log-concave sampling, Roy et al.~\cite[Theorem 3.2]{roy2022stochasticLangevin} obtain a Wasserstein-2 guarantee with $\widetilde{\mathcal{O}}(\epsilon^{-4})$ oracle complexity under the batch-size choice dependent on $\epsilon$. Under weaker smoothness assumptions than those of Roy et al.~\cite{roy2022stochasticLangevin}, Nakakita~\cite[Section 3.4]{nakakita2023non} reports a Wasserstein-2 guarantee with a batch size of $\widetilde{\mathcal{O}}(\epsilon^{-4})$ and $\widetilde{\mathcal{O}}(\epsilon^{-8})$ iterations. In contrast, Corollary~\ref{Theorem1_ZO_W_2} achieves $W_2(\mathcal{L}(x_K),\pi_\beta^{v_{\eta}})\le \epsilon$ with $\widetilde{\mathcal{O}}(\epsilon^{-4})$ iterations using a single zeroth-order estimator per iteration, with the batch size fixed at one independently of $\epsilon$. Thus, our method attains the same explicit $\epsilon^{-4}$ factor in total zeroth-order oracle complexity as Roy et al.~\cite{roy2022stochasticLangevin} without requiring the batch size to increase as the target accuracy improves. Moreover, our guarantee concerns the Gibbs measure associated with the curvature-regularized objective $v_{\eta}$ rather than the original objective $f$. For optimization, Naldi et al.~\cite[Corollary~4.8]{naldi2026langevin} derive an expected excess-risk guarantee for the original objective,
$\mathbb{E}[f(x_K)]-\inf_{x\in\mathbb{R}^d}f(x)\le \epsilon$, with $\widetilde{\mathcal{O}}(\epsilon^{-4})$ total function-evaluation complexity under an accuracy-dependent choice of the zeroth-order batch size. Theorem~\ref{Theorem3_ZO} establishes a guarantee
$\mathbb{E}[v_{\eta}(x_K)]-\inf_{x\in\mathbb{R}^d} v_{\eta}(x)\le \epsilon $, with the same $\epsilon^{-4}$ factor in total function-evaluation complexity, using a batch size fixed at one independently of $\epsilon$.   Finally, Theorem~\ref{Theorem1_ZO} gives $W_1 (\mathcal{L}(x_K),\pi_\beta^{v_{\eta}})\le \epsilon$,  after $\widetilde{\mathcal{O}}(\epsilon^{-2})$ iterations. This $\epsilon$-dependence is consistent with $\epsilon^{-2}$-type Wasserstein-1 convergence rates established for stochastic first-order Langevin dynamics in nonconvex settings; see, for example, \cite[Theorem~2.4]{zhang2023nonasymptotic}.

\section{Proof Strategy} \label{sec:preliminaries_section}
In this section, we explain the main idea behind the proofs of  Theorem \ref{Theorem1_ZO}, Corollary \ref{Theorem1_ZO_W_2}, and Theorem \ref{Theorem3_ZO}. First, we define some auxiliary continuous-time processes and a continuous-time interpolation of the ZOLD recursion \eqref{eq:ZO}, which are then used for the derivation of preliminary and main results.  

\textit{Langevin diffusion associated with the smoothed objective}. We consider the process $(Y_{t})_{t \in \mathbb{R}_{+}}$, which is the solution of the following overdamped Langevin stochastic differential equation (SDE) 
\begin{equation} \label{zeroth_order_Langevin_SDE}
%\begin{split}
  \rmd Y_t  = -    \nabla f_{\eta} (Y_t ) \rmd t +  \sqrt{2 \beta^{-1}} \ \rmd B_t, \qquad  Y_0  = x_0, \
%  \end{split}
\end{equation}
where $B_t$ is a standard $d$-dimensional Brownian motion.
Denote by $(\mathcal{H}_t)_{t \ge 0}$ the natural filtration of $(B_t)_{t \ge 0}$, by $\Sigma_{x_0}$ the sigma-algebra generated by $x_0$, by $\Xi_{\infty}$ the sigma-algebra generated by $\cup_{k \in \mathbb{N}} \Xi_k$, with $(\Xi_k)_{k \in \N}$ a given filtration representing the flow of past information, and by $\Upsilon_{\infty}$ the sigma-algebra generated by $(U_{k})_{k \in \mathbb{N}}$. We consider $(\xi_{k})_{k \in \mathbb{N}}$ to be an $\mathbb{R}^m$-valued $(\Xi_k)$-adapted process with $\xi_k \sim \mathcal{L}(\xi)$ for all $k \in \mathbb{N}$. We assume that $(\mathcal{H}_t)_{t \ge 0}$ is independent of $ \Xi_{\infty}  \vee   \Sigma_{x_0} \vee \Upsilon_{\infty}$. 
%where $ \Xi_{\infty}$ is defined in Section \ref{sec:ZO_Algorithm}.
%Moreover, denote by $\mathcal{F}_{\infty} $  the sigma-algebra of $  \bigcup_{t \ge 0}  \mathcal{F}_t $. 
\begin{remark}\label{semiconvexity_assumption_remark}
As a consequence of Assumption \ref{Lipschitz_assumption_full_gradient}, one obtains
\begin{equation} \label{linear_growth_Lipschitz_assumption_remark_inequality}
    | \nabla_{x} F(x, z) | \le L_1 \varphi(z)  | x |  + L_2 \widetilde{\varphi}(z) + | \nabla F(0, 0) |,
\end{equation}
where $\widetilde{\varphi}(z) := (\varphi (z) + \varphi (0))  | z|$. Using \eqref{biased_estimate_two_point_estimator}, Assumptions \ref{independence_assumption_ZO} and \ref{Lipschitz_assumption_full_gradient},  and \cite[Lemma 1]{nesterov2017random}, one deduces
\begin{equation} \label{linear_growth_flatness_bound}
\begin{split}
  | \nabla f_{\eta}(x) | \le  L_1 \E [\varphi (\xi)] | x | +  L_1 \E [\varphi (\xi)] \eta \sqrt{d} + L_2 \E [\widetilde{\varphi}(\xi)] + |\nabla F(0,0)|.
  \end{split}
\end{equation}
 By Assumptions \ref{independence_assumption_ZO} and \ref{Lipschitz_assumption_full_gradient}, one obtains
 %   \begin{equation*}
 %   \langle \nabla_{\theta} F(\theta + \eta v,x)-  \nabla_{\theta'} F(\theta' + \eta v,x) , \theta - \theta' \rangle \geq - S \varrho(x)|\theta - \theta'|^2,
%\end{equation*}
%and in addition, by , one have 
    \begin{equation} \label{semiconvexity_flatness_bound_remark}
    | \nabla f_{\eta}(x) -  \nabla f_{\eta}(x') | \le L_1  \E[\varphi(\xi)] |x  - x'|.
\end{equation}
Using Assumptions  \ref{independence_assumption_ZO} and \ref{strong_convexity_at_infinity_assumption_ZO}, one obtains
    \begin{equation} \label{strong_convexity_flatness_bound_remark}
    \langle \nabla f_{\eta}(x) -  \nabla f_{\eta}(x') , x - x' \rangle \geq \mu \E[\rho(\xi)] |x  - x'|^2,\  \text{ if } \  |x - x'|\geq R.
\end{equation}
As a consequence of \eqref{semiconvexity_flatness_bound_remark}, \eqref{strong_convexity_flatness_bound_remark}, and \cite[Appendix B]{johnston2025the}, the Langevin SDE \eqref{zeroth_order_Langevin_SDE}   has a unique strong solution and admits a unique invariant measure which is  characterized by the density $ \frac{e^{-\beta    f_\eta(x)  }}{\int_{\mathbb{R}^d} e^{-\beta    f_\eta(x)  } \rmd x}$. 
%This follows by using .
\end{remark}

\textit{Time-rescaled dynamics of \eqref{zeroth_order_Langevin_SDE}}. For each $\lambda>0$, let $Y_{t}^{\lambda}:= Y_{\lambda t}$, $t \in \mathbb{R}_{+}$, where $(Y_{ t})_{t \in \mathbb{R}_{+}}$ is the process defined in \eqref{zeroth_order_Langevin_SDE}, and let $\widetilde{B}^{\lambda}_{t} := B_{\lambda t} / \sqrt{\lambda}$, $t \ge 0$. We note that $(\widetilde{B}^{\lambda}_{t})_{t \ge 0} $ is also a standard Brownian motion and
\begin{equation} \label{Langevin_SDE_rescaled_proof_ZO}
%\begin{split}
%Y_{0}^{\lambda} & = x_0
    \rmd Y_{t}^{\lambda}  = - \lambda    \nabla f_{\eta}(Y_{t}^{\lambda}) \ \rmd t + \sqrt{2 \lambda \beta^{-1}} \ \rmd \widetilde{B}^{\lambda}_{t}, \qquad Y_{0}^{\lambda}  = x_0.
%    \end{split}
\end{equation}
We denote the natural filtration of $(\widetilde{B}^{\lambda}_{t})_{t \ge 0}$ by $(\mathcal{H}^{\lambda}_t)_{t \ge 0}$ with $\mathcal{H}^{\lambda}_t := \mathcal{H}_{\lambda t}$, $t \in \mathbb{R}_{+}$. This filtration  is independent of $ \Xi_{\infty}  \vee   \Sigma_{x_0} \vee \Upsilon_{\infty}$.

\textit{Continuous-time interpolation of \eqref{eq:ZO}}. We define $(\bar{X}_t)_{t \in \mathbb{R}_{+}}$ as the continuous-time interpolation of the ZOLD algorithm 
\begin{equation} \label{eq:ZO_cts_interpolation}
%\begin{split}
	  \rmd \bar{X}_{t}
   =
    - \lambda  G_{\eta}(\bar{X}_{\lfloor t \rfloor }, U_{\lceil t \rceil }, \xi_{\lceil t \rceil })  \ \rmd t
    + \sqrt{2 \lambda \beta^{-1}} \rmd \widetilde{B}^{\lambda}_{t}, \qquad \bar{X}_0  = x_0.
%\end{split}
\end{equation}
To prove the main results, we work with \eqref{eq:ZO_cts_interpolation} rather than \eqref{eq:ZO}, since  the law of this process coincides with that of the ZOLD algorithm at grid-points, i.e. $\mathcal{L}(\bar{X}_k) = \mathcal{L}(x_k)$, for each $k  \in \mathbb{N}$. 

\textit{Auxiliary process for the discretization error}. We introduce a continuous-time process $(\Phi_t^{s,\chi, \lambda})_{t \ge s}$, which denotes the solution of the SDE
\begin{equation} \label{before_continuous_time_process_started_interpolation_definition_ZO}
%\begin{split}
  \rmd \Phi_t^{s,\chi, \lambda }   = - \lambda   \nabla f_{\eta}(\Phi_t^{s,\chi, \lambda }) \ \rmd t + \sqrt{2 \lambda \beta^{-1}} \  \rmd \widetilde{B}^{\lambda}_t, \qquad \Phi_s^{s,\chi, \lambda }  = \chi. 
  %  \end{split}
\end{equation}
%Let $T>0$ denote a time horizon.
    Fix $k \in \mathbb{N}$. For any $t \ge kT $, define 
    \begin{equation} \label{continuous_time_process_started_interpolation_definition_ZO}
        \bar{\Phi}_t^{\lambda,k } := \Phi_t^{kT, \bar{X}_{kT},\lambda}, \qquad  \text{where} \qquad  T:= \lfloor 1/\lambda \rfloor.
    \end{equation}
    In other words, $\bar{\Phi}_t^{\lambda,k }$ is a process started from the value of the continuous-time interpolation of the ZOLD algorithm \eqref{eq:ZO_cts_interpolation} at time $kT$ and run until time $t \ge kT$ with the Langevin dynamics.
%Intuitively, $\bar{\Phi}_t^{\lambda,k,\text{ZO} }$ is a process started from the value of the continuous-time interpolation ZO process \eqref{eq:ZO_cts_interpolation} at
%time $kT$ and run until time $t \ge kT$ with the dynamics \eqref{before_continuous_time_process_started_interpolation_definition_ZO}.

\textit{Error bounds decomposition}. Let $K \in \mathbb{N}$ and set $k = \lfloor K/T \rfloor $, then $K \in [kT, (k+1)T]$. The proofs of Theorem \ref{Theorem1_ZO}, Corollary \ref{Theorem1_ZO_W_2}, and Theorem \ref{Theorem3_ZO} rely on the following decomposition
\begin{equation} \label{decomposition_triangle_inequality_Wasserstein}
\begin{split}
       W_p(\mathcal{L}(x_K), \pi_{\beta}^{v_{\eta}} )
     &\leq W_p (\mathcal{L}(\bar{X}_{K}),\mathcal{L}(\bar{\Phi}_K^{\lambda,k}))
    + W_p(\mathcal{L}(\bar{\Phi}_K^{\lambda,k}),\mathcal{L}(Y_{K}^{\lambda}))
   \\ & \quad  + W_p (\mathcal{L}(Y_{K}^{\lambda}), \pi_{\beta}^{f_{\eta}} )  + W_p( \pi_{\beta}^{f_{\eta}}, \pi_{\beta}^{v_{\eta}}) , \quad \text{for} \  p\in\{1,2\},
    \end{split}
\end{equation}
where we have used, in the first term on the right-hand side of \eqref{decomposition_triangle_inequality_Wasserstein}, the fact that the interpolated process and the ZOLD iterates have the same law at the grid points.

\subsection{Preliminary estimates}\label{sec:preliminary_estimates}
In this section, we present the auxiliary results used in the derivation of the main results stated in Theorem \ref{Theorem1_ZO}, Corollary \ref{Theorem1_ZO_W_2}, and Theorem \ref{Theorem3_ZO}.  
First, we provide a bound for the second moment of the zeroth-order estimator and show that $f_{\eta}$ is dissipative.
\begin{lemma} \label{linear_growth_two_point_estimator_lemma}
    Let Assumptions \ref{independence_assumption_ZO} and \ref{Lipschitz_assumption_full_gradient} hold. Then, for any $x  \in \mathbb{R}^d$,
  %  \begin{equation*}
  %  \begin{split}
  %          |G_\eta(\theta, v,x)|  \le   L \varphi (x) |v|^2 |\theta|  +  L \varphi (x) \eta |v |^3 +   \bar{\varphi}(x) |v|^2 ,
  %       \end{split}
  %  \end{equation*}
%In addition, we obtain
%and
    \begin{equation*} %\label{square_linear_growth_estimator_expectation_lemma}
    \begin{split}
       \E[    |G_\eta(x, U,\xi)|^2 ]  & \le  4  (4+d) L_1^2 \E[\varphi^2 (\xi) ]  |x|^2 
        + 8  (4+d)   (L_2^2 \E[\widetilde{\varphi}^2(\xi)] + | \nabla F(0, 0) |^2)  
        \\ & \quad + 2^{-1} \eta^2 (6+d)^3 L_1^2 \E[\varphi^2 (\xi)]  .
      \end{split}
    \end{equation*}
\end{lemma}
\begin{proof}
 %   As a consequence of Assumption \ref{Lipschitz_assumption_full_gradient}, we have
%\begin{equation} \label{linear_growth_Lipschitz_assumption_remark_inequality_supplementary}
%    | \nabla_{x} F(x, z) | \le L_1 \varphi(z)  | x |  + L_2 \widetilde{\varphi}(z) + | \nabla F(0, 0) |,
%\end{equation}
%where $\widetilde{\varphi}(z) := (\varphi (z) + \varphi (0))  | z|$.
By \cite[Lemma 6.2 (b) in Section 6.1.2.1]{lan2020first} and Assumption \ref{Lipschitz_assumption_full_gradient}, we have
\begin{equation*}
    \begin{split}
         \E[    |G_\eta(x, U,\xi)|^2 ]  & \le 2 (d+4) \E [|\nabla_{x} F(x, \xi)|^2] + 2^{-1} \eta^2 L_1^2 \E [\varphi^2(\xi)] (d+6)^3.
         %\\ & \le (d+4) [ 4 L_1^2 \E [\varphi^2(\xi)] |x|^2 + 8 L_2^2 \E [\widetilde{\varphi}^2(\xi)] + 8 |\nabla F(0,0)|^2 ]
        % \\ & \quad + \frac{\eta^2}{2} L_1^2 \E [\varphi^2(\xi)] (d+6)^3.
    \end{split}
\end{equation*}
We then conclude by using \eqref{linear_growth_Lipschitz_assumption_remark_inequality} to control the first term on the right-hand side above. 
\end{proof}

\begin{remark} \label{dissipativity_assumption_remark_ZO}
Let Assumptions  \ref{independence_assumption_ZO}, \ref{Lipschitz_assumption_full_gradient}, and \ref{strong_convexity_at_infinity_assumption_ZO} hold, then, for any $x \in \mathbb{R}^d$,
\begin{equation*} 
%\label{dissipativity_remark_inequality_flatness_full_gradient_ZO}
    \langle  \nabla f_{\eta}(x) ,x \rangle\geq a |x|^2- b,
\end{equation*}
where, for $\nu \in (0, \mu \left( \E [\rho(\xi)]  - 1/2 \right)   (L_1^{2} \E [\varphi^2(\xi)])^{-1} )$,
\begin{equation}
\begin{split} \label{constants_dissipativity_flatness_full_gradient_ZO}
    a & := \mu \left( \E [\rho(\xi)]  - 1/2 \right) - \nu L_1^2 \E [\varphi^2(\xi)],  
   \\  b & := \E [  \hat{b}(\xi)] +  2 \nu (L_2^2 \E [\widetilde{\varphi}^2(\xi)] + |\nabla F(0,0)|^2)  + \eta^2 d (\nu  L_1^2 \E [\varphi^2(\xi)] + (2 \nu)^{-1} ),
    \end{split}
\end{equation}
 with $ \hat{b}(z) := \max \{ \hat{b}_1(z),  \hat{b}_2(z) \}$, for $\hat{b}_1(z):= (L_2^2 \widetilde{\varphi}^2(z) + |\nabla F(0,0)|^2 )\mu^{-1}$ and $\hat{b}_2(z):=(L_1 \varphi (z) + \mu (\rho(z) - 2^{-1}) ) R^2 +  ( L_2 \widetilde{\varphi}(z)   +| \nabla F(0,0) |)R $.
%\begin{equation} \label{expression_b_hat_proof_remark}
%\begin{split}
%  \hat{b}(z) & :=   \max \Bigg \{  (L_2^2 \widetilde{\varphi}^2(z) + |\nabla F(0,0)|^2 )\mu^{-1}, 
%  \\ & \qquad \qquad \    
%  (L_1 \varphi (z) + \mu (\rho(z) - 2^{-1}) ) R^2 +  ( L_2 \widetilde{\varphi}(z)   +| \nabla F(0,0) |)R  \Bigg \}. 
%  \end{split}
%\end{equation}
%Therefore, all the critical points of $x \mapsto f_{\eta}(x)$ are contained in the ball of radius $\sqrt{b/a}$ centered at the origin.
\end{remark}
\begin{proof}
%By following the proof of \cite[Lemma 2]{johnston2025the}, we can derive the  dissipativity of $F(x,z)$. In other words,    %We follow the proof of \cite[Lemma 2]{johnston2025the} and adapt it to accommodate our assumptions. 
Using Assumption \ref{strong_convexity_at_infinity_assumption_ZO}, Young's inequality, and \eqref{linear_growth_Lipschitz_assumption_remark_inequality} in the proof of \cite[Lemma 2]{johnston2025the},  one obtains, for $| x | \ge R$,
\begin{equation} \label{dissipativity_remark_inequality_outside_ball_proof_ZO}
  %  \begin{split}
         \langle \nabla_{x} F(x,z) ,x \rangle 
         %& =  \langle \nabla_{x} F(x,z) - \nabla_{x^{'}} F(0,z)  , x   \rangle + \langle \nabla_{x^{'}} F(0,z)  ,x \rangle  
         %\\ & \ge \mu \rho(x) |\theta |^2 -  | \nabla_{\theta^{'}} F(0,x) | | \theta |
        %    \\ & \ge \mu (\rho(z) - 1/2) |x |^2 -  | \nabla_{x^{'}} F(0,z) |^2/(2 \mu)
             \ge \mu (\rho(z) - 1/2) |x |^2  - (L_2^2 \widetilde{\varphi}^2(z)+ |\nabla F(0,0)|^2)/\mu.  
   % \end{split}
\end{equation}
Using the Cauchy--Schwarz inequality and \eqref{linear_growth_Lipschitz_assumption_remark_inequality} in the proof of \cite[Lemma 2]{johnston2025the}, one derives, for $|x | < R$,
\begin{equation} \label{dissipativity_remark_inequality_inside_ball_proof_ZO}
    \begin{split}
         \langle \nabla_{x} F(x,z) ,x \rangle    
         %\ge  - | \nabla_{x} F(x,z) | | x |  
       %  \\ & \ge -L_1 \varphi (z) |x |^2 - ( L_2 \widetilde{\varphi}(z) + | \nabla F(0,0) |) |x |  
    %     \\ & \ge \mu (\rho(z) - 1/2) R^2 -L_1 \varphi (z) R^2
     %    \\ & \quad -  (L_2 \widetilde{\varphi}(z) + | \nabla F(0,0) |) R - \mu (\rho(z) - 1/2) R^2
             & \ge \mu (\rho(z) - 1/2) |x |^2  - (L_1 \varphi (z) + \mu (\rho(z) - 1/2) ) R^2 
            \\ & \quad -  ( L_2 \widetilde{\varphi}(z)   +| \nabla F(0,0) |) R.  
    \end{split}
\end{equation}
Combining \eqref{dissipativity_remark_inequality_outside_ball_proof_ZO} and \eqref{dissipativity_remark_inequality_inside_ball_proof_ZO} leads to
%the  dissipativity of $F(x,z)$, which is
\begin{equation} \label{dissipativity_remark_inequality_supplementary}
       \langle \nabla_{x} F(x,z) ,x \rangle \ge  \mu (\rho(z) - 1/2) |x |^2  - \hat{b}(z),
\end{equation}
%\eqref{dissipativity_remark_inequality}
%with $ \hat{b}(z) $ given in Remark \eqref{expression_b_hat_proof_remark}.
where $ \hat{b}(z)$ is given in the statement of Remark \ref{dissipativity_assumption_remark_ZO}. Using \eqref{dissipativity_remark_inequality_supplementary}, Young's inequality, and  \eqref{linear_growth_Lipschitz_assumption_remark_inequality}, one obtains, for $u \in \mathbb{R}^d$ and  $\nu \in (0, \mu \left( \E [\rho(\xi)]  - 1/2 \right)   (L_1^{2} \E [\varphi^2(\xi)])^{-1} )$
\begin{equation} \label{dissipativity_remark_proof_flatness_stochastic_gradient_ZO}
\begin{split}
  &  \langle \nabla_{x} F(x + \eta u ,z) , x \rangle 
  \\ & =     \langle \nabla_{x} F(x + \eta u ,z) , x + \eta u \rangle -  \langle \nabla_{x} F(x + \eta u ,z) , \eta u \rangle 
    \\ & \geq \mu \left(\rho(z)  - 1/2 \right)|x + \eta u  |^2-\hat{b}(z)
   - 2^{-1} \nu | \nabla_{x} F(x + \eta u ,z)|^2 - (2 \nu)^{-1} \eta^2 |  u|^2
       \\ & \geq \mu \left(\rho(z)  - 1/2 \right) (|x |^2 + 2 \langle x , \eta u \rangle  ) -\hat{b}(z) \\ & \quad -  \nu ( L_1^2 \varphi^2(z) (| x |^2 + 2 \langle x, \eta u \rangle  + \eta^2 | u |^2)) 
       \\ & \quad -   2 \nu  L_2^2 \widetilde{\varphi}^2(z) - 2 \nu  |\nabla F(0,0)|^2  - (2 \nu)^{-1} \eta^2 |  u|^2.
    \end{split}
\end{equation}
Using \eqref{dissipativity_remark_proof_flatness_stochastic_gradient_ZO} and \cite[Lemma 1]{nesterov2017random},  we obtain
%Taking the expectation over $X$ and $U \sim \mathcal{N}(0, I_d)$ on both sides of \eqref{dissipativity_remark_proof_flatness_stochastic_gradient_ZO}, using \cite[Lemma 1]{nesterov2017random}, we obtain
\begin{equation*}
\begin{split}
     \langle  \nabla f_{\eta}(x) ,x \rangle & =  \langle \E [\nabla_{x} F(x + \eta U ,\xi)] , x \rangle 
     \\ &  \ge (\mu \left( \E [\rho(\xi)]  - 1/2 \right)  - \nu L_1^2 \E [\varphi^2(\xi)]  ) |x |^2 - \E [\hat{b}(\xi)] 
     \\ & \quad -  2 \nu (L_2^2 \E [\widetilde{\varphi}^2(\xi)] + |\nabla F(0,0)|^2)
       - \eta^2 d (\nu  L_1^2 \E [\varphi^2(\xi)] + (2 \nu)^{-1} ).
     \end{split}
\end{equation*}
%Along the same lines of the proof of \cite[Lemma D.1]{bruno2026flatnessaware}, the ball 
%\begin{equation} 
%\label{definition_set_A_minimizers_ZO}
%\mathfrak{D} = \left\{ x \in \mathbb{R}^d : |x | \le \sqrt{b/a} \right\},
%\end{equation}
%contains all the critical points of $f_{\eta}(x)$.
\end{proof}
%The proof of Remark \ref{dissipativity_assumption_remark_ZO}, including the dissipativity condition on $f$,  is deferred to the supplementary material. 
%The following result is a consequence of the dissipativity condition of $f_{\eta}$. 
%    \begin{lemma} \label{set_A_remark_dissipativity_ZO}
%Let Assumptions \ref{Lipschitz_assumption_full_gradient} and \ref{strong_convexity_at_infinity_assumption_ZO} hold. Then, the following set 
%\begin{equation} 
%\label{definition_set_A_minimizers_ZO}
%\mathfrak{B} := \left\{ \theta \in \mathbb{R}^d : |\theta | \le \sqrt{b/a} \right\},
%\end{equation}
%contains all the minimizers of $f_{\eta}(\theta)$, where $a$ and $b$ are given in \eqref{constants_dissipativity_flatness_full_gradient_ZO}.
%\end{lemma}
%using the dissipativity 
%condition for $\nabla f_{\eta}$ established in Remark \ref{dissipativity_assumption_remark_ZO}.
We provide a second-moment bound for the continuous-time interpolation process  $(\bar{X}_{t})_{t \ge 0}$ of the ZOLD algorithm defined in \eqref{eq:ZO_cts_interpolation}.
\begin{lemma}\label{Lemma2_ZO}
    Let Assumptions \ref{independence_assumption_ZO}, \ref{Lipschitz_assumption_full_gradient}, and \ref{strong_convexity_at_infinity_assumption_ZO} hold. Then, for any $\lambda\in(0,\lambda_{\text{max}})$ with $\lambda_{\text{max}}$ given in \eqref{expression_lambda_max_ZO}, 
    \begin{equation*} %\label{eqProp2}
\sup_{t\geq0}\mathbb{E}[|\bar{X}_{t} |^2 ] \leq C_1 (1+\mathbb{E}[| x_0|^2] ),
    \end{equation*}
    where $C_1 = \mathcal{O} ( 1+ \eta^2 d^2 + d / \beta ) $ and is given explicitly in \eqref{eqProp2_constants}.
    %, whose explicit expression is given in %\eqref{eqProp2_constants}. %and 
%    \begin{equation} \label{eqProp2_constants}
%    \begin{split}
%      C_3 & =   \left(\frac{8 a^2}{ L^2 \E [\varphi^2(X)] (4+d)^2} + 2 \right)  C_2 \\ & \quad + \frac{8a}{(4+d)^2 L^2 \E[\varphi^2(X)]} \left(   \frac{2 \eta^2 a (6+d)^3}{(4+d)^2 }  + \frac{2 a  \E [\bar{\varphi}^2(X) ]) }{ L^2 \E[\varphi^2(X)]}+ \frac{d}{\beta} \right).
%          \end{split}
%    \end{equation}
%    with the constants $a$ and $b$ given in \eqref{constants_dissipativity_flatness_full_gradient_ZO}, and
%     \begin{equation} \label{eqProp2_constants_C_2}
%    \begin{split}
%         C_2 & =\left( b+  2 a \eta^2 \frac{(6+d)^3}{(4+d)^2}    + 2 a \frac{   \E[\bar{\varphi}^2(X)]}{L^2 \E[\varphi^2 (X) ]   }  + \frac{d}{\beta} \right) \frac{1}{ a -\lambda (4+d)^2 L^2 \E[\varphi^2(X)]} + 1.
%          \end{split}
%    \end{equation}
%    In particular, $C_2= \mathcal{O} \left( 1+ \eta^2 d + \frac{d}{\beta} \right)$ and $C_3 = \mathcal{O} \left( 1+ \eta^2 d + \frac{d}{\beta} \right) $.
    \end{lemma}
    \begin{proof}
     %   We begin by considering the ZO iterates $(\theta_{k}^{\text{ZO}})_{k \geq 0}$ in \eqref{eq:ZO} corresponding to interpolation scheme $(\bar{\theta}_t^{\text{ZO}})_{t \in \mathbb{R}_{+}}$ defined in \eqref{eq:ZO_cts_interpolation},
     Using the ZOLD recursion \eqref{eq:ZO}, one obtains
\begin{equation} \label{iterate_ZO_second_moment_new}
\begin{split}
|x_{k+1}|^2 
%& =\left|\theta_{k}^{\text{ZO}}-\lambda G^{(2)}_{\eta}(\theta^{\text{ZO}}_k, V_{k+1}, X_{k+1}) \right|^2+\frac{2\lambda}{\beta}|\xi_{k+1}|^2 \\ & \quad +2\langle \theta_{k}^{\text{ZO}} -\lambda G^{(2)}_{\eta}(\theta^{\text{ZO}}_k, V_{k+1}, X_{k+1}) , \sqrt{2 \lambda \beta^{-1}}\xi_{k+1}\rangle
 & = |x_{k} |^2 - 2 \lambda \langle x_{k} , G_{\eta}(x_k, U_{k+1}, \xi_{k+1}) \rangle 
  + \lambda^2 | G_{\eta}(x_k, U_{k+1}, \xi_{k+1}) |^2 
 \\ & \quad +\frac{2\lambda}{\beta}|\zeta_{k+1}|^2 
  +2\langle x_{k} -\lambda G_{\eta}(x_k, U_{k+1}, \xi_{k+1}) , \sqrt{2 \lambda \beta^{-1}} \zeta_{k+1}\rangle .
\end{split}
\end{equation}
%Since $\theta_{k}^{\text{ZO}}$ is independent of  $\xi_{k+1}$, $V_{k+1}$, and $X_{k+1}$, the last term vanishes under expectation. Thus 
Taking the conditional expectation on $x_{k}$  on both sides of \eqref{iterate_ZO_second_moment_new} and using Assumption \ref{independence_assumption_ZO},  \eqref{biased_estimate_two_point_estimator}, Remark \ref{dissipativity_assumption_remark_ZO}, Lemma \ref{linear_growth_two_point_estimator_lemma}, $\xi_k \sim \mathcal{L}(\xi)$, and $\lambda< a/(4(4+d) L_1^2 \E[\varphi^2(\xi)])$, we obtain
\begin{equation*}
\begin{split}
 \mathbb{E} \left[|x_{k+1}|^2  | x_{k} \right]   &= |x_{k}|^2 -2\lambda \mathbb{E}\left[\langle x_{k}, G_{\eta}(x_k, U_{k+1}, \xi_{k+1}) \rangle | x_{k} \right] 
\\ & \quad +\lambda^2  \mathbb{E}\left[|G_{\eta}(x_k, U_{k+1}, \xi_{k+1}) |^2 | x_{k}  \right] +2\lambda d/\beta
\\
    &\leq (1 - 2 \lambda a  + 4 \lambda^2   (4+d) L_1^2 \E[\varphi^2 (\xi) ]   ) | x_{k} |^2 + 2 \lambda b +2\lambda d/\beta 
    \\ & \quad  + \lambda^2  (6+d)^3 \frac{\eta^2}{2} L_1^2 \E[\varphi^2 (\xi)]   
     +  \lambda^2  8 (4+d) (L_2^2 \E[\widetilde{\varphi}^2(\xi)] + |\nabla F(0,0)|^2)   
    \\
    & \le  (1 -  \lambda a   ) | x_{k} |^2 
    \\ & \quad +  \lambda \left(  2 b+     \frac{a \eta^2(6+d)^3}{8(4+d)}    +  \frac{ 2 a(  L_2^2 \E[\widetilde{\varphi}^2(\xi)] + |\nabla F(0,0)|^2)}{L_1^2 \E[\varphi^2 (\xi) ]   }  + \frac{2 d}{\beta}   \right).
    \end{split}
\end{equation*}
This implies that
%Taking the expectation with respect to the iterate $x_{k}$ on both sides above, we arrive at
\begin{equation} \label{ULA_bound}
\begin{split}
 \mathbb{E}\left[| x_{k+1}|^2\right]
 &\leq \left(1- \lambda a \right)^{k+1} \mathbb{E}[|x_0|^2] 
\\ & \quad +   \frac{2 b}{a}+   \frac{ \eta^2(216+d^3)}{2 d}    +  \frac{2 (   L_2^2 \E[\widetilde{\varphi}^2(\xi)] + |\nabla F(0,0)|^2)}{L_1^2 \E[\varphi^2 (\xi) ]   }  + \frac{2 d}{a \beta}  
\\
    &\leq C_0\left(1+\mathbb{E}[|x_0|^2]\right),
    \end{split}
\end{equation}
where 
  \begin{equation} \label{eqProp2_constants_C_2}
    \begin{split}
         C_0 & = \frac{2 b}{a}+    \frac{108 \eta^2}{ d} + \frac{\eta^2 d^2}{2 }    +   \frac{   2 L_2^2 \E[\widetilde{\varphi}^2(\xi)] +2  |\nabla F(0,0)|^2}{L_1^2 \E[\varphi^2 (\xi) ]   }  + \frac{2 d}{a \beta}  + 1.
          \end{split}
    \end{equation}
%$C_2$ is given in \eqref{eqProp2_constants_C_2}. 
By \eqref{constants_dissipativity_flatness_full_gradient_ZO}, we have $b = \mathcal{O}(1+ \eta^2 d)$ and $a=\mathcal{O}(1)$.  Therefore,  $C_0= \mathcal{O} ( 1+ \eta^2 d^2 + d / \beta )$. Let $t\in[k,k+1)$ with $k  \in \mathbb{N}$. Using  the continuous-time interpolation of the ZOLD algorithm defined in \eqref{eq:ZO_cts_interpolation} and  H\"older's inequality,  one obtains
\begin{equation} \label{second_moment_interpolation_scheme}
    \begin{split}
 |\bar{X}_{t} |^2 
 &  \le  2|\bar{X}_{t} - \bar{X}_{\lfloor t\rfloor} |^2+2|\bar{X}_{\lfloor t\rfloor} |^2  
%\\ & \leq 4  \lambda^2 \left|\int_{\lfloor t \rfloor}^t G_{\eta}(\bar{X}_{\lfloor s \rfloor }, U_{\lceil s \rceil }, \xi_{\lceil s \rceil }) \ \rmd s \right|^2 +\frac{8\lambda}{\beta}|\widetilde{B}_t^{\lambda}-\widetilde{B}_{\lfloor t\rfloor}^{\lambda}|^2+2|\bar{X}_{\lfloor t\rfloor}|^2
\\ &\leq 4\lambda^2(t-\lfloor t\rfloor)\int_{\lfloor t\rfloor}^t|G_{\eta}(\bar{X}_{\lfloor s \rfloor }, U_{\lceil s \rceil }, \xi_{\lceil s \rceil })|^2 \rmd s 
    \\ & \quad + \frac{8\lambda}{\beta}|\widetilde{B}^{\lambda}_t-\widetilde{B}^{\lambda}_{\lfloor t\rfloor}|^2+2|\bar{X}_{\lfloor t\rfloor}|^2.
%    \\ & \leq 4\lambda^2\int_{\lfloor t\rfloor}^t (2 L_1^2 \varphi^2 (X_{\lceil s \rceil }) |U_{\lceil s \rceil }|^4 | \bar{\theta}_{\lfloor s \rfloor }|^2  +  4 L_1^2 \varphi^2 (X_{\lceil s \rceil }) \eta^2 |U_{\lceil s \rceil } |^6  ) \  \rmd s 
%    \\ & \quad + 4\lambda^2\int_{\lfloor t\rfloor}^t  8  |U_{\lceil s \rceil }|^4 ( L_2^2 \widetilde{\varphi}^2(X_{\lceil s \rceil }) + |\nabla F(0,0)|^2)  \  \rmd s  + \frac{8\lambda}{\beta}|\widetilde{B}^{\lambda}_t-\widetilde{B}^{\lambda}_{\lfloor t\rfloor}|^2+2|\bar{\theta}_{\lfloor t\rfloor}|^2.
        \end{split}
\end{equation}
%We observe that for any $s\in\left[\lfloor t\rfloor,t\right]$, we have $\lfloor s\rfloor=\lfloor t \rfloor$. Therefore, 
Taking expectations on both sides of \eqref{second_moment_interpolation_scheme}, using Lemma \ref{linear_growth_two_point_estimator_lemma}, $\xi_{\lceil t \rceil} \sim \mathcal{L}(\xi)$, and $\lambda< a/(4(4+d) L_1^2 \E[\varphi^2(\xi)])$, leads to
\begin{equation*}
     \begin{split}
& \mathbb{E} [|\bar{X}_{t} |^2]
\\ & \le  (16 \lambda^2 L_1^2 \E [\varphi^2(\xi)] (4+d) +2 )\mathbb{E} [|\bar{X}_{\lfloor t\rfloor}|^2] 
 \\ & \quad +  \lambda^2 ( 2 (6+d)^3 \eta^2  L_1^2 \E [\varphi^2(\xi)]  +   32 (4+d) (L_2^2 \E [\widetilde{\varphi}^2(\xi) ] + |\nabla F(0,0)|^2))  + \frac{8\lambda d}{\beta}  
\\ & \le \left( \frac{ a^2}{ L_1^2 \E [\varphi^2(\xi)] (4+d)} + 2 \right) \mathbb{E} [|\bar{X}_{\lfloor t\rfloor}|^2] 
+ \frac{a^2(216 + d^3) \eta^2}{2 (16+ d^2) L_1^2 \E [\varphi^2(\xi)]  }
\\ & \quad + \frac{2 a^2 (L_2^2 \E [\widetilde{\varphi}^2(\xi) ] + |\nabla F(0,0)|^2)}{(4+d) L_1^4 (\E [\varphi^2(\xi)])^2 } + \frac{2 a d }{\beta (4+d)   L_1^2 \E [\varphi^2(\xi)]  }.
%\\ & \quad +   \frac{a}{2(4+d) L_1^2 \E[\varphi^2(X)]} \left(   \frac{ \eta^2 a (6+d)^3}{(4+d) }  + \frac{16 a (L_2^2 \E [\widetilde{\varphi}^2(X) ] + |\nabla F(0,0)|^2) }{ L_1^2 \E[\varphi^2(X)]}+ \frac{8 d}{\beta} \right). 
%\\ & \textcolor{red}{OLD}\leq  \frac{16 a^2 \eta^2 (d+4)}{d(d+2)L^2 \E[\varphi^2(X_0)]} + \frac{16 a^2(2  \E [\bar{\varphi}^2(X_{0})] }{d(d+2) L^4 (\E [\varphi^2(X_{0 })])^2 }  + \frac{8 a}{\beta (d+2) L^2 \E [\varphi^2(X_0)]}
%\\ & \quad +\left(  \frac{8 a^2}{d (d+2) L^2 \E[\varphi^2(X_0)]}+2\right)\mathbb{E} [|\bar{\theta}_{\lfloor t\rfloor}^{\mathrm{ZO}}|^2].
\end{split}
\end{equation*}
Recall that  the law of the interpolation scheme \eqref{eq:ZO_cts_interpolation} coincides with the law
of the ZOLD algorithm at grid-points, i.e., $\mathcal{L}(\bar{X}_{ \lfloor t \rfloor }) = \mathcal{L}(x_k)$ for $t\in[k,k+1)$. Therefore, using the bound \eqref{ULA_bound} on the right-hand side above, we arrive at
\begin{equation*}
    \mathbb{E} [|\bar{X}_{t} |^2 ] \leq C_1\left(1+\mathbb{E}[|x_0|^2] \right),
\end{equation*}
where 
\begin{equation} \label{eqProp2_constants}
        \begin{split}
                  C_1 & =    \left(  \frac{2 b}{a}+    \frac{108 \eta^2}{ d} + \frac{\eta^2 d^2}{2 }    +   \frac{   2 L_2^2 \E[\widetilde{\varphi}^2(\xi)] +2  |\nabla F(0,0)|^2}{L_1^2 \E[\varphi^2 (\xi) ]   }  + \frac{2 d}{a \beta}  + 1 \right) 
         \\ & \qquad \times \left(\frac{ a^2}{d  L_1^2 \E [\varphi^2(\xi)] } + 2  \right)  
           + \frac{a^2(216 + d^3) \eta^2}{2 d^2 L_1^2 \E [\varphi^2(\xi)]  }
\\ & \quad + \frac{2 a^2 (L_2^2 \E [\widetilde{\varphi}^2(\xi) ] + |\nabla F(0,0)|^2)}{d L_1^4 (\E [\varphi^2(\xi)])^2 } + \frac{2 a  }{\beta    L_1^2 \E [\varphi^2(\xi)]  },
        \end{split}
    \end{equation}
    with $a$ and $b$ defined in \eqref{constants_dissipativity_flatness_full_gradient_ZO}. Thus, $ \sup_{t \ge 0} \mathbb{E} [|\bar{X}_{t} |^2 ] \leq C_1\left(1+\mathbb{E}[|x_0|^2] \right)$.
      \end{proof}
%$C_1$ given in \eqref{eqProp2_constants}.  In particular, $C_1 = \mathcal{O} \left( 1+ \eta^2 d^2 + \frac{d}{\beta} \right) $.
%$C_3$ is given in \eqref{eqProp2_constants}.
%    \begin{equation} \label{eqProp2_constants}
%    \begin{split}
%      C_1 & =    \Bigg( \left( b+  \frac{ 4 a \eta^2(6+d)^3}{(4+d)^2}    +  \frac{8 a (  L_2^2 \E[\widetilde{\varphi}^2(X)]+ |\nabla F(0,0)|^2)}{L_1^2 \E[\varphi^2 (X) ]   }  + \frac{d}{\beta} \right)
 %        \\ & \qquad \qquad  \times \left( a -\lambda (4+d)^2 L_1^2 \E[\varphi^2(X)]\right)^{-1} + 1 \Bigg) \left(\frac{8 a^2}{ L_1^2 \E [\varphi^2(X)] (4+d)^2} + 2 \right)  \\ & \quad + \frac{8a}{(4+d)^2 L_1^2 \E[\varphi^2(X)]} \left(   \frac{2 \eta^2 a (6+d)^3}{(4+d)^2 }  + \frac{4 a  (L_2 \E [\widetilde{\varphi}^2(X) ]) + |\nabla F(0,0)|^2) }{ L_1^2 \E[\varphi^2(X)]}+ \frac{d}{\beta} \right),
 %         \end{split}
%    \end{equation}
%    with the constant $C_0$ given explicitly in \eqref{eqProp2_constants_C_2}.
    %and the constant $a$ given in \eqref{constants_dissipativity_flatness_full_gradient_ZO}, and .
Then, we derive an estimate for the one-step error associated with  $(\bar{X}_{t})_{t \ge 0}$.  \begin{lemma}\label{one_step_error_lemma_ZO} Let Assumptions \ref{independence_assumption_ZO},  \ref{Lipschitz_assumption_full_gradient}, and \ref{strong_convexity_at_infinity_assumption_ZO} hold. Then, for any $t >0$ and   $\lambda\in(0,\lambda_{\text{max}})$  with $\lambda_{\text{max}}$ given in \eqref{expression_lambda_max_ZO},
\begin{equation*}
\mathbb{E} [|\bar{X}_t - \bar{X}_{\lfloor t \rfloor}|^2 ] \leq  C_2 \lambda(  1+   \mathbb{E}[|x_0|^2]  ),
\end{equation*}
where  $C_2=  \mathcal{O} ( 1+ \eta^2 d^2 + d / \beta )$ and is given explicitly in \eqref{constants_one_step_error_lemma_ZO}.
%\begin{equation} \label{constants_one_step_error_lemma_ZO}
%\begin{split}
 %   C_{2} & =2a C_1 +   \frac{ a \eta^2 (216+d^3)}{d} + \frac{4a (L_2^2  \E[\widetilde{\varphi}^2(X)] + |\nabla F(0,0)|^2) }{L_1^2 \E [\varphi^2 (X)]}  + \frac{4 d}{\beta},
%    \end{split}
%\end{equation}
%with $C_1$ given in \eqref{eqProp2_constants}.
%,whose explicit expression is given in \eqref{constants_one_step_error_lemma_ZO}.
%\begin{equation} \label{constants_one_step_error_lemma_ZO}
%\begin{split}
%    C_{2} & =4 a C_1+ 8 a \eta^2 (d+4) + \frac{16a   \E[\bar{\varphi}^2(X)]  }{L^2 \E [\varphi^2 (X)]} +8  d/\beta ,
%    \end{split}
%\end{equation}
%where $C_1$ is  defined in Lemma \ref{Lemma2_ZO}.
\end{lemma}
    \begin{proof}
%One considers the difference between $\bar{\theta}_{\lfloor t \rfloor}^{\text{fSGLD}}, \bar{\theta}_t^{\text{fSGLD}}$ to get the one-step error
Using the definition of the continuous-time interpolation of the ZOLD algorithm  \eqref{eq:ZO_cts_interpolation},  we obtain 
\begin{equation*}
%\begin{split}
  \E \left[  |\bar{X}_t - \bar{X}_{\lfloor t \rfloor}|^2 \right] = \E \left[ \left|\int_{\lfloor t \rfloor}^t - \lambda G_{\eta}(\bar{X}_{\lfloor s \rfloor }, U_{\lceil s \rceil }, \xi_{\lceil s \rceil })  \ \rmd s  +\sqrt{\frac{2\lambda}{\beta} } (\widetilde{B}_t^{\lambda}  - \widetilde{B}_{\lfloor t \rfloor}^{\lambda}) \right|^2 \right].
  %  \\ & \le 2 \lambda^2 \int_{\lfloor t \rfloor}^t  (2 L_1^2 \varphi^2 (X_{\lceil s \rceil }) |U_{\lceil s \rceil }|^4 |\bar{\theta}_{\lfloor s \rfloor }|^2  +  4 L_1^2 \varphi^2 (X_{\lceil s \rceil }) \eta^2 |U_{\lceil s \rceil } |^6  ) \ \rmd s 
%    \\ & \quad +  2 \lambda^2 \int_{\lfloor t \rfloor}^t    8  |U_{\lceil s \rceil }|^4 (L_2^2 \widetilde{\varphi}^2(X_{\lceil s \rceil }  ) + |\nabla F(0,0)|^2)   \ \rmd s   + \frac{4\lambda}{\beta}\left|\widetilde{B}_{\lfloor t \rfloor}^{\lambda}-\widetilde{B}_t^{\lambda} \right|^2.
 %   \end{split}
\end{equation*}
Using Assumption \ref{independence_assumption_ZO}, Lemma \ref{linear_growth_two_point_estimator_lemma}, Lemma \ref{Lemma2_ZO}, and $\lambda< a/(4(4+d) L_1^2 \E[\varphi^2(\xi)])$ leads to
\begin{equation*}
\begin{split}
 &    \mathbb{E}|\bar{X}_t - \bar{X}_{\lfloor t \rfloor}|^2 
 \\ & \leq 2\lambda^2 \int_{\lfloor t \rfloor}^t  4 L_1^2 \E [\varphi^2(\xi_{\lceil s \rceil })]    (4+d) C_1\left(1+\mathbb{E}[| x_0|^2] \right)  \rmd s
 \\ & \quad  + \lambda^2 \int_{\lfloor t \rfloor}^t   L_1^2 \E [\varphi^2 (\xi_{\lceil s \rceil })] \eta^2 (6+d)^3 \ \rmd s
 \\ &  \quad + 2\lambda^2 \int_{\lfloor t \rfloor}^t 8 (4+d)(L_2^2 \E[\widetilde{\varphi}^2(\xi_{\lceil s \rceil })] + |\nabla F(0,0)|^2   ) \ \rmd s
  + \frac{4\lambda d}{\beta} 
    \\ & \leq  \lambda \left( 2a C_1\left(1+\mathbb{E}[| x_0|^2] \right) + \frac{ a \eta^2 (216+d^3)}{ 4+d} + \frac{4 a (L_2^2  \E[\widetilde{\varphi}^2(\xi)] + |\nabla F(0,0)|^2)  }{L_1^2 \E [\varphi^2 (\xi)]}  + \frac{4 d}{\beta} \right) 
     \\ & \leq   \lambda C_2 \left(1+\mathbb{E}[|x_0|^2]\right), 
    \end{split}
\end{equation*}
where 
\begin{equation} \label{constants_one_step_error_lemma_ZO}
\begin{split}
    C_{2} & =2a C_1 +   \frac{ a \eta^2 (216+d^3)}{d} + \frac{4a (L_2^2  \E[\widetilde{\varphi}^2(\xi)] + |\nabla F(0,0)|^2) }{L_1^2 \E [\varphi^2 (\xi)]}  + \frac{4 d}{\beta},
    \end{split}
\end{equation}
with $a$ and $C_1$ defined in \eqref{constants_dissipativity_flatness_full_gradient_ZO} and \eqref{eqProp2_constants}, respectively.
%\begin{equation} \label{constants_one_step_error_lemma_ZO}
%\begin{split}
%    C_{2} & =4 a C_1 + \frac{8 a \eta^2 (6+d)^3}{(4+d)^2} + \frac{8a (L_2^2  \E[\widetilde{\varphi}^2(X)] + |\nabla F(0,0)|^2) }{L_1^2 \E [\varphi^2 (X)]}  + \frac{4 d}{\beta},
%    \end{split}
%\end{equation}
%with $C_1$ given in \eqref{}. In particular, $C_2 =\mathcal{O} \left( 1+ \eta^2 d + \frac{d}{\beta} \right)  $.
\end{proof}
The following result provides an upper bound for the second moment of the auxiliary process defined in  \eqref{continuous_time_process_started_interpolation_definition_ZO}.
\begin{lemma}\label{uniform_bound_auxiliary_lemma_ZO}  
Let Assumptions \ref{independence_assumption_ZO},  \ref{Lipschitz_assumption_full_gradient}, and \ref{strong_convexity_at_infinity_assumption_ZO} hold. Then, for any $\lambda\in(0,\lambda_{\text{max}})$  with $\lambda_{\text{max}}$ given in \eqref{expression_lambda_max_ZO}, 
\begin{equation*}
    \sup_{kT\leq t\leq (k+1)T}\mathbb{E} [| \bar{\Phi}_t^{\lambda,k }|^2] \leq  ( C_1 +2d/\beta+ 2b) \left(1+\mathbb{E} [|x_0|^2] \right),
\end{equation*}
with $C_1$ and $b$ defined in \eqref{eqProp2_constants} and  \eqref{constants_dissipativity_flatness_full_gradient_ZO}, respectively.
\end{lemma}
%where
%\begin{equation}
%    C_3 =  C_1 +2(d/\beta+b),
%\end{equation}
%, and whose explicit expression is provided in the proof.
    \begin{proof}
This follows along the same lines as the proof of  \cite[Lemma 6]{johnston2025the}, using the dissipativity of $f_{\eta}$ (Remark \ref{dissipativity_assumption_remark_ZO}) and Lemma \ref{Lemma2_ZO}.
    \end{proof}
Next, we obtain an upper bound for the first term on the right-hand side of \eqref{decomposition_triangle_inequality_Wasserstein}.    
\begin{lemma} \label{first_upper_bound_lemma}
    Let Assumptions \ref{independence_assumption_ZO}, \ref{Lipschitz_assumption_full_gradient}, and \ref{strong_convexity_at_infinity_assumption_ZO} hold. Then, for any $\lambda\in(0,\lambda_{\text{max}})$, with $\lambda_{\text{max}}$ given in \eqref{expression_lambda_max_ZO}, $t \in [kT, (k+1)T]$ and $k \in \N$,
    \begin{equation*}
        W_2 (\mathcal{L}(\bar{X}_t), \mathcal{L}(\bar{\Phi}_t^{\lambda,k } )) \le  \sqrt{\lambda} C_{3},
    \end{equation*}
    where $ C_{3} = \mathcal{O}(1+ d \beta^{- 1/2}+ \eta d^{3/2} + d^{1/2})$ and is given explicitly in \eqref{constant_first_inequality_split}.
\end{lemma}
\begin{proof}
By construction, the processes $(\bar{X}_t)_{t \ge kT}$ and $\bar{\Phi}_t^{\lambda,k }$ defined in \eqref{eq:ZO_cts_interpolation} and in  \eqref{continuous_time_process_started_interpolation_definition_ZO}, respectively, have the same initial condition $\bar{X}_{kT}$ and are driven by the same Brownian motion. We denote by $\mathcal{F}_{l,s }$, for $ s \in [l, l+1)$, the sigma-algebra generated by $(x_0, (U_{j})_{1 \le j \le  l }, (\xi_{j})_{1 \le j \le  l }, (\widetilde{B}^{\lambda}_{r})_{ 0 \le r \le s })$.   % and we recall that $(\mathcal{H}^{\lambda}_t)_{t \ge 0}$ is  the natural filtration of $(\widetilde{B}^{\lambda}_{t})_{t \ge 0}$.
 % $(\bar{X}_{ \lfloor  s \rfloor })_{kT \le   s < t}$ by, 
%In addition, we note that, for any $s \in [\lfloor t \rfloor , t]$, we have $\lfloor s \rfloor = \lfloor t \rfloor$ and $\lceil s \rceil = \lceil t \rceil$, and thus
%\begin{equation} \label{one_step_grid_ZOLD}
%\begin{split}
%    \bar{\theta}_{t} -  \bar{\theta}_{\lfloor t \rfloor } & = - \lambda \int_{\lfloor t \rfloor }^t   G_{\eta}(\bar{\theta}_{\lfloor s \rfloor }, U_{\lceil s \rceil }, X_{\lceil s \rceil }) \ \rmd s
%    + \sqrt{2 \lambda \beta^{-1}}  (\widetilde{B}^{\lambda}_{t} - \widetilde{B}^{\lambda}_{\lfloor t \rfloor})
%    \\ & = - \lambda (t - \lfloor t \rfloor )   G_{\eta}(\bar{\theta}_{\lfloor t \rfloor }, U_{\lceil t \rceil }, X_{\lceil t \rceil }) 
%    + \sqrt{2 \lambda \beta^{-1}}  (\widetilde{B}^{\lambda}_{t} - \widetilde{B}^{\lambda}_{\lfloor t \rfloor}).
 %   \end{split}
%\end{equation}
By using  It{\^o}'s formula, the definition of the   continuous-time interpolation of ZOLD algorithm \eqref{eq:ZO_cts_interpolation}, Assumption \ref{independence_assumption_ZO}, \eqref{biased_estimate_two_point_estimator}, and Young's inequality,  we have %for any $t \in (kT, (k+1)T]$
        \begin{equation*} %\label{second_moment_continuous_interpolation_auxiliary_process_ZOLD}
        \begin{split}
         & \E [|\bar{X}_t -\bar{\Phi}_t^{\lambda,k } |^2 ] 
         \\ &=-2\lambda \int_{kT}^t  \E [ \langle \bar{X}_s -\bar{\Phi}_s^{\lambda,k }, G_{\eta}(\bar{X}_{\lfloor s \rfloor }, U_{\lceil s \rceil }, \xi_{\lceil s \rceil }) - \nabla f_{\eta}(\bar{\Phi}_s^{\lambda,k }) \rangle ]  \rmd s
          \\ &=-2\lambda \int_{kT}^t \E [\langle \bar{X}_{ s }  - \bar{X}_{\lfloor s \rfloor} , G_{\eta}(\bar{X}_{\lfloor s \rfloor }, U_{\lceil s \rceil }, \xi_{\lceil s \rceil }) - \nabla f_{\eta}(\bar{\Phi}_s^{\lambda,k }) \rangle ]  \rmd s
            \\ & \quad -2\lambda \int_{kT}^t \E [\langle \bar{X}_{\lfloor s \rfloor}   - \bar{\Phi}_s^{\lambda,k }, \E [ G_{\eta}(\bar{X}_{\lfloor s \rfloor }, U_{\lceil s \rceil }, \xi_{\lceil s \rceil }) - \nabla f_{\eta}(\bar{\Phi}_s^{\lambda,k }) | \mathcal{F}_{kT,  s  } ] \rangle ]  \rmd s
  \\ & =  2\lambda \int_{kT}^t (s -\lfloor s \rfloor )\E [\langle  \lambda    G_{\eta}(\bar{X}_{\lfloor s \rfloor }, U_{\lceil s \rceil }, \xi_{\lceil s \rceil })  , G_{\eta}(\bar{X}_{\lfloor s \rfloor }, U_{\lceil s \rceil }, \xi_{\lceil s \rceil }) \rangle ]  \rmd s
   \\ & \quad - 2\lambda \int_{kT}^t (s -\lfloor s \rfloor )\E [\langle  \lambda    G_{\eta}(\bar{X}_{\lfloor s \rfloor }, U_{\lceil s \rceil }, \xi_{\lceil s \rceil })  , \nabla f_{\eta}(\bar{\Phi}_s^{\lambda,k }) \rangle ]  \rmd s
   \\ & \quad  -2\lambda \int_{kT}^t \E [\langle  \sqrt{2 \lambda \beta^{-1}}  (\widetilde{B}^{\lambda}_{s} - \widetilde{B}^{\lambda}_{\lfloor s \rfloor}) , G_{\eta}(\bar{X}_{\lfloor s \rfloor }, U_{\lceil s \rceil }, \xi_{\lceil s \rceil })  \rangle ]  \rmd s
    \\ & \quad  + 2\lambda \int_{kT}^t \E [\langle  \sqrt{2 \lambda \beta^{-1}}  (\widetilde{B}^{\lambda}_{s} - \widetilde{B}^{\lambda}_{\lfloor s \rfloor}) ,  \nabla f_{\eta}(\bar{\Phi}_s^{\lambda,k }) - \nabla f_{\eta}(\bar{\Phi}_{\lfloor s \rfloor }^{\lambda,k })   \rangle ]  \rmd s
   \\ & \quad - 2\lambda \int_{kT}^t \E[ \langle \bar{X}_{\lfloor s \rfloor }  -\bar{\Phi}_s^{\lambda,k }, \E [  (\nabla f_{\eta}(\bar{X}_{\lfloor s \rfloor }) - \nabla f_{\eta}(\bar{\Phi}_s^{\lambda,k }) )|  \mathcal{F}_{kT,  s  } ] \rangle ]   \rmd s 
            \end{split}
        \end{equation*}
           \begin{equation*} %\label{second_moment_continuous_interpolation_auxiliary_process_ZOLD}
        \begin{split}
      & \le 3\lambda^2      \int_{kT}^t  \E [| G_{\eta}(\bar{X}_{\lfloor s \rfloor },   U_{\lceil s \rceil }, \xi_{\lceil s \rceil })|^2   ]  \ \rmd s
            + \lambda^2      \int_{kT}^t  \E [|\nabla f_{\eta}(\bar{\Phi}_s^{\lambda,k }) |^2 ] \ \rmd s 
            \\    & \quad   +   (2 \lambda^2 d \beta^{-1})  \int_{kT}^t  (s - \lfloor s \rfloor)  \ \rmd s + \lambda L_1^2 \E [\varphi^2(\xi)]  \int_{kT}^t  \E [| \bar{\Phi}_s^{\lambda,k } - \bar{\Phi}_{\lfloor s \rfloor}^{\lambda,k } |^2 ]     \ \rmd s 
%  \\ &=2\lambda^2 \int_{kT}^t   \E [|   G_{\eta}(\bar{\theta}_{\lfloor s \rfloor },   U_{\lceil s \rceil }, X_{\lceil s \rceil })  |^2  ]   \rmd s
 %             \\ & \quad - 2 \lambda^2  \int_{kT}^t  \E [\langle G_{\eta}(\bar{\theta}_{\lfloor s \rfloor },   U_{\lceil s \rceil }, X_{\lceil s \rceil }) , \nabla f_{\eta}(\bar{\Phi}_s^{\lambda,k } ) \rangle ] \ \rmd s
 %            \\ & \quad -2\lambda \int_{kT}^t \E[ \langle \bar{\theta}_{\lfloor s \rfloor }  -\bar{\Phi}_s^{\lambda,k }, \E [  (\nabla f_{\eta}(\bar{\theta}_{\lfloor s \rfloor }) - \nabla f_{\eta}(\bar{\Phi}_s^{\lambda,k }) )| \mathcal{F}^{(1)}_{\lfloor  s \rfloor } ] \rangle ]   \rmd s 
     % \end{split}
     %   \end{equation*}
     %     \begin{equation*}
     %   \begin{split}
       \\      & \quad + \lambda  \int_{kT}^t \E[ | \bar{X}_{\lfloor s \rfloor }  -\bar{\Phi}_s^{\lambda,k } |^2 \  \rmd s + \lambda \int_{kT}^t  \E [  |\nabla f_{\eta}(\bar{X}_{\lfloor s \rfloor }) - \nabla f_{\eta}(\bar{\Phi}_s^{\lambda,k }) |^2  ]   \   \rmd s.
                  \end{split}
        \end{equation*}
Applying  Lemma \ref{linear_growth_two_point_estimator_lemma}, Assumption \ref{independence_assumption_ZO}, the linear growth condition \eqref{linear_growth_flatness_bound}, the definition of the auxiliary process \eqref{continuous_time_process_started_interpolation_definition_ZO}, Lipschitz condition \eqref{semiconvexity_flatness_bound_remark}, $T= \lfloor 1/\lambda \rfloor$, and Lemma \ref{one_step_error_lemma_ZO}  on the right-hand side above, we obtain
%\newpage
%, for any $t \in (kT, (k+1)T]$     
         \begin{equation} \label{second_moment_continuous_interpolation_auxiliary_process_ZOLD_continuation}
        \begin{split}
          & \E [|\bar{X}_t -\bar{\Phi}_t^{\lambda,k } |^2 ]     
          \\   & \le 3\lambda^2       \int_{kT}^t   4 (4+d)   L_1^2  \E[\varphi^2 (\xi_{\lceil s \rceil }) ] \E [ |\bar{X}_{\lfloor s \rfloor }|^2 ]  \ \rmd s 
       \\ & \quad +  3\lambda^2 \int_{kT}^t   8 (4+d)(  L_2^2 \E[\widetilde{\varphi}^2(\xi_{\lceil s \rceil })] +  |\nabla F(0,0)|^2 )  \ \rmd s 
    %                         \end{split}
    %    \end{equation}
    %        \begin{equation*} %\label{second_moment_continuous_interpolation_auxiliary_process_ZOLD_continuation}
     %   \begin{split}
     \\     & \quad + 3\lambda^2  \int_{kT}^t  2  \eta^2  (216+d^3) L_1^2  \E[\varphi^2 (\xi_{\lceil s \rceil })]  \ \rmd s 
          \\   & \quad + \lambda^2 \int_{kT}^t  (2 L_1^2 \E [\varphi^2 (\xi_{\lceil s \rceil })]  \E [| \bar{\Phi}_s^{\lambda,k } |^2] +  4  L_1^2 \E [\varphi^2 (\xi_{\lceil s \rceil })] \eta^2 d ) \ \rmd s  
        \\ & \quad + \lambda^2 \int_{kT}^t 8     (L_2^2 \E[\widetilde{\varphi}^2(\xi_{\lceil s \rceil })] + | \nabla F(0, 0) |^2)     \ \rmd s
             \\ &  \quad + 2 \lambda d \beta^{-1} + \lambda L_1^2 \E [\varphi^2(\xi)]  \int_{kT}^t  (2 \lambda^2 (s - \lfloor s \rfloor  ) \int_{\lfloor s \rfloor }^s \E [|\nabla f_{\eta}(\bar{\Phi}_r^{\lambda,k }) |^2 ] \ \rmd r )
             \\ & \quad +  \lambda L_1^2 \E [\varphi^2(\xi)]  \int_{kT}^t 4 \lambda d \beta^{-1} (s - \lfloor s \rfloor)  \ \rmd s  
            \\ & \quad +  \lambda \int_{kT}^t  (2+ 2 L^2_1 \E [\varphi^2 (\xi_{\lceil s \rceil })])   \E[ | \bar{X}_{\lfloor s \rfloor }  -\bar{X}_{ s  }  |^2 ] \  \rmd s  
             \\ & \quad + \lambda  \int_{kT}^t (2+ 2 L^2_1 \E [\varphi^2 (\xi_{\lceil s \rceil })])   \E[ | \bar{X}_{ s  }  -\bar{\Phi}_s^{\lambda,k } |^2 ] \  \rmd s
            \\ & \le \lambda         (12(4+d)   \bar{L}_1   C_1(1+\mathbb{E} [ |x_0|^2 ] ) 
            + 24(4+d)\bar{L}_{2, \nabla F} + 6   \bar{L}_1  \eta^2 (216+d^3)) 
            \\ & \quad  + \lambda ( \bar{L}_1  (2 C_1 +4d/\beta+4b) (1+\mathbb{E} [ |x_0|^2 ] ) +  4 \bar{L}_1  \eta^2 d  + 2 d \beta^{-1})
            \\ & \quad + \lambda ( 4  \bar{L}_1^2 
(C_1+2d/\beta+2b)(1+\E[|x_0|^2]) + 2  \bar{L}_1(4 \bar{L}_1 \eta^2d
+8 \bar{L}_{2, \nabla F}))  
%\\ & \quad + \lambda (2  \bar{L}_1(4 \bar{L}_1 \eta^2d
%+8 \bar{L}_{2, \nabla F}))
            \\ & \quad +  \lambda (2 + 2 \bar{L}_1)    C_2 (  1+   \mathbb{E}[|x_0|^2]  ) +   \lambda (2 + 2 \bar{L}_1)  \int_{kT}^t \E[ | \bar{X}_{ s  }  -\bar{\Phi}_s^{\lambda,k } |^2 ] \ \rmd s,
    %        \\ & \textcolor{red}{OLD}
  %          \\ & \le 3\lambda         (4 (4+d)   L_1^2 \E[\varphi^2 (X) ]  C_1\left(1+\mathbb{E} [ |\theta_0|^2 ] \right)  + 2   L_1^2 \eta^2 (216+d^3) \E[\varphi^2 (X)] )
  %          \\ & \quad + 3\lambda (8(4+d)( L_2^2   \E[\widetilde{\varphi}^2(X)] +  |\nabla F(0,0)|^2)) 
 %           \\ & \quad + \lambda   (2 L_1^2 \E [\varphi^2 (X)]  (C_1 +2d/\beta+2b) \left(1+\mathbb{E} [ |\theta_0|^2 ] \right) +  4 L_1^2 \E [\varphi^2 (X)] \eta^2 d) 
 %           \\ & \quad + 2 \lambda (1+ L^2_1 \E [\varphi^2 (X)]) \left(  C_2 (  1+   \mathbb{E}[|\theta_0|^2]  ) + \int_{kT}^t \E[ | \bar{\theta}_{ s  }  -\bar{\Phi}_s^{\lambda,k } |^2 ] \right) .
           \end{split}
        \end{equation}
  where
  \begin{equation} \label{L_1_bar_L_2_bar_definition}
     \bar{L}_1 := L^2_1 \E [\varphi^2 (\xi)], \qquad  \bar{L}_{2, \nabla F} := L_2^2   \E[\widetilde{\varphi}^2(\xi)] +  |\nabla F(0,0)|^2.
  \end{equation}
     Using Gr{\"o}nwall's inequality in \eqref{second_moment_continuous_interpolation_auxiliary_process_ZOLD_continuation} together with $T= \lfloor 1/\lambda \rfloor$, we obtain 
\begin{equation*}
      W_2(\mathcal{L}(\bar{X}_t),\mathcal{L}(\bar{\Phi}_t^{\lambda,k }))  \leq \sqrt{\E [ | \bar{X}_t- \bar{\Phi}_t^{\lambda,k }|^2] }   \le C_{3} \sqrt{\lambda},
\end{equation*} 
where 
    \begin{equation} \label{constant_first_inequality_split}
    \begin{split}
        C_{3} & = e^{   1+ \bar{L}_1 }  
        \Bigg( 12(4+d)   \bar{L}_1   C_1\left(1+\mathbb{E} [ |x_0|^2 ] \right)
                   + 24(4+d) \bar{L}_{2, \nabla F}
            \\ & \qquad \qquad \   +     \bar{L}_1   (2 C_1 +4d/\beta+4b) (1+\mathbb{E} [ |x_0|^2 ] ) + 2 d \beta^{-1} + 4 \bar{L}_1   \eta^2 d  
 %                       \end{split}
 %   \end{equation}  
 %   \begin{equation*} 
    %\begin{split}
     \\ & \qquad \qquad  \     
              + 6  \bar{L}_1  \eta^2 (216+d^3) 
              + 4  \bar{L}_1^2  (C_1+2d/\beta+2b)(1+\E[|x_0|^2])
\\ & \qquad  \qquad  \   +  2  \bar{L}_1  (4\bar{L}_1 \eta^2d
+8 \bar{L}_{2, \nabla F})   
            +  ( 2 + 2 \bar{L}_1 )   C_2 (  1+   \mathbb{E}[|x_0|^2]  ) \Bigg)^{\frac{1}{2}},
               \end{split}
    \end{equation}
    with $\bar{L}_1$ and $\bar{L}_{2, \nabla F}$ given in \eqref{L_1_bar_L_2_bar_definition}, $C_1 = \mathcal{O} ( 1+ \eta^2 d^2 + d /\beta )$ defined in \eqref{eqProp2_constants} and $C_2 =  \mathcal{O}( 1+ \eta^2 d^2 + d /\beta  )$ defined in \eqref{constants_one_step_error_lemma_ZO}, respectively, and $b$ given in \eqref{constants_dissipativity_flatness_full_gradient_ZO}.
\end{proof}

 Recall that $(\mathcal{H}_t)_{t \ge 0}$ is the natural filtration of standard $d$-dimensional Brownian motion $(B_t)_{t \ge 0}$ in \eqref{zeroth_order_Langevin_SDE}, and we denote by $\mathcal{H}_{\infty}$, the sigma-algebra generated by $\cup_{t \ge 0} \mathcal{H}_t$. The upper bounds for the second and third term on the right-hand side of \eqref{decomposition_triangle_inequality_Wasserstein} are derived using the following contraction estimate for the overdamped Langevin SDE \eqref{zeroth_order_Langevin_SDE}. %established in \cite[Proposition 4]{johnston2025the} and \cite[Proposition 5]{johnston2025the}.
\begin{proposition}\label{prop2_ZO_before}
 Let Assumptions \ref{independence_assumption_ZO},  \ref{Lipschitz_assumption_full_gradient}, and \ref{strong_convexity_at_infinity_assumption_ZO} hold. Let $(\widetilde{Y}_t)_{t \in \mathbb{R}_{+}}$ be the solution of \eqref{zeroth_order_Langevin_SDE} with initial condition $\widetilde{Y}_0 = \widetilde{x}_0$, which is independent of $\mathcal{H}_{\infty}$ and assume $ \E [|\widetilde{x}_0|^2 ] < \infty $. Then, for $p \in \{ 1,2\}$, \begin{equation*}
    W_p(\mathcal{L}(Y_t),\mathcal{L}(\widetilde{Y}_t))\leq C_{W_p}e^{-C_{r_p}t} \Delta_0^{(p)}, 
\end{equation*}
where
\begin{equation} \label{prop2_constant_CW1_Cr1}
\begin{split}
 C_{W_1} &  =2e^{\beta L_1 \E[\varphi(\xi)] R^2/8},  
  \\  C_{r_1} & = 
      \left((4 /3) (\beta e)^{-1}) \min(1/R^2, \mu \E [\rho(\xi)] \beta/8) \right) \mathbbm{1}_{ \beta L_1 \E[\varphi(\xi)] R^2\leq8}
      \\ & \quad + \Bigg(\frac{4 \sqrt{2 \beta \pi} ((\beta L_1 \E[\varphi(\xi)])^{-1}+(\beta \mu  \E [\rho(\xi)])^{-1}) e^{ \beta L_1 \E[\varphi(\xi)] R^2 8^{-1}}}{R ( L_1 \E[\varphi(\xi)])^{1/2}} 
      \\ & \qquad \qquad  +16 \beta(\beta\mu\E [\rho(\xi)] R)^{-2}\Bigg)^{-1} \mathbbm{1}_{ \beta L_1 \E[\varphi(\xi)] R^2 > 8},
      \\   \Delta_0^{(1)}  & = W_1(\mathcal{L}(x_0),\mathcal{L}(\widetilde{x}_0)), 
%\frac{4}{3 \beta e}\min(1/R^2, \mu \E [\rho(X)] \beta/8) \1_{\beta L_1 \E[\varphi(X)] R^2\leq8} \\ & \quad +
%2 \beta^{-1} \Bigg(8\sqrt{2\pi}R^{-1}(\beta L_1 \E[\varphi(X)])^{-1/2}((\beta L_1 \E[\varphi(X)])^{-1}+(\beta \mu  \E [\rho(X)])^{-1}) \\ & \qquad  \quad \times \exp(\beta L_1 \E[\varphi(X)] R^2/8)+32(\beta\mu  \E [\rho(X)] R)^{-2}\Bigg)^{-1} \1_{\beta L_1 \E[\varphi(X)] R^2\geq8},
% \\  C_0^{''}(\widetilde{\zeta}) & = \left(\frac{2e^2}{\widetilde{\zeta}}\left(1+\frac{2}{\sqrt{\widetilde{\zeta}}}\right)\sqrt{\frac{2}{\sqrt{\beta/8} \mu \E [\rho(X)] -\widetilde{\zeta}}} \right) 
 %   \\ & \quad \vee  \left( \frac{2+\sqrt{\widetilde{\zeta}}}{\widetilde{\zeta}(1-e^{-2})}\left[\frac{2\sqrt{2}e^2}{\sqrt{\widetilde{\zeta}\sqrt{\beta/8}(\mu \E [\rho(X)] -\widetilde{\zeta})}}+\frac{1}{\sqrt{\beta/8}\mu \E [\rho(X)] -\widetilde{\zeta}}\right]\right)
%\\    C_{W_2} & =2\max\{1,R^{-1/2}\}C_0^{''}(\widetilde{\zeta} )e^{\left(\sqrt{\beta/32}(\mu \E[\rho(X)] + L_1 \E[\varphi(X)])+ \widetilde{\zeta} /2\right)\beta R^2/2} \\ & \quad \times \sqrt{(2/\beta)\max\{4/\widetilde{\zeta}  +2,8/(e \widetilde{\zeta} ^2)\}/(\sqrt{\beta/2}R+1)},
%    \\ C_{r_2} & =2\min\{1,1/\widetilde{\zeta}  \}e^{-(1/4)\sqrt{(\beta/2)^3}(\mu \E[\rho(X)]+ L_1 \E[\varphi(X)])R^2}/C_0^{''}(\widetilde{\zeta} ), 
%    \\ \epsilon_0^{(1)} & = W_1(\mathcal{L}(\theta_0),\mathcal{L}(\widetilde{\theta}_0)),
%    \quad  \epsilon_0^{(2)}  = \max \{W_2(\mathcal{L}(\theta_0),\mathcal{L}(\widetilde{\theta}_0)),(W_1(\mathcal{L}(\theta_0),\mathcal{L}(\widetilde{\theta}_0)))^{\frac{1}{2}} \}.
\end{split}
    \end{equation}
    and, for any $\delta \in(0,\sqrt{\beta/8}\mu \E[\rho(\xi)])$, \begin{equation} \label{prop2_constant_CW1_Cr1_continuation}
\begin{split}
       C_{W_2} & =2\max\{1,R^{-1/2}\}C_0^{''}(\delta )e^{\left(\sqrt{\beta/32}(\mu \E[\rho(\xi)] + L_1 \E[\varphi(\xi)])+ \delta /2\right)\beta R^2/2} \\ & \quad \times \sqrt{(2/\beta)\max\{4/\delta  +2,8/(e  \delta^2)\}/(\sqrt{\beta/2}R+1)},
    \\ C_{r_2} & =2\min\{1,1/\delta  \}e^{-(1/4)\sqrt{(\beta/2)^3}(\mu \E[\rho(\xi)]+ L_1 \E[\varphi(\xi)])R^2}/C_0^{''}(\delta ),
       \end{split}
    \end{equation}
    \begin{equation*} 
\begin{split}
      C_0^{''}(\delta) & = \max \Bigg \{ \frac{2e^2}{\delta}\left(1+\frac{2}{\delta^{1/2}}\right) \left(\frac{2}{(\beta/8)^{1/2} \mu \E [\rho(\xi)] -\delta}\right)^{\frac{1}{2}}, 
    \\ & \qquad \qquad    \frac{2\sqrt{2}e^2 (2+\delta^{\frac{1}{2}})}{\delta(1-e^{-2}) (\delta (\beta/8)^{\frac{1}{2}} \mu \E [\rho(\xi)] -\delta^2)^{\frac{1}{2}}}
    \\ & \qquad \qquad +\frac{2+\delta^{\frac{1}{2}}}{\delta(1-e^{-2})((\beta/8)^{1/2} \mu \E [\rho(\xi)] -\delta)} \Bigg \},
    \\  \Delta_0^{(2)}  & = \max \{W_2(\mathcal{L}(x_0),\mathcal{L}(\widetilde{x}_0)),(W_1(\mathcal{L}(x_0),\mathcal{L}(\widetilde{x}_0)))^{\frac{1}{2}} \}.
\end{split}
    \end{equation*}
    In particular, $C_{W_1} =\mathcal{O}(e^{ R^2\beta})$, $C_{r_1} =\mathcal{O}(\beta^{1/2}e^{-R^2\beta})$, $C_{W_2}=\mathcal{O}(\beta^{-1}e^{R^2\beta^{3/2}})$,  $C_{r_2}= \mathcal{O}(e^{-R^2\beta^{3/2}})$.
    \end{proposition}
    \begin{proof}
The Wasserstein-$1$ bound follows along the same lines as the proof of \cite[Proposition 4]{johnston2025the}, with the modifications that the semiconvexity (or one-sided Lipschitz) constant $K$ in \cite[Assumption 4]{johnston2025the} is replaced by the Lipschitz constant $L_1\E[\varphi(\xi)]$ for $\nabla f_{\eta}$ established in \eqref{semiconvexity_flatness_bound_remark}, and the strong convexity-at-infinity constant $\mu$ in \cite[Assumption 2]{johnston2025the} is replaced by the corresponding strong convexity-at-infinity constant $\mu\E[\rho(\xi)]$ for $\nabla f_{\eta}$ established in \eqref{strong_convexity_flatness_bound_remark}. The Wasserstein-$2$ bound follows similarly by replacing the corresponding constants in the proof of \cite[Proposition 5]{johnston2025the}.
%The bound in Wasserstein-1 follows along the same lines of the proof  of \cite[Proposition 4]{johnston2025the} with the difference that we replace the semi-convexity constant $K$ in \cite[Assumption 4]{johnston2025the} with the Lipschitz constant $ L_1  \E[\varphi(\xi)] $ of $\nabla f_{\eta}$ established in \eqref{semiconvexity_flatness_bound_remark} and the strong convexity-at-infinity constant $\mu $ in \cite[Assumption 2]{johnston2025the} with the same type of constant $ \mu  \E[\rho(\xi)] $ of $\nabla f_{\eta}$ established in \eqref{strong_convexity_flatness_bound_remark}.  
%The bound in Wasserstein-2 follows replacing the aforementioned constants in the proof of  \cite[Proposition 5]{johnston2025the}.
%The estimate for $p=1$ follows by using \eqref{strong_convexity_flatness_bound_remark} and \eqref{semiconvexity_flatness_bound_remark} in \cite[Theorem 1, Corollary 2, and Lemma 1]{eberle2016reflection}. The bound for $p=2$ is obtained by using \eqref{strong_convexity_flatness_bound_remark} and \eqref{semiconvexity_flatness_bound_remark} in \cite[Theorem 1.3]{luo2016exponential}.
    \end{proof}

%Denote by $(\mathcal{F}_t)_{t \ge 0}$ the natural filtration of $(B_t)_{t \ge 0}$ and by $\Sigma_{\theta_0}$ the sigma-algebra generated by $\theta_0$, and we assume that $(\mathcal{F}_t)_{t \ge 0}$ is independent of $ \mathcal{X}_{\infty}  \vee   \Sigma_{\theta_0}$, where $ \mathcal{X}_{\infty}$ is the sigma-algebra generated by $\cup_{n \in \mathbb{N}} \mathcal{X}_n$ with $\mathcal{X}_n$ being a given filtration representing the flow of past information.

    \begin{lemma}\label{Lemma7_ZO_before}
Let Assumptions \ref{independence_assumption_ZO},  \ref{Lipschitz_assumption_full_gradient}, and \ref{strong_convexity_at_infinity_assumption_ZO} hold. Then, for any $\lambda\in(0,\lambda_{\text{max}})$  with $\lambda_{\text{max}}$ given in \eqref{expression_lambda_max_ZO},  $t\in [kT,(k+1)T]$, $k \in\mathbb{N}$, and $p \in \{1,2\}$,
 \begin{equation*}
     W_p(\mathcal{L}(\bar{\Phi}_t^{\lambda,k}),\mathcal{L}(Y_{t}^{\lambda}))\leq C_{4,p} \lambda^{1/(2p)},
 \end{equation*}
        where
        \begin{equation} \label{Lemma7_ZO_constants}
        \begin{split}
            C_{4,1}  =  \frac{ C_{3} C_{W_1}}{1-e^{-C_{r_1}/2}},
               \qquad     C_{4,2}  = \frac{\max\{C_{3},\sqrt{C_{3}}\}C_{W_2}}{1-e^{-C_{r_2}/2}},
             \end{split}
        \end{equation}
 with the constants $C_{W_1} =\mathcal{O}(e^{ R^2\beta})$, $C_{r_1} =\mathcal{O}(\beta^{1/2}e^{-R^2\beta})$, $C_{W_2}=\mathcal{O}(\beta^{-1}e^{R^2\beta^{3/2}})$,  $C_{r_2}= \mathcal{O}(e^{-R^2\beta^{3/2}})$ given explicitly in Proposition \ref{prop2_ZO_before}, and the constant $C_3 =  \mathcal{O}(1+ d \beta^{- 1/2}+ \eta d^{3/2} + d^{1/2})$ with its full expression provided in \eqref{constant_first_inequality_split}.
        \end{lemma}
\begin{proof}
It follows by using the auxiliary process defined in \eqref{continuous_time_process_started_interpolation_definition_ZO} together with the estimates in Proposition \ref{prop2_ZO_before} and Lemma \ref{first_upper_bound_lemma}, by adapting the arguments in the proofs of \cite[Lemma~8]{johnston2025the} and \cite[Lemma~9]{johnston2025the}.
%It follows using the auxiliary process defined in \eqref{continuous_time_process_started_interpolation_definition_ZO} and the estimates in Proposition \ref{prop2_ZO_before} and Lemma \ref{first_upper_bound_lemma} in the proofs of \cite[Lemma 8]{johnston2025the} and \cite[Lemma 9]{johnston2025the}.   
%In particular, $C_{4,1}  = \mathcal{O}\left(1+  \eta \sqrt{d} + \sqrt{d/\beta} +  \beta^{-1/4} \sqrt{d \eta}\right)$ and $C_{4,2}  = \mathcal{O}\left(1+  \eta \sqrt{d} + \sqrt{d/\beta} +  \beta^{-1/4}\sqrt{d \eta}\right)$.
\end{proof}
%Next, under Assumption \ref{higher_order_regularity_objective_assumption}, we control 
The distance between $f_{\eta}$ and $v_{\eta}$ is controlled using Assumption \ref{higher_order_regularity_objective_assumption}.
%in a similar way to \cite[Lemma 3.5]{zhang2025zeroth}.
%Lastly, we derive an upper bound for controlling the distance between $\pi_{\beta}$, the invariant measure of the Langevin SDE \eqref{Langevin_SDE_equation}, and the target measure  \eqref{Hessian-trace_regularized_objective} using the following result.
%\begin{lemma}[Theorem 1 in \cite{nesterov2017random}] \label{remainder_lemma_objective}
 %    Let Assumptions \ref{independence_assumption_ZO} and \ref{Lipschitz_assumption_full_gradient} hold. Then,
 %      \begin{equation*}
%        | f_{\eta}(\theta) - f(\theta)| \le \frac{L_1 \E [\varphi(X)] \eta^2 d}{2} , \qquad  \ \theta \in \R^d.
%    \end{equation*}
%\end{lemma}
%\begin{proof}
%    It follows by using Assumptions \ref{independence_assumption_ZO} and \ref{Lipschitz_assumption_full_gradient} in the proof of Theorem 1 in \cite{nesterov2017random}.
%\end{proof}
%The proof of Theorem \ref{remainder_lemma_objective} is provided in Appendix in \cite{nesterov2017random}.
\begin{lemma}\label{remainder_lemma}
    Let Assumptions \ref{independence_assumption_ZO} and \ref{higher_order_regularity_objective_assumption} hold. Then, for any $x \in \mathbb{R}^d$,
    \begin{equation*}
        | f_{\eta}(x) - v_{\eta}(x)| \le \frac{L_3 \E[\bar{\varphi}(\xi)]}{24} \eta^4 (d+4)^2.
    \end{equation*}
\end{lemma}
\begin{proof}
It follows by using Assumptions \ref{independence_assumption_ZO} and \ref{higher_order_regularity_objective_assumption} in the proof of \cite[Lemma 3.5]{zhang2025zeroth}.
 %  The result follows by adapting the proof of \cite[Lemma 3.5]{zhang2025zeroth} under Assumptions \ref{independence_assumption_ZO} and \ref{higher_order_regularity_objective_assumption}.
\end{proof}
%The proof of Lemma \ref{remainder_lemma} is provided in Appendix B.2 in \cite{zhang2025zeroth}.
The upper bound for the last term on the right-hand side of \eqref{decomposition_triangle_inequality_Wasserstein} is derived using Lemma \ref{remainder_lemma}.
%Using Lemma \ref{remainder_lemma}, we control the last upper bound of the right-hand side of \eqref{decomposition_triangle_inequality_Wasserstein} between the invariant measure of the overdamped Langevin SDE \eqref{Langevin_SDE_equation} and the Gibbs measure  \eqref{Hessian-trace_regularized_objective}.
%The last upper bound of the decomposition \eqref{decomposition_triangle_inequality_Wasserstein} is controlled using 
%corresponding to  We can use Lemma \ref{remainder_lemma} to derive 
\begin{proposition} \label{Wasserstein_KL_convergence_two_measures_theorem_ZO}
 %   Let Assumptions \ref{independence_assumption_ZO}, \ref{Lipschitz_assumption_full_gradient}, and \ref{strong_convexity_at_infinity_assumption_ZO} hold. Then
%      \begin{equation*}
%    \begin{split}
%       W_2( \pi_{\beta}, \pi_{\beta}^{f})   \le \underline{D}^{f},
%    \end{split}
%\end{equation*}
%where $\underline{D}^{f}=   \mathcal{O}(\sqrt{\beta} \eta \sqrt{d} + \sqrt{\beta} \eta^2 d + \eta d^{\frac{3}{2}})$, whose explicit expression is given in \eqref{Wassserstein_L_divergence_corollary_final_objective}. Furthermore, suppose that Assumption
%\ref{higher_order_regularity_objective_assumption}
%also holds. Then, 
Let Assumptions \ref{independence_assumption_ZO}, \ref{Lipschitz_assumption_full_gradient}, \ref{strong_convexity_at_infinity_assumption_ZO}, and
\ref{higher_order_regularity_objective_assumption} hold. Then,
    \begin{equation*}
    \begin{split}
       W_2( \pi_{\beta}^{f_{\eta}}, \pi_{\beta}^{v_{\eta}})   \le C_5, 
    \end{split}
\end{equation*}
where $C_5=   \widetilde{\mathcal{O}}(   \beta^{\frac{1}{4}}(\eta d^{\frac{1}{2}} + \eta^3 d^{\frac{3}{2}}) + \beta^{\frac{1}{2}}(\eta^2 d + \eta^4 d^{2}) + \eta^2 d^{\frac{3}{2}}  + \beta^{ - \frac{1}{4}} \eta d )$, and its explicit expression is given in \eqref{Wassserstein_L_divergence_corollary_final}. 
\end{proposition}
\begin{proof}
Let  $\text{KL}(\pi_{\beta}^{f_{\eta}} | | \pi_{\beta}^{v_{\eta}}) $ be the Kullback-Leibler divergence between $\pi_{\beta}^{f_{\eta}}$ and $\pi_{\beta}^{v_{\eta}}$.  Let $ Z_{\beta, f_{\eta}} :=   \int_{\mathbb{R}^d} e^{- \beta f_{\eta}(x)} \ \rmd x$ %$Z_{\beta, f} :=   \int_{\mathbb{R}^d} e^{- \beta f(\theta)} \ \rmd \theta $, 
 and  $Z_{\beta, v_{\eta}} :=   \int_{\mathbb{R}^d} e^{- \beta v_{\eta}(x)} \ \rmd x  $. 
 %$\pi_{\beta}^{f}$ and between $\pi_{\beta}$ and respectively.
Using Lemma \ref{remainder_lemma}, we have
%	\begin{equation} \label{KL_divergence_pi_beta_objective}
%		\begin{split}
%			\text{KL}(\pi_{\beta} | | \pi_{\beta}^{f}) 
%			& =    \int_{\mathbb{R}^d}   \log \left( \frac{Z_{\beta}^{-1} e^{- \beta f_{\eta}(\theta)} }{ Z_{\beta, f }^{-1}  e^{- \beta f(\theta)}}  \right)  \pi_{\beta} ( \rmd \theta)
%			\\ & = \log \left( \frac{\int_{\mathbb{R}^d} e^{- \beta f_{\eta}(\theta) + \beta(f_{\eta}(\theta)-f(\theta))} \ \rmd \theta}{  \int_{\mathbb{R}^d} e^{- \beta f_{\eta}(\theta)} \ \rmd \theta}  \right)
%          \\ & \quad    - \beta   \int_{\mathbb{R}^d} (f_{\eta}(\theta) - f(\theta))   \pi_{\beta} ( \rmd \theta)     
%            \\ & \le  L_1 \E[\varphi(X)] \beta \eta^2 d =: \underline{B}^{f}.
%		\end{split}
%	\end{equation}
%Similarly, we obtain using Lemma \ref{remainder_lemma},
	\begin{equation} \label{KL_divergence_pi_beta}
		\begin{split}
			\text{KL}(\pi_{\beta}^{f_{\eta}} | | \pi_{\beta}^{v_{\eta}}) 
			& =    \int_{\mathbb{R}^d}   \log \left( \frac{Z_{\beta, f_{\eta}}^{-1} e^{- \beta f_{\eta}(x)} }{ Z_{\beta, v_{\eta} }^{-1}  e^{- \beta v_{\eta}(x)}}  \right)  \pi_{\beta}^{f_{\eta}} ( \rmd x)
			\\ & = \log \left( \frac{\int_{\mathbb{R}^d} e^{- \beta f_{\eta}(x) + \beta(f_{\eta}(x)-v_{\eta}(x))} \ \rmd x}{  \int_{\mathbb{R}^d} e^{- \beta f_{\eta}(\theta)} \ \rmd x}  \right)
            \\ & \quad  - \beta   \int_{\mathbb{R}^d} (f_{\eta}(x) - v_{\eta}(x))   \pi_{\beta}^{f_{\eta}} ( \rmd x)     
             \\ & \le  \frac{L_3 \E[\bar{\varphi}(\xi)] \beta \eta^4 (d+4)^2}{12} .
		\end{split}
	\end{equation}
As a consequence of Remark \ref{dissipativity_assumption_remark_ZO}, all the critical points of $x \mapsto f_{\eta}(x)$ are contained in the ball $\mathfrak{D} := \{ x \in \mathbb{R}^d : |x | \le \sqrt{b/a} \}$, see, e.g., \cite[Lemma D.1]{bruno2026flatnessaware}. We denote the complement of this set by $\mathfrak{D}^c$. Using the dissipativity of $f_{\eta}$, we obtain, for $x \in \mathfrak{D}^c$ and  $\widetilde{c} \in (0,1)$, 
	\begin{equation} \label{quadratic_lower_bound_dissipative_g_epsilon}
		\begin{split}
			f_{\eta}(x) & =  f_{\eta}(\widetilde{c} x) + \int_{\widetilde{c}}^1 \langle x, \nabla f_{\eta} (t x)\rangle \ \rmd t
			\\ &  \ge f_{\eta}( x_{f_{\eta}}^{\star}) + \int_{\widetilde{c}}^1 t^{-1} \langle t x, \nabla f_{\eta} (t x)\rangle \ \rmd t
            	\end{split}
	\end{equation}
    	\begin{equation*} %\label{quadratic_lower_bound_dissipative_g_epsilon}
		\begin{split}
			 &  \ge f_{\eta}( x_{f_{\eta}}^{\star})  + \int_{\widetilde{c}}^1 t^{-1} (a |t x|^2 - b)  \ \rmd t
		%	\\ &  \ge  \frac{a (1-\widetilde{c}^2)}{2} | \theta|^2 + b \log \widetilde{c} + f_{\eta}( \theta_{f_{\eta}}^{\star}) 
			\\ &   \ge   \bar{c} | x|^2 +  b \log \widetilde{c} + f_{\eta}( x_{f_{\eta}}^{\star}) ,
		\end{split}
	\end{equation*}
	where $\bar{c}:= \frac{a (1-\widetilde{c}^2)}{2}>0$ and $x_{f_{\eta}}^{\star}$ is a minimizer of $f_{\eta}$. Let $\gamma \in (0, \beta \bar{c})$ and $c_{v_{\eta}, \beta}:= \min_{x \in \mathfrak{D}} e^{  - \beta   v_{\eta}(x)  }$. %Using Lemma \ref{remainder_lemma}, \eqref{quadratic_lower_bound_dissipative_g_epsilon}, and Lemma \ref{set_A_remark_dissipativity_ZO}, we have  
   Using Lemma \ref{remainder_lemma}, the inequality \eqref{quadratic_lower_bound_dissipative_g_epsilon}, and $| x | \le \sqrt{b/a}$ for $x \in \mathfrak{D}$, we have 
    \begin{equation} \label{exponential_moments_target_measure}
        \begin{split}
&   \int_{\mathbb{R}^d}  e^{ (\beta \bar{c} - \gamma)| x |^2}  \pi_{\beta}^{v_{\eta}}( \rmd x ) 
%\\ & =   \frac{\int_{\mathbb{R}^d}   e^{ (\beta \bar{c} - \gamma)| \theta |^2 -  \beta    v(\theta)  }   \ \rmd \theta }{ \int_{\mathbb{R}^d}    e^{  - \beta   v(\theta)  }   \ \rmd \theta} 
\\ & \le    \frac{\int_{\mathfrak{D}}   e^{ (\beta \bar{c} - \gamma)| x |^2 -  \beta    v_{\eta}(x)  }   \ \rmd x }{ \int_{\mathfrak{D}}    e^{  - \beta   v_{\eta}(x)  }   \ \rmd x} + \frac{\int_{\mathfrak{D}^{c}}   e^{ (\beta \bar{c} - \gamma)| x |^2 -  \beta  f_{\eta}(x) + \beta ( f_{\eta}(x) -   v_{\eta}(x))  }   \ \rmd x }{ \int_{\mathfrak{D}}    e^{  - \beta   v_{\eta}(x)  }   \ \rmd x} 
  \\ & \le e^{ \frac{(\beta \bar{c} - \gamma)b}{a}  }  +   e^{\frac{L_3 \E[\bar{\varphi}(\xi)]}{24} \beta \eta^4 (d+4)^2 }   \frac{ \int_{\mathfrak{D}^c}   e^{ (\beta \bar{c} - \gamma)| x |^2 - \beta   f_{\eta}(x) }   \ \rmd x}{ \int_{\mathfrak{D}}    e^{  - \beta   v_{\eta}(x)  }   \ \rmd x}
   \\ & \le e^{ \frac{(\beta \bar{c} - \gamma)b}{a}  }  +   e^{\frac{L_3  \E[\bar{\varphi}(\xi)]}{24} \beta \eta^4 (d+4)^2 }  \frac{ e^{- \beta   b \log \widetilde{c} - \beta f_{\eta}( x_{f_{\eta}}^{\star})}  \int_{\mathfrak{D}^c}   e^{ - | x |^2  \gamma  }   \ \rmd x}{ c_{v_{\eta}, \beta} (\frac{b}{a})^{\frac{d}{2}} \frac{\pi^{d/2}}{\Gamma(d/2+1)} }
 %   \\ & \le e^{ \frac{(\beta \bar{c} - \gamma)b}{a}  }    + \frac{ e^{- \beta   b \log \widetilde{c} - \beta f_{\eta}( \theta_{f_{\eta}}^{\star}) + \frac{L_3 \E[\bar{\varphi}(X)]}{24} \beta \eta^4 (d+4)^2 } \int_{\mathfrak{B}^c}   e^{ - | \theta |^2  \gamma } \ \rmd \theta   }{ c_{v_{\beta}} (\frac{b}{a})^{\frac{d}{2}} \frac{\pi^{d/2}}{\Gamma(d/2+1)} } 
   \\ & \le e^{ \frac{(\beta \bar{c} - \gamma)b}{a}  }    + \frac{ e^{- \beta   b \log \widetilde{c} - \beta f_{\eta}( x_{f_{\eta}}^{\star}) + \frac{L_3 \E[\bar{\varphi}(\xi)]}{24} \beta \eta^4 (d+4)^2 } a^{\frac{d}{2}}   \Gamma(\frac{d}{2}+1)  }{ c_{v_{\eta},\beta}  \gamma^{\frac{d}{2}} b^{\frac{d}{2}}  } < \infty.
        \end{split}
    \end{equation}
Applying the Bolley-Villani inequality \cite[Corollary 2.3]{bolley2005weighted} with \eqref{KL_divergence_pi_beta} and \eqref{exponential_moments_target_measure}, and using the inequality $\log(w+ \bar{w}) \le \log(2)+ \max \{ \log(w), \log(\bar{w})\}$ for $w, \bar{w}>0$, we have
\begin{equation*} 
\begin{split}
		W_2(\pi_{\beta}^{f_{\eta}} , \pi_{\beta}^{v_{\eta}}) 
     & \le \left[ \frac{6}{\beta \bar{c} - \gamma }  + \frac{4}{\beta \bar{c} - \gamma}    \log \left( \int_{\mathbb{R}^d}   e^{ (\beta \bar{c} - \gamma)| x |^2}  \pi_{\beta}^{v_{\eta}}( \rmd x )   \right) \right]^{\frac{1}{2}} 
    \\ & \quad \times \left[ \sqrt{\text{KL}(\pi_{\beta}^{f_{\eta}} || \pi_{\beta}^{v_{\eta}})} + 2^{-\frac{1}{4}}(\text{KL}(\pi_{\beta}^{f_{\eta}} || \pi_{\beta}^{v_{\eta}}))^{\frac{1}{4}}  \right]
     \\ & \le  C_5 ,
           \end{split}
	\end{equation*}
    where 
        \begin{equation} \label{Wassserstein_L_divergence_corollary_final} 
\begin{split}
    C_5 & := \Bigg[ \frac{6 + 4 \log(2)}{\beta \bar{c} - \gamma}  
    %+  \frac{4 \log(2)}{\beta \bar{c} - \gamma}
	 +      \max \Bigg \{ \frac{ 4 b}{a},  \frac{4 \beta   b \log \frac{1}{\widetilde{c}}}{\beta \bar{c} - \gamma} - \frac{4 \beta f_{\eta}( x_{f_{\eta}}^{\star})}{\beta \bar{c} - \gamma}    +\frac{2 d \log(a)}{ \beta \bar{c} - \gamma} 
  %                      \end{split}
%	\end{equation}
 %           \begin{equation*} %\label{Wassserstein_L_divergence_corollary_final} 
%\begin{split}  
   \\          & \qquad \qquad \qquad \qquad \qquad \quad  \  +  \frac{ L_3 \E[\bar{\varphi}(\xi)] \beta \eta^4 (d+4)^2 }{6 (\beta \bar{c} - \gamma)} -  \frac{4 \log (c_{v_{\eta}, \beta})}{\beta \bar{c} - \gamma}
               \\ & \qquad \qquad \qquad \qquad \qquad \quad \ + \frac{4 \log( \Gamma(d/2+1))}{\beta \bar{c} - \gamma} - \frac{2 d \log\left( \gamma b \right)}{\beta \bar{c} - \gamma} \Bigg \}  \Bigg]^{\frac{1}{2}}
 \\  & \quad \times \Bigg[ \sqrt{  \frac{L_3 \E[\bar{\varphi}(\xi)]\beta \eta^4 (d+4)^2}{12}  } 
 +  \left(  \frac{L_3 \E[\bar{\varphi}(\xi)]  \beta \eta^4 (d+4)^2}{24} \right)^{\frac{1}{4}} \Bigg]  .
         \end{split}
	\end{equation}
We note that 
\begin{equation} \label{estimate_last_two_terms_last_bound}
    \begin{split}
   &  \frac{-4 \beta f_{\eta}( x_{f_{\eta}}^{\star}) - 4 \log (c_{v_{\eta}, \beta})}{\beta \bar{c} - \gamma} 
     \\ & \le \frac{8 \beta L_1 \E [\varphi(\xi)] b a^{-1} +   6^{-1}  L_3 \E[\bar{\varphi}(\xi)]  \beta \eta^4 (d+4)^2  }{\beta \bar{c} - \gamma}.
    \end{split}
\end{equation}
Indeed, by \eqref{semiconvexity_flatness_bound_remark}, we have
\begin{equation*}
    f_{\eta} (x) - f_{\eta} (x_{f_{\eta}}^{\star}) \le 2^{-1} L_1 \E[\varphi(\xi)] | x - x_{f_{\eta}}^{\star}|^2 \le 2 L_1 \E[\varphi(\xi)]  b a^{-1} , \qquad x \in \mathfrak{D}, 
\end{equation*}
and by Lemma \ref{remainder_lemma}
\begin{equation*}
    \max_{x \in \mathfrak{D}} v_{\eta}(x) - f_{\eta} (x_{f_{\eta}}^{\star}) \le  2 L_1 \E[\varphi(\xi)]  b a^{-1}  + L_3 \E[\bar{\varphi}(\xi)]  \eta^4 (d+4)^2 (24)^{-1}.
\end{equation*}
Using $  \min_{x \in \mathfrak{D}} e^{-\beta v_{\eta}(x)}  =  e^{-\beta \max_{x \in \mathfrak{D}} v_{\eta}(x)}$ in $c_{v_{\eta}, \beta}$, we obtain \eqref{estimate_last_two_terms_last_bound}.
Next, we choose $\widetilde{c} = \frac{1}{2}$  in \eqref{quadratic_lower_bound_dissipative_g_epsilon}, thus $\bar{c} = \frac{3a}{8} $, and choose $\gamma = \frac{3a \beta}{16} $ in \eqref{exponential_moments_target_measure}, thus $ \beta \bar{c} - \gamma  =\frac{3 a \beta}{16} $. Therefore, the expression in the first bracket on the right-hand side of \eqref{Wassserstein_L_divergence_corollary_final} is of $\widetilde{\mathcal{O}}( 1+ \sqrt{d/\beta} + \eta^2 d )$,  since by \eqref{constants_dissipativity_flatness_full_gradient_ZO}, we have $b = \mathcal{O}(1+ \eta^2 d)$ and $a=\mathcal{O}(1)$. The expression in the second bracket on the right-hand side of \eqref{Wassserstein_L_divergence_corollary_final} is of $\mathcal{O}( \beta^{\frac{1}{2}} \eta^2 d + \beta^{\frac{1}{4}} \eta d^{\frac{1}{2}} )$.
    \end{proof}
    \begin{remark} \label{Wasserstein_KL_convergence_two_measures_theorem_ZO_remark}
    Under the coupling $ \eta = \beta^{-\frac{1+\tau}{4}}d^{-1}$,  for   $\tau \in (0,1)$, the constant
$  C_5 $ in Proposition \ref{Wasserstein_KL_convergence_two_measures_theorem_ZO} is of 
\begin{equation*}
\begin{split}
      \widetilde{\mathcal{O}}( \beta^{- \frac{\tau}{4}} d^{-\frac{1}{2}} + \beta^{- \frac{1}{2}- \frac{3 \tau}{4}} d^{- \frac{3}{2}} + \beta^{- \frac{\tau}{2}} d^{- 1} + \beta^{-\frac{1}{2} - \tau} d^{-2} + \beta^{-\frac{(1+\tau)}{2}} d^{- \frac{1}{2}} + \beta^{-  \frac{(2 + \tau)}{4}} ).
      \end{split}
\end{equation*}
    \end{remark}
\section{Proofs of the main results} \label{sec:proofs_main_results}
In this section, we provide the proofs of Theorem \ref{Theorem1_ZO}, Corollary \ref{Theorem1_ZO_W_2}, and Theorem \ref{Theorem3_ZO}.
 
\begin{proof}[Proof of Theorem \ref{Theorem1_ZO} and Corollary \ref{Theorem1_ZO_W_2}]
    Let \( K \in \mathbb{N} \) and set \( k = \lfloor K/T \rfloor \), then \( K \in [kT, (k+1)T] \). Fix  $\lambda\in(0,\lambda_{\text{max}})$, and \( t \in [kT, (k+1)T] \). 
    %Then, for $p\in\{1,2\}$, by the triangle inequality,
%\begin{equation*}
%\begin{split}
%   & W_p\left(\mathcal{L}(\theta_K^{\text{ZO}}), \pi^{\star}_{\beta,\eta}\right)
%    \\ &\leq W_p\left(\mathcal{L}(\bar{\theta}_K^{\text{ZO}}),\mathcal{L}(\bar{\Phi}_K^{\lambda,k , \text{ZO}})\right)
%    + W_p\left(\mathcal{L}(\bar{\Phi}_K^{\lambda,k , \text{ZO}}),\mathcal{L}(Y_{K}^{\lambda, \text{ZO}})\right)
%   \\ & \quad  + W_p\left(\mathcal{L}(Y_{K}^{\lambda, \text{ZO}}), \pi_{\beta}^{\text{ZO}} \right) 
%     + W_p(\pi_{\beta}^{\text{ZO}}, \pi^{\star}_{\beta,\eta} ).
%    \end{split}
%\end{equation*}
By Lemma~\ref{first_upper_bound_lemma}, Lemma \ref{Lemma7_ZO_before},  Proposition \ref{prop2_ZO_before}, Proposition \ref{Wasserstein_KL_convergence_two_measures_theorem_ZO} together with Remark \ref{Wasserstein_KL_convergence_two_measures_theorem_ZO_remark}, we have the following bounds for $p\in\{1,2\}$,
\begin{equation} \label{bounds_theorem_corollary_proofs}
\begin{split}
   & W_p(\mathcal{L}(\bar{X}_K),\mathcal{L}(\bar{\Phi}_K^{\lambda,k }))\leq C_3\lambda^{1/2}, 
    \qquad  W_p(\mathcal{L}(\bar{\Phi}_K^{\lambda,k }),\mathcal{L}(Y_{K}^{\lambda}))\leq C_{4,p}\lambda^{1/(2p)} 
\\
   & W_p(\mathcal{L}(Y_{K}^{\lambda}), \pi_{\beta}^{f_{\eta}}  )\leq C_{W_p}e^{-C_{r_p}\lambda K}\Delta_0^{(p)} , 
    \quad     W_p( \pi_{\beta}^{f_{\eta}}, \pi_{\beta}^{v_{\eta}})   \le C_5,
    \end{split}
\end{equation}
where $\Delta_0^{(1)}  =W_1( \mathcal{L}(x_0),\pi_{\beta}^{f_{\eta}})$, $\Delta_0^{(2)}  =\max\{ W_2(\mathcal{L}(x_0),\pi_{\beta}^{f_{\eta}}), (\Delta_0^{(1)})^{\frac{1}{2}} \}$. Under the coupling \eqref{coupling}, we have 
\begin{equation} \label{C_3_constant_explicit_eta}
     C_{3} = \mathcal{O}(1+ d \beta^{- 1/2}+ d^{1/2} (\beta^{-\frac{1+\tau}{4}}  + 1)),
\end{equation}
while the order of  $ C_5$ is given in Remark \ref{Wasserstein_KL_convergence_two_measures_theorem_ZO_remark}.
Using \eqref{bounds_theorem_corollary_proofs} in \eqref{decomposition_triangle_inequality_Wasserstein}, we obtain
% and 
% \begin{equation} 
% \label{value_underline_D_eta}
% \begin{split}
%    \underline{B}  =   \widetilde{\mathcal{O}}( \beta^{- \frac{\tau}{4}} d^{-\frac{1}{2}} + \beta^{- \frac{1}{2}- \frac{3 \tau}{4}} d^{- \frac{3}{2}} + \beta^{- \frac{\tau}{2}} d^{- \frac{3}{2}} + \beta^{-\frac{1}{2} - \tau} d^{-2} + \beta^{-\frac{(1+\tau)}{2}} d^{- \frac{1}{2}} + \beta^{- \frac{3}{4} - \frac{\tau}{4}} ). 
%     \end{split}
% \end{equation}
%Recall that the constants $ C_{W_1}, C_{W_2}, C_{r_1}$, and $C_{r_2}$ defined in Proposition \ref{prop2_ZO_before} do not depend on $\eta$.
%Using  the constants in \eqref{C_3_constant_explicit_eta} and \eqref{value_underline_D_eta} in the estimates \eqref{bounds_theorem_corollary_proofs}, we obtain the following bound for the right-hand side of \eqref{decomposition_triangle_inequality_Wasserstein}
%\eqref{decomposition_triangle_inequality_Wasserstein_objective} and
%where $\epsilon_0^{(1)}  =W_1( \mathcal{L}(\theta_0),\pi_{\beta})$ and $\epsilon_0^{(2)}  =\max\{ W_2(\mathcal{L}(\theta_0),\pi_{\beta}), \sqrt{\epsilon_0^{(1)}} \}$.
\begin{equation} \label{bounds_theorem_corollary_proofs_continuation}
%\begin{split}
%   W_1\left(\mathcal{L}(\theta_K), \pi_{\beta}^{f}\right) & \leq C_{W_1}e^{-C_{r_1}\lambda K}\epsilon_0^{(1)}  +C_{6,1} \lambda^{1/2} +  \underline{D}^{f},
%   \\ W_2\left(\mathcal{L}(\theta_K),\pi_{\beta}^{f} \right) & \leq C_{W_2}e^{-C_{r_2}\lambda K} \epsilon_0^{(2)} +C_{6,2}\lambda^{1/4} +  \underline{D}^{f},
%   \\
   W_p(\mathcal{L}(x_K), \pi_{\beta}^{v_{\eta}})  \leq C_{W_p}e^{-C_{r_p}\lambda K}\Delta_0^{(p)}  +C_{6,p} \lambda^{1/(2p)} +  C_5, \quad p \in \{ 1, 2\},
 %  \\ W_2(\mathcal{L}(x_K),\pi_{\beta}^{v_{\eta}} ) & \leq C_{W_2}e^{-C_{r_2}\lambda K} \Delta_0^{(2)} +C_{6,2}\lambda^{1/4} +  C_5, 
  % \end{split}
\end{equation}
with
%$\epsilon_0^{(1)}  =W_1( \mathcal{L}(\theta_0),\pi_{\beta})$, $\epsilon_0^{(2)}  =\max\{ W_2(\mathcal{L}(\theta_0),\pi_{\beta}), \sqrt{\epsilon_0^{(1)}} \}$, 
\begin{align} 
  %  \begin{split}
       C_{6,1} & := C_3 ( 1+   C_{W_1} (1-e^{-C_{r_1}/2})^{-1} ), \label{constant_C_6_1_main_proof} \\
        C_{6,2} & :=C_3+  \max\{C_{3},\sqrt{C_{3}}\} C_{W_2}(1-e^{-C_{r_2}/2})^{-1}, \label{constant_C_6_2_main_proof}
       %\left( e^{  (1 + 2  L_1^2 \E [\varphi^2(X)] )} \sqrt{ (3   + 4  L_1^2 \E [\varphi^2(X)])  C_2 (  1+   \mathbb{E}[|\theta_0|^2]  )} \right)
    %   \\ & \qquad  \times \left( 1+   C_{W_1} (1-e^{-C_{r_1}/2})^{-1}\right) 
     %  \\ & \quad \times \Bigg( 2  S \E[\varrho(X)]  C_2 (  1+   \mathbb{E}[|\theta_0|^2]   )  
%        \\ & \qquad \quad  +   \sqrt{C_2 (  1+   \mathbb{E}[|\theta_0|^2]  )}
%      \\ & \qquad  \quad \times  \Bigg[ 2 L  \E [\varphi (X)]  \sqrt{(C_1+C_3)\left(1+\mathbb{E}|\theta_0|^2\right)}  + 4  L \E [\varphi (X)] \eta \sqrt{d}+ 4  \E [\bar{\varphi}(X)]    \Bigg]\Bigg)^{\frac{1}{2}}, 
        %\\ & \quad  \times  
 %   \end{split}
%\end{equation}
%$C_{6,1}=C_4+C_{5,1}$, 
%\begin{equation}
%
%    \begin{split}
  % \end{split}
\end{align}
where the constants $C_{W_1}=\mathcal{O}(e^{ R^2\beta})$, $C_{r_1}=\mathcal{O}(\beta^{1/2}e^{-R^2\beta})$,  $C_{W_2}=\mathcal{O}(\beta^{-1}e^{R^2\beta^{3/2}})$,  $C_{r_2}= \mathcal{O}(e^{-R^2\beta^{3/2}})$, and the order of $C_3$ is given in \eqref{C_3_constant_explicit_eta}.
Therefore,
\begin{align} 
%\begin{split}
    C_{6,1}  & = \mathcal{O}(1+ e^{R^2 \beta }(1+ \beta^{-\frac{1}{2}} + d( \beta^{-\frac{1}{2}} +  \beta^{- 1}) + d^{\frac{1}{2}} (\beta^{- \frac{1+\tau}{4}} + \beta^{- \frac{3+\tau}{4}} +1))  \label{constant_C_6_1_main_proof_asymptotic_notation}
    \\ & \qquad \ + d^{\frac{1}{2}}(d^{\frac{1}{2}} \beta^{-\frac{1}{2}} + \beta^{-\frac{1+\tau }{4}} +1)), \nonumber
  \\  C_{6,2} & = \mathcal{O}(1+ e^{R^2 \beta^{\frac{3}{2}} }(\beta^{-1}+ d \beta^{- \frac{3}{2}}+ d^{\frac{1}{2}} \beta^{-\frac{(5+ \tau)}{4} } + d^{\frac{1}{2}}\beta^{-1})  \label{constant_C_6_2_main_proof_asymptotic_notation}
  \\ & \qquad \ + d^{\frac{1}{2}}(\beta^{-\frac{(1+\tau)}{4}} + d^{\frac{1}{2}} \beta^{-\frac{1}{2}} +1) ). \nonumber
  %  \end{split}
\end{align}
 For any  $\epsilon >0$, if we choose $\beta $ such that 
$ C_5  \le \frac{\epsilon}{3}$, then choose $\lambda$  such that $\lambda < \lambda_{\text{max}}$ and $C_{6,1}\lambda^{1/2} \le \frac{\epsilon}{3}$, and finally choose $K$ such that  $ C_{W_1}e^{-C_{r_1}\lambda K}\Delta_0^{(1)}  \le \frac{\epsilon}{3}$, consequently, we obtain $W_1(\mathcal{L}(x_K),\pi_{\beta}^{v_{\eta}}) \le \epsilon$.
 %  \begin{equation*}
 %      C_{W_1}e^{-C_{r_1}\lambda K}\epsilon_0^{(1)}  \le \frac{\epsilon}{3}, \qquad C_{6,1}\lambda^{1/4} \le \frac{\epsilon}{3}, 
 %      \qquad \underline{D}  \le \frac{\epsilon}{3},
 %  \end{equation*}
%   then  $W_1(\mathcal{L}(\theta_K^{\text{ZOLD}}),\pi_{\beta}^{v}) \le \epsilon$. 
   This yields 
   \begin{equation} \label{definition_beta_bar_W_1_ZOLD}
   \begin{split}
    \beta & \ge \beta_{\epsilon} :=   C_5^{0} \max\Bigg \{ \left(\frac{18}{\epsilon d^{\frac{1}{2}}} \right)^{\frac{4}{\tau}}, \left(\frac{18}{\epsilon d^{\frac{3}{2}}} \right)^{\frac{4}{2+ 3 \tau}}, \left(\frac{18}{\epsilon d} \right)^{\frac{2}{\tau}}, 
    \\ & \qquad \qquad \qquad \qquad  \left(\frac{18}{\epsilon d^{2}} \right)^{\frac{2}{1+ 2 \tau}},  \left(\frac{18}{\epsilon d^{\frac{1}{2}}} \right)^{\frac{2}{1+  \tau}}, \left(\frac{18}{\epsilon } \right)^{\frac{4}{2+  \tau}} \Bigg \},
     \end{split}
\end{equation}
 where $C_5^{0}>0$ is a constant independent of $\lambda$ and $K$,  $0< \lambda < \min \left \{ \frac{\epsilon^2}{9 C_{6,1}^2}, \lambda_{\text{max}} \right \}$, $\lambda_{\text{max}}$ given in \eqref{expression_lambda_max_ZO} and choose $\lambda = \lambda_{\star} := \frac{1}{2}   \min \left \{ \frac{\epsilon^2}{9 C_{6,1}^2}, \lambda_{\text{max}} \right \}$, then
  \begin{equation*}
       K  \ge  \left \lceil \lambda_{\star}^{-1} C_{r_1}^{-1} \ln \left( \frac{3 C_{W_1} \Delta_0^{(1)}  }{\epsilon} \right) \right \rceil  ,
 \end{equation*}
 with  $C_{6,1}$ given in \eqref{constant_C_6_1_main_proof} and its order provided in  \eqref{constant_C_6_1_main_proof_asymptotic_notation}. 
 %This concludes the proof of Theorem \ref{Theorem1_ZO}.
  % \begin{equation*}
%   \begin{split}
 %   \beta & \ge   \underline{D}^0 \max \left(  \left( \frac{12 \sqrt{d}}{\epsilon}\right)^{\frac{4}{\tau}},  \left(\frac{12 d }{\epsilon}\right)^{\frac{2}{\tau}},  \left(\frac{12 d^2}{\epsilon} \right)^{\frac{2}{1+\tau}} \right)
%    \\ & = \underline{D}^0 12^{\frac{4}{\tau}} \frac{d^{\frac{2}{\tau}}}{\epsilon^{\frac{2}{1+\tau}}},
%     \end{split}
%\end{equation*}
% where $\underline{D}^{0}>0$ is a constant independent of $\lambda$ and $K$, and
% If we choose $\beta>0$ such that the upper bound in Proposition \ref{Wasserstein_KL_convergence_two_measures_theorem_ZO} is 
%\begin{equation} \label{bound_W_2_delta_proof_distance_measure}
%    \begin{split}
%       W_2( \pi_{\beta}^{\text{ZO}}, \pi_{\beta}^{v})  & \le \underline{D} 
%       \\ & \le \frac{\epsilon}{3}, 
%    \end{split}
%\end{equation}
%then $W_1(\mathcal{L}(\theta_K^{\text{ZO}}),\pi^{\star}_{\beta,\eta} ) \le \epsilon$. 
For Corollary \ref{Theorem1_ZO_W_2}, if we choose $\beta $ as in \eqref{definition_beta_bar_W_1_ZOLD}, then choose $\lambda$ such that $\lambda < \lambda_{\text{max}}$ and $C_{6,2}\lambda^{1/4} \le \frac{\epsilon}{3}$, and finally choose $K$ such that $C_{W_2}e^{-C_{r_2}\lambda K}\Delta_0^{(2)} \le \frac{\epsilon}{3}$, thus, we obtain  $W_2(\mathcal{L}(x_K),\pi_{\beta}^{v_{\eta}}) \le \epsilon$. This yields $ 0< \lambda < \min \left\{ \frac{\epsilon^4}{81 C_{6,2}^4}, \lambda_{\text{max}} \right \}$, and choose $\lambda = \lambda_{\triangle}:= \frac{1}{2} \min \left\{ \frac{\epsilon^4}{81 C_{6,2}^4}, \lambda_{\text{max}} \right \}$ with $\lambda_{\text{max}}$ given in \eqref{expression_lambda_max_ZO}, then
 \begin{equation*}
      K  \ge \left \lceil  \lambda_{\triangle}^{-1} C_{r_2}^{-1} \ln \left( \frac{3 C_{W_2} \Delta_0^{(2)} }{\epsilon} \right) \right \rceil, 
 \end{equation*}
 with $C_{6,2}$ given in \eqref{constant_C_6_2_main_proof} with its order provided in  \eqref{constant_C_6_2_main_proof_asymptotic_notation}.
\end{proof}
\begin{proof}[Proof of Theorem \ref{Theorem3_ZO}]
%We follow \cite[Section 3.5]{pmlr-v65-raginsky17a} and apply it to the Gaussian smoothed-objective $f_{\eta}$. 
Let   $x^{\pi_\beta}$ be a random variable with law
 $ \mathcal{L}(x^{\pi_\beta})=\pi_\beta^{f_{\eta}}$. We decompose the expected excess risk using the following splitting
    \begin{equation} \label{splitting_expected_excess_risk_ZO}
\begin{split}
  &  \E  [f_{\eta}(x_K)] - \inf_{x \in \mathbb{R}^d} f_{\eta}(x)
     \\ &  = (\E  [f_{\eta}(x_K)] - \E  [f_{\eta}(x^{\pi_\beta})]) + (\E  [f_{\eta}(x^{\pi_\beta})] - \inf_{x \in \mathbb{R}^d} f_{\eta}(x)).
    \end{split}
\end{equation}
The estimate of the first difference on the right-hand side of \eqref{splitting_expected_excess_risk_ZO} follows using   \cite[Lemma 6]{pmlr-v65-raginsky17a}, the growth of $f_{\eta}$ in  \eqref{linear_growth_flatness_bound},  Lemma \ref{Lemma2_ZO}, and the first three bounds in $W_2$ in \eqref{bounds_theorem_corollary_proofs} with the coupling \eqref{coupling}, 
        \begin{equation} \label{lemmaT1_inequality}
        \begin{split}
         &  \E  [f_{\eta}(x_K)] - \E  [f_{\eta}(x^{\pi_\beta})] \\ & \leq  (  L_1 \E[\varphi(\xi)]    ( \sqrt{\max \{\mathbb{E}[|x_{K} |^2],  \E [|x^{\pi_\beta}|^2] \} })
       \\ & \qquad  \   +  \beta^{-\frac{1+\tau}{4}}d^{-\frac{1}{2}}  L_1 \E[\varphi(\xi)] 
           +    L_2 \E [\widetilde{\varphi}(\xi)] +  |\nabla F(0,0)| ) W_2(\mathcal{L}(x_K),\pi_{\beta}^{f_{\eta}})
           \\ & \leq  C_1^{\#} e^{-C_{r_2}\lambda K}  +  C_2^{\#} \lambda^{1/4}, 
 %        \\ & \leq  \Bigg(  L_1 \E[\varphi(X)]    (C_1(1+\mathbb{E}[|\theta_0|^2]) +   1+ b a^{-1} + d (\beta a)^{-1})  
%       \\ & \qquad \quad    +  \frac{ L_1 \E[\varphi(X)] }{   \beta^{\frac{1+\tau}{4}}d^{1/2}}
 %          +    L_2 \E [\widetilde{\varphi}(X)] +  |\nabla F(0,0)| \Bigg) 
%           \\ & \qquad \quad  \times (C_{W_2}e^{-C_{r_2}\lambda K} \epsilon_0^{(2)} +C_{5,2}\lambda^{1/4}),
          %\\ & \quad \times ,
     %     \\ & \quad  \times
%        \left(C_{W_2}e^{-C_{r_2}\lambda K} \max\{ W_2(\mathcal{L}(\theta_0),\pi_{\beta}), \sqrt{W_1( \mathcal{L}(\theta_0),\pi_{\beta})} \} +C_{6,2}\lambda^{1/4} \right),
            \end{split}
        \end{equation}
        %where  $\epsilon_0^{(2)}  =\max\{ W_2(\mathcal{L}(\theta_0),\pi_{\beta}^{\text{ZO}}), \sqrt{W_1( \mathcal{L}(\theta_0),\pi_{\beta}^{\text{ZO}})} \}$ and 
        where
\begin{align} 
  %  \begin{split}
     C_1^{\#}  & :=  (  L_1 \E[\varphi(\xi)]    (C_1(1+\mathbb{E}[|x_0|^2]) +   1+ b a^{-1} + d (\beta a)^{-1})  \label{constants_excess_risk_bounds_C_1}
       \\ & \qquad \quad    + \beta^{-\frac{1+\tau}{4}}d^{-\frac{1}{2}}   L_1 \E[\varphi(\xi)]  
           +    L_2 \E [\widetilde{\varphi}(\xi)] +  |\nabla F(0,0)| ) C_{W_2} \Delta_0^{(2)}, \nonumber
\\ C_2^{\#}  & := (  L_1 \E[\varphi(\xi)]    (C_1(1+\mathbb{E}[|x_0|^2]) +   1+ b a^{-1} + d (\beta a)^{-1})  \label{constants_excess_risk_bounds}
       \\ & \qquad \quad    +   \beta^{-\frac{1+\tau}{4}}d^{-\frac{1}{2}}  L_1 \E[\varphi(\xi)] 
           +    L_2 \E [\widetilde{\varphi}(\xi)] +  |\nabla F(0,0)| ) C_{6,2}, \nonumber
          % \end{split}
           \end{align}
        with $C_1 =  \mathcal{O} ( 1+  \beta^{- \frac{1+\tau}{2}} + d \beta^{-1} )$ given in  Lemma \ref{Lemma2_ZO},  $b= \mathcal{O}(1+ \beta^{- \frac{1+\tau}{2}} d^{-1})$ and $a=\mathcal{O}(1)$ defined in \eqref{constants_dissipativity_flatness_full_gradient_ZO},  $\Delta_0^{(2)}  =\max\{ W_2(\mathcal{L}(x_0),\pi_{\beta}^{f_{\eta}}), (W_1(\mathcal{L}(x_0),\pi_{\beta}^{f_{\eta}}))^{\frac{1}{2}} \}$, $C_{W_2} =  \mathcal{O}(\beta^{-1}e^{R^2\beta^{3/2}})$ and $C_{r_2}= \mathcal{O}(e^{-R^2\beta^{3/2}})$ defined in Proposition \ref{prop2_ZO_before}, and $C_{6,2} $ given in \eqref{constant_C_6_2_main_proof}. The second moment of $x^{\pi_\beta}$ in \eqref{lemmaT1_inequality} was estimated using Remark \ref{dissipativity_assumption_remark_ZO}, and $\int_{\mathbb{R}^d}\mathcal{G}\widetilde{V}_2(x) \pi_\beta^{f_{\eta}}(\rmd x)=0$, with $\mathcal{G}$ being the infinitesimal generator  associated with the Langevin SDE \eqref{zeroth_order_Langevin_SDE}, and $\widetilde{V}_2(x) := 1+ |x|^2$. In other words, 
        \begin{equation*} %\label{control_second_moment_invariant_measure_proof}
            \begin{split}
\E [|x^{\pi_\beta}|^2] & \leq\int_{\mathbb{R}^d}\widetilde{V}_2(x) \pi_\beta^{f_{\eta}}(\rmd x)
\\ & \leq -(2 a)^{-1} \int_{\mathbb{R}^d}\mathcal{G}\widetilde{V}_2(x) \pi_\beta^{f_{\eta}}(\rmd x) +(a+ b +d \beta^{-1})/a
    \\ & \le  1+ b a^{-1} + d (\beta a)^{-1}.
    %(a+ b +d \beta^{-1})/a,
            \end{split}
        \end{equation*}
Therefore,
     \begin{align} 
C_1^{\#}  & = \mathcal{O}(e^{R^2 \beta^{\frac{3}{2}}} (d + \beta^{-1}+ \beta^{-\frac{3+\tau}{2}} + \beta^{-\frac{5+\tau}{4}} )), \label{constants_excess_risk_bounds_C_1_non_asymptotic_notation}
\\  C_2^{\#}   & = \mathcal{O}( d^{\frac{1}{2}} (1+ d \beta^{-1} )^{\frac{3}{2}} \beta^{-5/4} e^{R^2 \beta^{\frac{3}{2}}}). \label{constants_excess_risk_bounds_C_2_non_asymptotic_notation}
        \end{align}
The estimate of the second term on the right-hand side of \eqref{splitting_expected_excess_risk_ZO} follows by applying  \cite[Proposition 11]{pmlr-v65-raginsky17a} with the restriction that $\beta \ge 2/a$, 
    \begin{equation} \label{lemmaT2_inequality} 
%\begin{split}
    \E  [f_{\eta}(x^{\pi_\beta})] - \inf\limits_{x \in \mathbb{R}^d} f_{\eta}(x) 
   \leq \frac{d}{2\beta}\log\left(\frac{ e L_1 \E [\varphi(\xi)]}{a } \left( \frac{b \beta}{d} +1 \right) \right) .
  %  \end{split}
\end{equation}
%- \frac{1}{\beta} \log \left( \frac{2 \pi^{\frac{d}{2}}}{d \  \Gamma (\frac{d}{2})} \right) + \frac{2}{\beta} + \frac{d}{\beta} \log(\beta M)
Using Lemma \ref{remainder_lemma} with the coupling \eqref{coupling} and the estimates \eqref{lemmaT1_inequality} and \eqref{lemmaT2_inequality} on the right-hand side of \eqref{splitting_expected_excess_risk_ZO}, we obtain
\begin{equation} \label{expected_excess_risk_g_epsilon_proof_semiconvexity_ZO_regularized_objective}
    \begin{split}
         \E  [v_{\eta}(x_K)] - \inf_{x \in \mathbb{R}^d} v_{\eta}(x) \le C_1^{\#} e^{-C_{r_2}\lambda K}  +  C_2^{\#} \lambda^{1/4} + C_3^{\#},
     %    \\ & \le  \Bigg( 2 L_1 \E[\varphi(X)]    (C_1(1+\mathbb{E}[|\theta_0|^2]) +   1+ b a^{-1} + d (\beta a)^{-1})  
 %      \\ & \qquad \quad    +  2 L_1 \E[\varphi(X)]   \eta \sqrt{d}  
 %          +  2  L_2 \E [\widetilde{\varphi}(X)] + 2 |\nabla F(0,0)| \Bigg) 
 %          \\ & \qquad \quad  \times (C_{W_2}e^{-C_{r_2}\lambda K} \epsilon_0^{(2)} +C_{5,2}\lambda^{1/4}) \\ &  \qquad \quad + \frac{d}{2\beta}\log\left(\frac{ e L_1 \E [\varphi(X)]}{a } \left( \frac{b \beta}{d} +1 \right) \right) +  \frac{  3 L_3 \E [\bar{\varphi}(X)]}{ \beta^{1+\tau} d^2}
    \end{split}
\end{equation}
where
\begin{equation} \label{constants_excess_risk_bounds_C_3}
    \begin{split}
 %    C_1^{\#}  & :=  \Bigg(  L_1 \E[\varphi(X)]    (C_1(1+\mathbb{E}[|\theta_0|^2]) +   1+ b a^{-1} + d (\beta a)^{-1})  
 %      \\ & \qquad \quad    +   L_1 \E[\varphi(X)]   \eta \sqrt{d}  
 %          +    L_2 \E [\widetilde{\varphi}(X)] +  |\nabla F(0,0)| \Bigg) C_{W_2} \epsilon_0^{(2)},
%\\ C_2^{\#}  & := \Bigg(  L_1 \E[\varphi(X)]    (C_1(1+\mathbb{E}[|\theta_0|^2]) +   1+ b a^{-1} + d (\beta a)^{-1})  
 %      \\ & \qquad \quad    +   L_1 \E[\varphi(X)]   \eta \sqrt{d}  
 %          +    L_2 \E [\widetilde{\varphi}(X)] +  |\nabla F(0,0)| \Bigg) C_{5,2},
  %  \\
    C_3^{\#}  & := \frac{d}{2\beta}\log\left(\frac{ e L_1 \E [\varphi(\xi)]}{a } \left( \frac{b \beta}{d} +1 \right) \right) +  \frac{  3 L_3 \E [\bar{\varphi}(\xi)]}{ \beta^{1+\tau} d^2}.
    \end{split}
\end{equation}
For $\epsilon \in (0,1]$, if we choose $\beta $ such that $C_3^{\#} \le \epsilon /3$, then choose $\lambda$ such that $\lambda < \lambda_{\text{max}}$ with $\lambda_{\text{max}}$ given in \eqref{expression_lambda_max_ZO} and $C_2^{\#}  \lambda^{1/4} \le \epsilon /3$, and finally choose $K$ such that $C_1^{\#} e^{-C_{r_2}\lambda K} \le \epsilon / 3 $, consequently, we obtain  $\E  [v_{\eta}(x_K)] - \inf_{x \in \mathbb{R}^d} v_{\eta}(x) \le \epsilon $. This yields  
\begin{equation} \label{beta_delta_expression_excess_risk_bound}
\begin{split}
  \beta & \geq \beta_{\epsilon}  := \max \Bigg \{ \beta_{\epsilon}^{\#}, \frac{9d}{2 \epsilon}\log\left(\frac{e L_1\mathbb{E}[\varphi(\xi)]}{ad}\left(b+1\right)\left(d+1\right)\right), \\ & \qquad \qquad \qquad  \ \  \left( \frac{27 L_3 \E [\bar{\varphi}(\xi)]}{ d^2 \epsilon}  \right)^{\frac{1}{1+\tau}}, \frac{2}{a} \Bigg \},   
  \end{split}
\end{equation}
 where $\beta_{\epsilon}^{\#}$ is the root of the function $g^\#(\beta) = \frac{\log\left( \beta+1\right)}{ \beta}-\frac{2 \epsilon}{9d}$, $\beta>0$, i.e. $g^\#(\beta_{\epsilon}^{\#})=0 $, and the last restriction in \eqref{beta_delta_expression_excess_risk_bound} is given by the use of \cite[Proposition 11]{pmlr-v65-raginsky17a} in the derivation of \eqref{lemmaT2_inequality}. Furthermore, since 
\begin{equation*}
C^{\#}_3 \leq \frac{d}{2\beta}\log\left(\frac{e L_1\mathbb{E}[\varphi(\xi)]}{ad}\left(b+1\right)\left(d+1\right)\left(\beta+1\right)\right) +  \frac{  3 L_3 \E [\bar{\varphi}(\xi)]}{ \beta^{1+\tau} d^2},
\end{equation*}
by setting $\frac{d}{2\beta}\log\left(\frac{e L_1\mathbb{E}[\varphi(\xi)]}{ad}\left(b+1\right)\left(d+1\right)\right) \leq \epsilon/9$, $\frac{d}{2\beta}\log\left( \beta+1\right)  \leq \epsilon/9$, together with $\frac{  3 L_3 \E [\bar{\varphi}(\xi)]}{ \beta^{1+\tau} d^2} \leq \epsilon/9  $, yields $C^{\#}_3 \leq \epsilon/3$. The desired result follows since $\frac{\log(\beta+1)}{\beta}$ is decreasing in $\beta$. Moreover, we obtain $ 0 < \lambda < \min \left \{ \frac{\epsilon^4}{81(C^{\#}_2)^4},  \lambda_{\max} \right \}$ with $\lambda_{\text{max}}$ given in \eqref{expression_lambda_max_ZO} and choose $ \lambda = \lambda_{\circ}:= \frac{1}{2} \min \left \{ \frac{\epsilon^4}{81(C^{\#}_2)^4},  \lambda_{\max} \right \}$, then
\begin{equation*}
    K =  \left \lceil \lambda_{\circ}^{-1} C_{r_2}^{-1} \ln \left( \frac{3C^{\#}_1}{\epsilon} \right) \right \rceil,
\end{equation*}
 with  the order of $C^{\#}_1$ and $C^{\#}_2$ given in \eqref{constants_excess_risk_bounds_C_1_non_asymptotic_notation} and \eqref{constants_excess_risk_bounds_C_2_non_asymptotic_notation}, respectively, and $C_{r_2}= \mathcal{O}(e^{-R^2\beta^{3/2}})$.
    \end{proof}
%\begin{remark}
% Our expected excess risks established with respect to the Hessian-trace regularized objective $v$ have a mild polynomial dependence on the dimension. In contrast, the expected excess risk established for the original objective $f$ for zeroth-order ULA established in  \cite[Corollary 4.8]{naldi2026langevin} depend on the log-Sobolev constant, which leads to exponential dependence in the dimension in the noncovex setting. As in  \cite[Corollary 4.8]{naldi2026langevin}, the inverse temperature parameter $\beta = \widetilde{O}(d \epsilon^{-1})$ which is unavoidable, see e.g. \cite[Remark 2.15]{naldi2026langevin}. This is a consequence of the classical concentration of the Gibbs measure required for bounding  $ \E  [f_{\eta}(Y_{\infty})] - \inf\limits_{\theta \in \mathbb{R}^d} f_{\eta}(\theta) $. 
%\end{remark}
%Comparison with convergence to non-regularized minima, for the expected excess risk for the original objective $f$ established in \cite{naldi2026langevin}. \textcolor{red}{Note that \cite[Corollary 4.8]{naldi2026langevin} requires $\beta = \widetilde{\Theta}(d / \bar{\epsilon})$. Our theorem and corollary give stronger distributional guarantees in $W_1$ and $W_2$ than \cite[Corollary 2.8]{naldi2026langevin}.}

\section{Conclusions} \label{sec:conclusions}

In this work, we developed a non-asymptotic analysis of stochastic zeroth-order Langevin dynamics for nonconvex optimization. By exploiting the implicit curvature-regularization effect entailed in Gaussian smoothing, we derived excess-risk bounds for the resulting curvature-regularized objective and  Wasserstein guarantees for a Gibbs measure associated with this objective, without relying on log-Sobolev or Poincar{\'e} inequalities.  Our framework also accommodates stochastic oracles with non-uniform dependence on the data variable. Future work includes extending the framework to constrained optimization and more general derivative-free Langevin schemes.

\section*{Acknowledgements}
This work was supported by the InnoCORE program of the Ministry of Science and ICT (Grant No.~1.260007.01) and the 2026 Research Fund of UNIST (Grant No.~1.260001.01).

\bibliographystyle{plain}
\bibliography{references}

\begin{thebibliography}{10}

\bibitem{ahn2023escape}
Kwangjun Ahn, Ali Jadbabaie, and Suvrit Sra.
\newblock How to escape sharp minima with random perturbations.
\newblock In {\em Proceedings of the 41st International Conference on Machine
  Learning}, 2024.

\bibitem{andriushchenko2020square}
Maksym Andriushchenko, Francesco Croce, Nicolas Flammarion, and Matthias Hein.
\newblock Square {A}ttack: {A} query-efficient black-box adversarial attack via
  random search.
\newblock In {\em European Conference on Computer Vision}, pages 484--501,
  2020.

\bibitem{arora2022understanding}
Sanjeev Arora, Zhiyuan Li, and Abhishek Panigrahi.
\newblock Understanding gradient descent on the edge of stability in deep
  learning.
\newblock In {\em International Conference on Machine Learning}, pages
  948--1024. PMLR, 2022.

\bibitem{balasubramanian2022zeroth}
Krishnakumar Balasubramanian and Saeed Ghadimi.
\newblock Zeroth-order nonconvex stochastic optimization: Handling constraints,
  high dimensionality, and saddle points.
\newblock {\em Foundations of Computational Mathematics}, 22(1):35--76, 2022.

\bibitem{bolley2005weighted}
Fran{\c{c}}ois Bolley and C{\'e}dric Villani.
\newblock Weighted {C}sisz{\'a}r-{K}ullback-{P}insker inequalities and
  applications to transportation inequalities.
\newblock In {\em Annales de la Facult{\'e} des sciences de Toulouse:
  Math{\'e}matiques}, pages 331--352, 2005.

\bibitem{bruno2026flatnessaware}
Stefano Bruno, Youngsik Hwang, Jaehyeon An, Sotirios Sabanis, and Dongyoung
  Lim.
\newblock Flatness-aware stochastic gradient {L}angevin dynamics.
\newblock In {\em Proceedings of the 43rd International Conference on Machine
  Learning}, volume 306 of {\em Proceedings of Machine Learning Research},
  pages 10017--10050. PMLR, 2026.

\bibitem{bruno2025wasserstein}
Stefano Bruno and Sotirios Sabanis.
\newblock Wasserstein convergence of score-based generative models under
  semiconvexity and discontinuous gradients.
\newblock {\em Transactions on Machine Learning Research}, 2025.

\bibitem{chen2019zo}
Xiangyi Chen, Sijia Liu, Kaidi Xu, Xingguo Li, Xue Lin, Mingyi Hong, and David
  Cox.
\newblock {ZO}-{AdaMM}: {Z}eroth-order adaptive momentum method for black-box
  optimization.
\newblock {\em Advances in Neural Information Processing Systems}, 32, 2019.

\bibitem{choromanski18a}
Krzysztof Choromanski, Mark Rowland, Vikas Sindhwani, Richard Turner, and
  Adrian Weller.
\newblock Structured evolution with compact architectures for scalable policy
  optimization.
\newblock In {\em Proceedings of the 35th International Conference on Machine
  Learning}, volume~80, pages 970--978. PMLR, 2018.

\bibitem{eberle2016reflection}
Andreas Eberle.
\newblock Reflection couplings and contraction rates for diffusions.
\newblock {\em Probability Theory and Related Fields}, 166(3):851--886, 2016.

\bibitem{flaxman_online}
Abraham~D. Flaxman, Adam~Tauman Kalai, and H.~Brendan McMahan.
\newblock Online convex optimization in the bandit setting: gradient descent
  without a gradient.
\newblock In {\em Proceedings of the Sixteenth Annual ACM-SIAM Symposium on
  Discrete Algorithms}, page 385–394. Society for Industrial and Applied
  Mathematics, 2005.

\bibitem{ghadimi2013stochastic}
Saeed Ghadimi and Guanghui Lan.
\newblock Stochastic first-and zeroth-order methods for nonconvex stochastic
  programming.
\newblock {\em SIAM Journal on Optimization}, 23(4):2341--2368, 2013.

\bibitem{hochreiter1997flat}
Sepp Hochreiter and J{\"u}rgen Schmidhuber.
\newblock Flat minima.
\newblock {\em Neural Computation}, 9(1):1--42, 01 1997.

\bibitem{hwang1980laplace}
Chii-Ruey Hwang.
\newblock Laplace's method revisited: weak convergence of probability measures.
\newblock {\em The Annals of Probability}, pages 1177--1182, 1980.

\bibitem{ji2019improved}
Kaiyi Ji, Zhe Wang, Yi~Zhou, and Yingbin Liang.
\newblock Improved zeroth-order variance reduced algorithms and analysis for
  nonconvex optimization.
\newblock In {\em International Conference on Machine Learning}, pages
  3100--3109. PMLR, 2019.

\bibitem{johnston2025the}
Tim Johnston, Iosif Lytras, Nikolaos Makras, and Sotirios Sabanis.
\newblock The performance of the unadjusted {L}angevin algorithm without
  smoothness assumptions.
\newblock {\em Transactions on Machine Learning Research}, 2025.

\bibitem{lan2020first}
Guanghui Lan.
\newblock {\em First-order and stochastic optimization methods for machine
  learning}, volume~1.
\newblock Springer, 2020.

\bibitem{li2026zeroth}
Hongbo Li and James~C Spall.
\newblock Zeroth-order {L}angevin {M}onte {C}arlo via {SPSA} under noisy
  function measurements.
\newblock {\em arXiv preprint arXiv:2608.07837}, 2026.

\bibitem{liuone}
Lewis Liu and Zhaoran Wang.
\newblock One-point gradient estimators for zeroth-order stochastic gradient
  {L}angevin dynamics, 2020.
\newblock OPT2020: 12th Annual Workshop on Optimization for Machine Learning
  Poster \#96. https://opt-ml.org/oldopt/papers/2020/paper\_96.pdf.

\bibitem{malladi2023fine}
Sadhika Malladi, Tianyu Gao, Eshaan Nichani, Alex Damian, Jason~D Lee, Danqi
  Chen, and Sanjeev Arora.
\newblock Fine-tuning language models with just forward passes.
\newblock {\em Advances in Neural Information Processing Systems},
  36:53038--53075, 2023.

\bibitem{nakakita2023non}
Shogo Nakakita.
\newblock Non-asymptotic analysis of {L}angevin-type {M}onte {C}arlo
  algorithms.
\newblock {\em Advances in Applied Probability}, page 1–38, 2025.

\bibitem{naldi2026langevin}
Emanuele Naldi, Marco Rando, Lorenzo Rosasco, and Silvia Villa.
\newblock Langevin for nonconvex optimization: Exact, inexact and zeroth-order.
\newblock {\em arXiv preprint arXiv:2607.22353}, 2026.

\bibitem{nesterov2017random}
Yurii Nesterov and Vladimir Spokoiny.
\newblock Random gradient-free minimization of convex functions.
\newblock {\em Foundations of Computational Mathematics}, 17:527--566, 2017.

\bibitem{pmlr-v65-raginsky17a}
Maxim Raginsky, Alexander Rakhlin, and Matus Telgarsky.
\newblock Non-convex learning via stochastic gradient {L}angevin dynamics: {A}
  nonasymptotic analysis.
\newblock In {\em Proceedings of the 2017 Conference on Learning Theory},
  volume~65, pages 1674--1703. PMLR, 2017.

\bibitem{roy2022stochastic}
Abhishek Roy, Krishnakumar Balasubramanian, Saeed Ghadimi, and Prasant
  Mohapatra.
\newblock Stochastic zeroth-order optimization under nonstationarity and
  nonconvexity.
\newblock {\em Journal of Machine Learning Research}, 23(64):1--47, 2022.

\bibitem{roy2022stochasticLangevin}
Abhishek Roy, Lingqing Shen, Krishnakumar Balasubramanian, and Saeed Ghadimi.
\newblock Stochastic zeroth-order discretizations of {L}angevin diffusions for
  {B}ayesian inference.
\newblock {\em Bernoulli}, 28(3):1810--1834, 2022.

\bibitem{sahin2026zerothorder}
M.~Berk Sahin, Behzad Sharif, and Abolfazl Hashemi.
\newblock Zeroth-order non-log-concave sampling with variance reduction and
  applications to inverse problems.
\newblock In {\em Forty-third International Conference on Machine Learning},
  2026.

\bibitem{salimans2017evolution}
Tim Salimans, Jonathan Ho, Xi~Chen, Szymon Sidor, and Ilya Sutskever.
\newblock Evolution strategies as a scalable alternative to reinforcement
  learning.
\newblock {\em arXiv preprint arXiv:1703.03864}, 2017.

\bibitem{Shamir13}
Ohad Shamir.
\newblock On the complexity of bandit and derivative-free stochastic convex
  optimization.
\newblock In {\em Proceedings of the 26th Annual Conference on Learning
  Theory}, volume~30. PMLR, 2013.

\bibitem{silveri2025beyond}
Marta~Gentiloni Silveri and Antonio Ocello.
\newblock Beyond log-concavity and score regularity: {I}mproved convergence
  bounds for score-based generative models in {$W_2$}-distance.
\newblock In {\em Forty-second International Conference on Machine Learning},
  2025.

\bibitem{snoek2012practical}
Jasper Snoek, Hugo Larochelle, and Ryan Adams.
\newblock Practical {B}ayesian optimization of machine learning algorithms.
\newblock {\em Advances in neural information processing systems}, 25, 2012.

\bibitem{welling2011bayesian}
Max Welling and Yee~W Teh.
\newblock Bayesian learning via stochastic gradient {L}angevin dynamics.
\newblock In {\em Proceedings of the 28th International Conference on Machine
  Learning}, pages 681--688, 2011.

\bibitem{wen2023how}
Kaiyue Wen, Tengyu Ma, and Zhiyuan Li.
\newblock How sharpness-aware minimization minimizes sharpness?
\newblock In {\em The Eleventh International Conference on Learning
  Representations}, 2023.

\bibitem{zhang2024dpzero}
Liang Zhang, Bingcong Li, Kiran~Koshy Thekumparampil, Sewoong Oh, and Niao He.
\newblock {DPZero}: {P}rivate fine-tuning of language models without
  backpropagation.
\newblock In {\em International Conference on Machine Learning}. PMLR, 2024.

\bibitem{zhang2025zeroth}
Liang Zhang, Bingcong Li, Kiran~Koshy Thekumparampil, Sewoong Oh, Michael
  Muehlebach, and Niao He.
\newblock Zeroth-order optimization finds flat minima.
\newblock In {\em The Thirty-ninth Annual Conference on Neural Information
  Processing Systems}, 2025.

\bibitem{zhang2023nonasymptotic}
Ying Zhang, {\"O}mer~Deniz Akyildiz, Theodoros Damoulas, and Sotirios Sabanis.
\newblock Nonasymptotic estimates for stochastic gradient {L}angevin dynamics
  under local conditions in nonconvex optimization.
\newblock {\em Applied Mathematics \& Optimization}, 87(2):25, 2023.

\end{thebibliography}
\end{document}